\documentclass{memol}

\usepackage{mathdots, bbm,mathabx,stmaryrd}
\usepackage[pdftex]{color,graphicx}
\usepackage[all,cmtip]{xy}
\makeindex 
\newcommand\Tr{\mathrm{Tr}}
\newcommand\esp{\mathbb E}
\newcommand\etc{,\ldots ,}
\newcommand\one{\mathbbm{1}}
\newcommand\toN{^{(N)}}
\newcommand\limN{\underset{N \rightarrow \infty}\longrightarrow}
\newcommand\Nlim{\underset{N \rightarrow \infty}\lim}
\newcommand\lara{\langle \mathbf x, \mathbf x^*\rangle}
\newcommand\eq{\begin{eqnarray*}}
\newcommand\qe{\end{eqnarray*}}
\newcommand\eqa{\begin{eqnarray}}
\newcommand\qea{\end{eqnarray}}
\newcommand{\circlearrowleftRot}{\mathbin{\text{\rotatebox[origin=c]{90}{$\circlearrowleft$}}}}
\newcommand{\circlearrowleftRotBis}{\mathbin{\text{\rotatebox[origin=c]{-90}{$\circlearrowright$}}}}
\newtheorem{theorem}{Theorem}[chapter]
\newtheorem{proposition}[theorem]{Proposition}
\newtheorem{corollary}[theorem]{Corollary}
\newtheorem{lemma}[theorem]{Lemma}
\theoremstyle{definition}
\newtheorem{definition}[theorem]{Definition}
\newtheorem{example}[theorem]{Example}
\theoremstyle{remark}
\newtheorem{remark}[theorem]{Remark}
\numberwithin{section}{chapter}
\numberwithin{equation}{chapter}

\begin{document}

\frontmatter
\title[Traffic Distributions and Independence]{Traffic Distributions and Independence: \\Permutation Invariant Random Matrices and the Three Notions of Independence}
\author{Camille MALE}
\address{Institut de Math\'ematiques de Bordeaux UMR5251 CNRS}
\curraddr{}
\email{camille.male@math.u-bordeaux.fr}
\urladdr{https://camillemale.com}
\thanks{Partially supported by the ANR project ANR-08-BLAN-0311-01, the Fondation Sciences Math\'ematiques de Paris, and the MAP5 at the university Paris Descartes.}
\date{on Nov 2016 at {\it Memoirs of the AMS}. Accepted on Oct 2017}
\maketitle

\tableofcontents
\begin{abstract}
	Voiculescu's notion of asymptotic free independence is known for a large class of random matrices including independent unitary invariant matrices. This notion is extended for independent random matrices invariant in law by conjugation by permutation matrices. This fact leads naturally to an extension of free probability, formalized under the notions of \emph{traffic probability}. 

We first establish this construction for random matrices. We define the traffic distribution of random matrices, which is richer than the $^*$-distribution of free probability. The knowledge of the individual traffic distributions of independent permutation invariant families of matrices is sufficient to compute the limiting distribution of the join family. Under a factorization assumption, we call traffic independence the asymptotic rule that plays the role of independence with respect to traffic distributions. Wigner matrices, Haar unitary matrices and uniform permutation matrices converge in traffic distributions, a fact which yields new results on the limiting $^*$-distributions of several matrices we can construct from them. 

Then we define the abstract traffic spaces as non commutative probability spaces with more structure. We prove that at an algebraic level, traffic independence in some sense unifies the three canonical notions of tensor, free and Boolean independence. A central limiting theorem is stated in this context, interpolating between the tensor, free and Boolean central limit theorems.
\end{abstract}
\mainmatter

\chapter*{Introduction}
\specialsection*{Presentation of the Article}

\par Motivated by the study of von Neumann algebras of free groups, Voiculescu introduces  \cite{VOI6} free probability theory as a non commutative probability theory equipped with the well known notion of free independence. 
Voiculescu shows \cite{VOI5}  that free independence describes the global asymptotic behavior of eigenvalues of independent unitarily invariant matrices and Wigner matrices. Asymptotic free independence holds for a large class of matrices, see for instance \cite{BeG09,CC,HP,ScSc05,Shl98}. 
\par Recall that in free probability theory, a non commutative probability space $(\mathcal A, \Phi)$ consists in is a unital (non commutative) algebra $\mathcal A$ equipped with a linear form $\Phi$ playing the role of the expectation and called the sate (satisfying mild assumptions, see Definition \ref{Def:NCPS}). This context is considered together with the notion of free independence which replaces the classical notion of tensor independence. It is a canonical rule that, given two probability spaces $(\mathcal A_1,\Phi_1)$ and $(\mathcal A_2, \Phi_2)$, associates a state $\Phi=\Phi_1*\Phi_2$ on the free product algebra $\mathcal A_1*\mathcal A_2$ whose restriction on $\mathcal A_i$ is $\Phi_i$ for $i=1,2$. Speicher proves in \cite{SPI97} that tensor independence and free independence are the only universal notions of independence (symmetric and associative, and satisfying a universal calculation rule), and that in the non unital case the third notion of Boolean independence appears.
\\
\par  In the first part of this article we introduce a general method to study random matrices which are not invariant by conjugation by unitary matrices: we start by considering independent families $\mathbf A_N^{(1)} \etc \mathbf A_N^{(L)}$ of random matrices which are invariant by conjugation by permutation matrices (Definition \ref{Def:Intro}). We consider the problem of characterizing the limiting joint $^*$-distribution of the collection of all matrices $ \mathbf A_N^{(1)} \cup \dots \cup \mathbf A_N^{(L)}$ when the size of the matrices goes to infinity (the $^*$-distribution is the expectation of normalized trace of polynomials in the matrices, see Section \ref{Sec:AsymStarFrennes}). If the families are not unitarily invariant, in general the knowledge of the limiting $^*$-distribution of each $\mathbf A_N^{(\ell)}$ is not sufficient to characterize the limiting $^*$-distribution of $\mathbf A_N^{(1)} \cup \dots \cup \mathbf A_N^{(L)}$. 
\par This question is studied in \cite{RYA} when the matrices are independent and have independent and identically distributed entries. In particular if the matrices have size $N\times N$ and their entries are Bernouilli random variables with parameter $\frac p N$, $p$ fixed, then the matrices are not asymptotically free.  In another direction, an interpolation of the classical and free convolutions is introduced and studied in \cite{BGL11}: given two distributions $\mu$ and $\nu$ on the real line, the distribution $ \mu \boxplus_t \nu, t\geq 0$ is the limit of the empirical spectral distribution of matrices $A_N + U_N^{(t)} B_N U_N^{(t)}$ where $U_N^{(t)}$ is a diffusion on the unitary group starting from a uniform permutation matrix, and $A_N, B_N$ are diagonal matrices that converge in moments to $\mu$ and $\nu$ respectively. 
\par The strategy therein is to consider more observables on matrices, and then to define a more general setting that free probability theory where we need to enrich the notion of $^*$-distribution. For a family $\mathbf A_N$ of matrices, we define the \emph{traffic distribution} as a linear functional on finite connected directed graphs whose edges are labelled by matrices (Definition \ref{Def:StarGraphsMon}). It contains more information than the data of the normalized trace of monomials of $\mathbf A_N$. Hence the convergence in traffic distribution of a family of random matrices implies its convergence in $^*$-distribution, but it is possible that two random matrices have the same $^*$-distribution but not the same traffic distribution. Equivalently, the convergence in traffic distribution of $\mathbf A_N$ is the convergence in $^*$-distribution of the matrices of the so-called \emph{traffic space} generated by $\mathbf A_N$ (Lemma \ref{Lem:EquivCVTraf}).

\par  In the main result Theorem \ref{MainTh}, we prove an analogue of Voiculescu's theorem, by giving a formula for the limiting traffic distribution of independent permutation invariant families of matrices $\mathbf A_N^{(1)} \cup \dots \cup \mathbf A_N^{(L)}$,
knowing the limiting traffic distribution of each one. Our method extends the moment method in a similar fashion as in the first chapter of \cite{GUI} and applies for a large class of random matrices.  In Theorem \ref{MainTh}, the limiting traffic distribution of $\mathbf A_N^{(1)} \cup \dots \cup \mathbf A_N^{(L)}$ is called the product of the limiting traffic distributions of the $\mathbf A_N^{(\ell)}$'s. We will also say that the random matrices $ \mathbf A_N^{(1)} \cup \dots \cup \mathbf A_N^{(L)}$ are asymptotically traffic independent.
\par With Theorem \ref{MainTh} and several propositions stated in the article, we illustrate new ways to test the asymptotically freeness of independent random matrices. We consider Wigner, unitary Haar and uniform permutation matrices, as well as matrices obtained from several operations on them (transpose, entry-wise product). We also point out new examples of matrices that are not asymptotically free independent. Other models are considered in \cite{MP14,DLM15,MAL122}. 
\\
\par In the second part of the article, the notion of non commutative probability space is enhanced to consider the limits in traffic distribution of random matrices. The notion of operads (see e.g. \cite{May97}) formalizes the idea that a given set of operations with good compatibility conditions defines an algebraic structure. We define the operad of \emph{graph operations}, which is a slight modification of Jones' planar algebras \cite[Example 2.6]{Jon99} and Spivak's wiring diagram algebras \cite{Spi13}. 
An algebraic traffic space is an algebra over the operad of graph operations endowed with a linear form playing the role of the expectation (satisfying mild assumptions). An element of a traffic space are simply called a \emph{traffic}. In particular, a traffic space is a non commutative probability space. Roughly speaking, traffics have the property that we can compose them not only by linear operations, but thanks to schemes given by finite connected graphs with an input and an output. 

The interest in defining traffic spaces is that we can consider traffic independence in a universal way and state the usual limit theorems of probability. Traffic independence is a notion more general that the notions of independence in non commutative probability. As the $^*$-distribution is just a part of the traffic distribution, traffic independence encodes a large class of relations for non commutative random variables. In particular, it encodes both the tensor and the free independence of $^*$-distributions: we exhibit two different classes of traffics for which the traffic independence is equivalent to the tensor or free independence. Moreover, the additional structure of observables of traffic spaces implies the existence of another linear form than the usual trace. For a third class of traffics, with respect to this linear map, traffic independence encodes the notion of Boolean independence.

This connection with Boolean independence is actually deduced from the limit theorems. Let $(a_n)_{n\geq 1}$ be a family of independent identically distributed self-adjoint traffics. We prove a law of large number for the sum $\frac{a_1+ \dots + a_n}{n}$, and assuming the variables are centered in a certain sense, we prove the convergence in traffic distribution of $\frac{a_1+ \dots + a_n}{\sqrt{n}}$ and interpret the limit. We prove that $\frac{a_1+ \dots + a_n}{\sqrt{n}}$ converges in $^*$-distribution to the sum $x +z$ of a semicircular variable $x$ free independent from a Gaussian variable $z$ (the two limits in the central limit theorems for free and tensor independence). The degree of freedom between the classical and the free worlds is possible thanks to the wealthy of information contained in the traffic distribution.

Nonetheless, the limiting traffic distribution of $\frac{a_1+ \dots + a_n}{\sqrt{n}}$ is actually the distribution of the sum of three variables $x + y+ z$. Seen as a non commutative random variable, the variable $y$ has variance zero, which explains why it does not appear in the convergence in $^*$-distribution. But it has a nonzero traffic distribution, and can actually be interpreted as the limit in the central limit theorem for Boolean independence. Remarkably, the variables $x$, $y$ and $z$ in the central limit theorem are not traffic independent in general. We give a matrix model $(X_N, Y_N, X_N)$ which converges in traffic distribution to $(x,y,z)$.
\\
\par In parallel to the evolution of this paper, Gabriel gives in \cite{GAB15,GAB152, GAB153} another answer to Theorem \ref{MainTh} motivated by the Schur-Weyl duality. His approach involves the \emph{partition algebras} instead of the test graphs defined in this article. Gabriel's theory is essentially equivalent to the one introduced in this article. A dictionary \cite{CDGM} between the two approaches is in preparation.
\\
\par The organization of the article is the following. In Chapter \ref{Theme1}, we define the traffic distributions of large random matrices and present our main result (Theorem \ref{MainTh}) and its applications. In Chapter \ref{Sec:DefIndepTraffic} we  define traffic independence and prove a simple criterion to prove that random matrices are not asymptotically free independent. We also prove the main result Theorem \ref{MainTh}. Chapter \ref{Sec:ApplWignerMatrices} is dedicated to applications of this theorem for Wigner matrices,  uniform permutation matrices and unitary Haar matrices. 
In Chapter \ref{Sec:TrafficSpaces} we introduce the general traffic spaces. In Chapter \ref{Sec:ThreeIndep}, we relate traffic independence with the three universal notions of independence.  The law of large numbers and the central limit theorem for traffic independence are stated in Chapter 6.

\specialsection*{Notations and Preliminaries}

While considering a matrix $A_N$, we implicitly mean a sequence $(A_N)_{N\geq 1}$, the matrix $A_N$ being of size $N$. We study large matrices and the term ''asymptotic`` refers to the limit when $N$ goes to infinity. We consider random matrices $A_N$ whose entries admit moments of all orders, that is $\forall K\geq 1, \forall n,m=1\etc N$, $\esp\big[|A_N(n,m)|^{K}\big]<\infty$. We denote by $\mathbb I_N$ the identity matrix and for $A_N$ a complex matrix we denote by $A_N^*$ its conjugate transpose. We often consider families of matrices, denoted by bold characters e.g. $\mathbf A_N = (A_j)_{j\in J}$. For a family $\mathbf A_N$ of matrices, we denote by $\mathbf A_N^*$ the family of their complex transpose.

\begin{definition}[Unitary random matrices and invariances]\label{Def:Intro} 
~

\begin{enumerate}
	\item A unitary matrix $U_N$ is a matrix such that $U_NU_N^* =U_N^*U_N=\mathbb I_N$. {\bf A unitary Haar matrix} $U_N$ is a random unitary matrix distributed according to the Haar distribution on the unitary group, that is the unique probability measure on the unitary group invariant by right and left multiplication of elements of the group.
	\item A permutation matrix $V_N$ of size $N$ is a unitary matrix for which there is a permutation $\sigma$ of $\{1\etc N\}$ such that the entry $(i,j)$ of $V_N$ is one if $i=\sigma(j)$ and zero otherwise. {\bf A uniform permutation matrix} is a random permutation matrix uniformly chosen among all the $N!$ choices.
	\item  A family of random matrices $\mathbf A_N = (A_j)_{j\in J}$ is said to be {\bf unitarily invariant} whenever it is invariant in law by conjugation by unitary matrix, that is for any unitary matrix $U_N$
		$$\mathbf A_N \overset{ \mathcal Law} =U_N \mathbf A_NU_N^*:=(U_N A_{j}U_N^*)_{j \in J}.$$
 Equivalently, $\mathbf A_N$ has the same law as $U_N \mathbf A_N U_N^*$, where $U_N$ is a unitary Haar matrix independent of $\mathbf A_N$.
	\item  A family of random matrices is said to be {\bf permutation invariant} whenever it is invariant in law by conjugation by any permutation matrix. Equivalently, $\mathbf A_N$ has the same law as $V_N \mathbf A_N V_N^*$ where $V_N$ is a uniform permutation matrix independent of $\mathbf A_N$.
\end{enumerate}
\end{definition}

\begin{definition}[Wigner matrices]\label{Def:WigMatrices} A (centered) complex Wigner matrix is a Hermitian matrix $X_N =\big( \frac{x_{i,j}}{\sqrt n}\big)$ whose sub-diagonal entries are independent and centered random variables such that: 
\begin{enumerate}
	\item the diagonal entries $(x_{i,i})_{i=1\etc N}$, respectively the sub-diagonal entries $(x_{i,j})_{j<i}$, are identically distributed,
	\item the distribution of $x_{i,j}$ does not depend on $N$, has finite moments of all orders ($\esp[|x_{i,j}|^k]<\infty$ for any $k\geq 1$) and $\esp[x_{i,j}]= 0$ for any $i,j=1\etc N$.
\end{enumerate}
We call parameter of $X_N$ the common value $(\alpha,\beta)$ of $(\esp[|x_{i,j}|^2], \esp[x_{i,j}^2])$ for $i\neq j$. A real Wigner matrix is a complex matrix with real entries.
\end{definition}
Note that a complex Wigner matrix is almost surely a real Wigner matrix if and only if $\alpha=\beta$.

\begin{lemma} A unitary Haar matrix, a uniform permutation matrix and a real Wigner matrix are permutation invariant. A complex Wigner matrix is permutation invariant if and only if the entries have the same distribution as their complex conjugate (i.e. $x_{i,j} \overset{\mathrm{\mathcal Law}}=\overline{x_{i,j}}$).
\end{lemma}
\begin{proof} Let $U_N$ be a unitary Haar matrix. For any permutation matrix $V_N$, the matrix $V_NU_NV_N^*$ has the same distribution as $U_N$ since $V_N$ is unitary and $U_N$ is Haar distributed.
\par Let $W_N = \big((\one ( i=\sigma_W(j) ) \big)_{i,j}$ be a uniform permutation matrix, associated to a uniform permutation $\sigma_W$ of $\{1\etc N\}$. Then for any permutation matrix $V_N$ associated to a permutation $\sigma$, one has
	\eq
		V_NW_NV_N^* & = & \bigg( \sum_{k,\ell} \one\Big( i= \sigma(k)  , k= \sigma_W(\ell) ,  j=\sigma(\ell)  \Big) \bigg)_{i,j}\\
		& = &  \bigg(  \one\Big(  \sigma^{-1}(i)= \sigma_W( \sigma^{-1}(j)   \Big) \bigg)_{i,j} \\
		 & = & \Big( \one\big( {\sigma} \circ \sigma_W \circ  {\sigma^{-1}}(i)=j\big)\Big)_{i,j} \overset{\mathcal Law} = \big( \one\big(  \sigma_W  (i)=j\big)\big)_{i,j},
	\qe
so $V_NW_NV_N^*$ has the same distribution as $W_N$.
\par Let $X_N$ be a Wigner matrix. Denote by $E_{i,j}=\big(\one(k=i,\ell=j)\big)_{k,\ell=1\etc N}$,  the elementary matrix for each $i,j=1\etc N$. For any permutation matrix $V_N$ associated to a permutation $\sigma$, the Hermitian matrix $V_NX_NV_N^*$ has independent sub-diagonal entries since
	\eq
		\lefteqn{V_NX_NV_N^*  }\\
		& = & \sum_{i,j}\Big(  X_N(\sigma^{(-1)}(i),\sigma^{(-1)}(j))  \Big) E_{i,j}\\
		& = &  \sum_{i} E_{\sigma (i),\sigma (i)} X_N(i,i) +  \sum_{i<j}\Big(E_{\sigma (i),\sigma (j)} X_N(i,j) + E_{\sigma (j),\sigma (i)}\overline{X_N(j,i)}\Big).
	\qe
 The distribution of its diagonal entries is the distribution of those of $X_N$. The non diagonal entry $(k,\ell)$ for $k<\ell$ is distributed as $X_{\sigma(k),\sigma(\ell)}$ if $\sigma(k)<\sigma(\ell)$ and as $\overline{X_{\sigma(k),\sigma(\ell)}}$ if $\sigma(k)>\sigma(\ell)$. The invariance in distribution by complex conjugation of the non diagonal entries of a Wigner matrix is then a necessary and sufficient condition for the permutation invariance of $X_N$.\end{proof}

\part{The Asymptotic Traffic Distributions of Random Matrices}
\chapter{Statement of the Main Theorem and Applications} \label{Theme1}


\section{Asymptotic free independence of large matrices}\label{Sec:AsymStarFrennes}

\par The random matrices under consideration are assumed to have entries with finite moments of all orders, that is $\esp[|A_N(i,j)|^K]<\infty$ for any $K\geq 1$. The (mean) empirical spectral distribution of a random matrix $H_N$ is the probability measure
		$$  {\mathcal L}_{H_N} : f \mapsto \esp \Big[ \frac 1 N \sum_{i=1}^N f(\lambda_i) \Big],$$
where $\lambda_1 \etc \lambda _N$ are the eigenvalues of $H_N$ and $f : \mathbb C \to \mathbb C$ are integrable functions in $\lambda$. 
\par In random matrix and free probability, a common notion is the $^*$-distribution, which extends the notion of empirical spectral distribution for several matrices. Denote by $\mathbb C\lara$ the space of $^*$-polynomials, i.e. finite complex linear combinations of words in symbols $\mathbf x=(x_j)_{j\in J}$ and $\mathbf x^*=(x_j^*)_{j\in J}$. The $^*$-distribution of a family of random matrices $\mathbf A_N$ is defined by
	\eq
		\Phi_{\mathbf A_N} : P \in \mathbb C \lara \mapsto \esp \Big[ \frac 1 N \Tr \big( P(\mathbf A_N) \big) \Big] \in \mathbb C,
	\qe
where $\Tr$ is the usual trace of matrices. The family $\mathbf A_N$ converges in $^*$-distribution whenever $\Phi_{\mathbf A_N}$ converges pointwise as $N$ goes to infinity. The convergence in $^*$-distribution of $\mathbf A_N$ is actually equivalent to the convergence in moments of the empirical spectral distribution of any Hermitian random matrix $H_N=P(\mathbf A_N)$. 
\par Voiculescu's asymptotic freeness theorem and its extensions \cite{AGZ, COL,DYK,VOI5,VOI4} give a characterization of the limiting $^*$-distribution of a collection of unitarily invariant independent families of matrices.

 \begin{theorem}[Asymptotic free independence]~\label{Th:AsyFree}
\\Let $\mathbf A_N^{(1)} \etc \mathbf A_N^{(L)}$ be independent families of $N \times N$ random matrices. Make the following hypotheses:
\begin{enumerate}
	\item Each family, except possibly one, is unitarily invariant (Definition \ref{Def:Intro}).
	\item Each family converges in $^*$-distribution, namely 
\eq
	\Phi_\ell  ( P  ) := \Nlim \esp\Big[ \frac 1 N \Tr \big[ P(\mathbf A_N^{(\ell)}) \big] \Big]
\qe
 exists for each $\ell\in \{1\etc L\}$ and any $^*$-polynomial $P$.
	\item For each $\ell\in \{1\etc L\}$, either the matrices of $\mathbf A_N^{(\ell)}$ are uniformly bounded in operator norm and $ \frac 1 N \Tr \big[ P(\mathbf A_N^{(\ell)}) \big]  \limN \Phi_\ell  ( P  )$ almost surely for any $\ell$ and $P$, or $\mathbf A_N^{(\ell)}$ is of the form $U_N \tilde {\mathbf A}_N^{(\ell)} U_N^*$ where $U_N$ is a Haar unitary random matrix and $\tilde {\mathbf A}_N^{(\ell)}$ is deterministic.
\end{enumerate}
Then the families $\mathbf A_N^{(1)} \etc \mathbf A_N^{(L)}$ are asymptotically free independent, namely:
	\begin{enumerate}
		\item They have a limiting joint $^*$-distribution: denoting $\mathbf A_N = \mathbf A_N^{(1)} \cup \dots \cup \mathbf A_N^{(L)}$, for any $^*$-polynomial $P$ the limit 
			\eqa\label{Eq:CVDist}
				\Phi  ( P  ) := \Nlim \esp\Big[ \frac 1 N \Tr \big[ P(\mathbf A_N) \big] \Big] \textrm{ exists}.
			\qea
		\item They limiting joint $^*$-distribution $\Phi$ is the free product of the limiting distributions $\Phi_\ell$: for any $n\geq 1$ and any indices $\ell_1,\ell_2, \dots, \ell_n$ in $\{1\etc L\}$ such that  $\ell_j\neq \ell_{j+1}, \forall j=1\etc n-1$, for any $^*$-polynomials $P_1, P_2, \dots$ where $P_j \in \mathbb C\langle \mathbf x_{\ell_j},\mathbf x_{\ell_j}^*\rangle$ satisfies $\Phi_{\ell_j}( P_j) = 0, \forall j\geq 1$, one has
		\eqa\label{Eq:Free}
			\Phi ( P_1  \times \dots  \times P_n) = 0.
		\qea
	\end{enumerate}
	Moreover, the conclusion of the theorem remains valid if a family consists in independent Wigner matrices.

\end{theorem}
  
Nica \cite{NIC93} and Neagu \cite{NEA} proved that independent permutation matrices uniformly distributed are asymptotically free independent and asymptotically free from independent Wigner matrices. See \cite{AGZ,NS} for more details.
   
\section{Convergence in traffic distribution}\label{SecCVTraffics}

Our motivation is to state an analogue of Theorem \ref{Th:AsyFree} where the independent families of matrices $\mathbf A_1 \etc \mathbf A_L$ are only assumed to be permutation invariant rather than unitarily invariant (Definition \ref{Def:Intro}). For that task, we need more than the $^*$-distributions of the families to compute their limiting joint $^*$-distributions. 

\subsection{Definition}

We introduce a generalization of $^*$-polynomials in matrices. These operations are given by combinatorial graphs for which we now fix the notations. The considered graphs are directed and can have multiple edges in any directions and loops. Formally, a graph is a couple $(V,E)$, where $V$ is a non empty set, called the set of vertices, and $E$ is a multi-set (elements appear with a certain multiplicity) of ordered pairs of vertices $e=(v,w)\in V^2$. The set $E$ is possibly empty and is called the set of edges. We point out that graphs with no edges are allowed, but a graph must always have at least one vertex.
 A graph $(V,E)$ is finite if both $V$ and $E$ are finite. 
Graphs are considered up to isomorphisms. Two graphs $G_1=(V_1,E_1)$ and $G_2=(V_2,E_2)$ are identified whenever there exists a bijection $\phi : V_1 \to V_2$ preserving the adjacency of vertices, the orientation of edges and their multiplicity. 
 
\begin{figure}
      \includegraphics{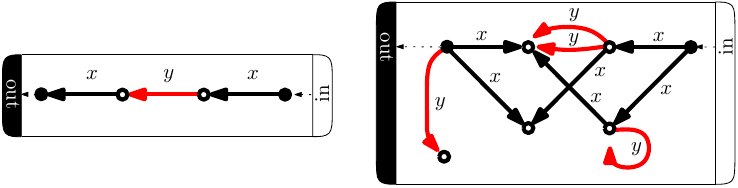}
     \caption{Left: the linear composition of operators. Right: the composition according to a graph monomial.}
    \label{fig:PolyGraphPoly.}
  \end{figure}

\begin{definition}[Graph polynomials and traffic distributions] \label{Def:StarGraphsMon} Let $J$ be an index set and consider two families of formal variables $\mathbf x=(x_j)_{j\in J}$ and $\mathbf x^*=(x^*_j)_{j\in J}$, i.e. pairwise distinct symbols. \begin{enumerate}
	\item A $^*$-test graph in the variables $\mathbf x$ is a finite, connected, oriented graph whose edges are labelled by $x_j$ and $x_j^*$, $j\in J$. Formally, it consists in a quadruple $T=(V,E,\gamma,\varepsilon)$, where 
\begin{enumerate}
	\item $(V,E)$ is a finite connected graph,
	\item $\gamma$ is a map $E \to J$,
	\item $\varepsilon$ is a map $E \to \{1 ,*\}$.
\end{enumerate}
The maps $\gamma$ and $\varepsilon$ indicate that an edge $e \in E$ has the label $x_{\gamma(e)}^{\varepsilon(e)}$. For multiple edges, each edge has its own label. We simply call $T$ a test graph when $\varepsilon$ is constant to one and denote it $T=(V,E,\gamma)$. Similarly, a $^*$-test graph in one variable is denoted $T=(V,E,\varepsilon)$.

\item A $^*$-graph monomial is a $^*$-test graph with two given vertices: an input and an output. It consists in a triplet $g=(T,in,out)$ where $T$ is as above and $in,out\in V$ are possibly equal. We denote by $\mathbb C \mathcal G\lara$ the space of finite linear complex combinations of $^*$-graph monomials in variable $\mathbf x$ (graphs are considered up to isomorphisms of graphs preserving labeled and in/outputs). 
\item Denote $[N]:=\{1 \etc N\}$. For any $^*$-graph monomial $g=(T,in, out)$, $T=(V,E,\gamma,\varepsilon)$, in the variables $\mathbf x=(x_j)_{j\in J}$, and for any family $\mathbf A_N=(A_j)_{j\in J}$ of matrices, we define the matrix $g(\mathbf A_N)$ whose $(i,j)$-entry is 
\eqa\label{DefIntroGraphMonMatrix}
	g(\mathbf A_N)(i,j) = \sum_{\substack{ \phi:V \to [N] \\ \phi(out)=i, \phi(in)=j}} \prod_{e=(v,w)\in E} A_{\gamma(e)}^{\varepsilon(e)}\big( \phi(w), \phi(v)\big).
	\qea
	This definition is extended for $g\in \mathbb C \mathcal G \lara$ by linearity. Note that for the graph monomial consisting in a simple line $g=(\underset{out}\cdot \overset{x_{\gamma_1}^{\varepsilon_1}}\leftarrow \cdot \, \dots \, \cdot \overset{x_{\gamma_K}^{\varepsilon_K}}\leftarrow  \underset{in} \cdot)$, then $g(\mathbf A_N)$ is the product of matrices $A_{\gamma_1}^{\varepsilon_1} \dots A_{\gamma_K}^{\varepsilon_K}$, see Figure \ref{fig:PolyGraphPoly.}.

	\item The traffic space generated by $\mathbf A_N$ is the sub-space of $\mathrm M_N(\mathbb C)$ consisting in all matrices $g(\mathbf A_N)$ for $g \in \mathbb C \mathcal G\lara$. The traffic distribution of $\mathbf A_N$ is the map 
		$$\bar\Phi_{\mathbf A_N} : g  \in \mathbb C \mathcal G\lara \mapsto \esp\Big[ \frac 1 N \Tr \big[ g(\mathbf A_N) \big] \Big] \in \mathbb C.$$
	We say that a family $\mathbf A_N$ of random matrices converges in traffic distribution (implicitly when $N\rightarrow \infty$) whenever  $\Nlim \bar\Phi_{\mathbf A_N}(\,g\,)$ exists for any $^*$-graph polynomial $g$.	
\end{enumerate}
\end{definition}

The term \emph{traffics} refers to the matrices $\mathbf A_N$ and is formalized abstractly in Section \ref{Sec:TrafficSpaces}. Definition \eqref{DefIntroGraphMonMatrix} is considered in \cite{MS09}. 
The traffic distribution of family $\mathbf A_N$ is denoted $\bar \Phi_{\mathbf A_N}$, to not be confused with the $^*$-distribution $\Phi_{\mathbf A_N}$. 

\begin{figure}
     \includegraphics{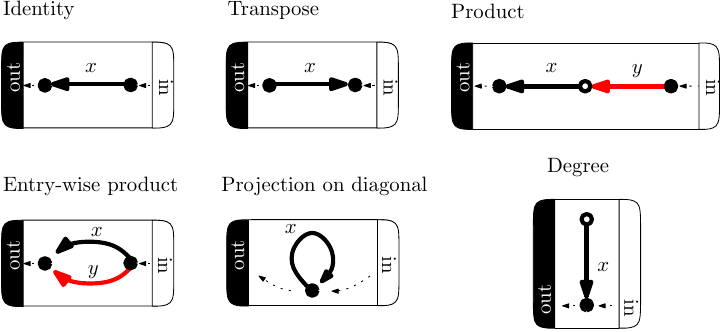}
    \caption{Example of graph monomials.}
    \label{fig:ExGraphPoly}
  \end{figure}

\begin{example}\label{Ex:GraphPolyMatrices}Consider $\mathbf A_N=(A_j)_{j\in J}$ a family of matrices and let us describe $g(\mathbf A_N)$ for the following $^*$-graph monomials $g$ in variables $\mathbf x = (x_j)_{j\in J}$ (see Figure \ref{fig:ExGraphPoly}).
 \begin{enumerate}

	\item {\bf Identity:} Let $g$ be the graph monomial with two vertices $in$ and $out$ and one edge from $in$ to $out$ labelled $x_j$. Then $g(\mathbf A_N) = A_j$. We denote in short $g=( \underset{out}\cdot \overset{x_j}{\leftarrow}  \underset{in}\cdot)=( \cdot \overset{x_j}{\leftarrow}  \cdot)$, with the convention that the vertex $out$ is on the left and the vertex $in$ is on the right. This is the same convention as for a matrix entry $A_{k,\ell}$, where the indices $\ell$ and $k$ indicate the element of the basis for the source and target spaces respectively of the associated linear map $\mathbb C^N \to \mathbb C^N$. 
	\item {\bf Constant:} For $(\cdot )$ the graph with a single vertex $in=out$ and no edge, one has that $(\cdot)(\mathbf A_N)=\mathbb I_N$ the identity matrix. 
	\item {\bf Product:} For $(\cdot \overset{x_{j'}}{\leftarrow} \cdot \overset{x_{j}}{\leftarrow} \cdot)$ the graph with three vertices $out, v, in$, one edge from $in$ to $v$ labelled $x_j$ and another one from $v$ to $out$ labelled $x_{j'}$, one has $(\cdot \overset{x_{j'}}{\leftarrow} \cdot \overset{x_{j}}{\leftarrow} \cdot)(\mathbf A_N) = A_{j'} \times A_j$.
			
	\item  {\bf Transpose:} For $(\cdot \overset{x_j}{\rightarrow} \cdot)$ the graph with two vertices $in$ and $out$ and one edge from $out$ to $in$ labelled $x_j$, $(\cdot \overset{x_j}{\rightarrow} \cdot)(\mathbf A_N)$ is the transpose $A_j^t$ of $A_j$.
					 
	\item  {\bf Projection on the diagonal:} For $(\circlearrowleft^{x_j})$ consisting in a single vertex $in=out$ and one loop labelled $x_j$, then $(\circlearrowleft^{x_j})(\mathbf A_N)$ is the diagonal matrix of diagonal elements of $A_j$. We denote it $\Delta(A_j)$. Note that $\Delta$ is a projection.
		
	\item  {\bf Degree:} For $\big( \underset {^.} \downarrow^{x_j}\big)$ consisting in two vertices $in=out$ and $v$ and an edge from $v$ to $in=out$ labelled $x_j$, then 
		$$\big( \underset {^.} \downarrow^{x_j}\big)(\mathbf A_N)=\mathrm{diag}(\sum_{m=1}^N A_{n,m})_{n=1\etc N}$$
 is the diagonal matrix, that we denote $deg(A_j)$, whose $n$-th diagonal element is the sum of the entries of $A_j$ on its $n$-th row. It is also a projection.
		
	\item {\bf Entry-wise products:} For $(\cdot \overset{x_j}{\underset{x_{j'}}\leftleftarrows} \cdot)$ with two vertices $in$ and $out$, two edges from $in$ to $out$, one labelled $x_j$ and the other one labelled $x_{j'}$, one has $(\cdot \overset{x_j}{\underset{x_{j'}}\leftleftarrows} \cdot)(\mathbf A_N) = A_{j'} \circ A_j$, where $\circ$ denotes the entry-wise product of matrices (also known as Hadamard or Schur product)

	\item {\bf Complex transpose:} For $( \cdot \overset{x_j^*}{\leftarrow}  \cdot)$ with two vertices $in$ and $out$ and one edge from $in$ to $out$ labelled $x_j^*$ then $( \cdot \overset{x_j^*}{\leftarrow}  \cdot)(\mathbf A_N) = A_j^*$. Note also that for any $^*$-graph monomial $g$, we can write $\big( g(\mathbf A_N) \big)^* = g^*(\mathbf A_N)$ where $g^*$ is obtained by interchanging the input and the output of $g$, reversing the orientation of its edges and reversing labels $x_j$ and $x_j^*$. We also denote $g^t$ obtained similarly without reversing labels $x_j$ and $x_j^*$, which satisfies $\big( g(\mathbf A_N) \big)^* = g^t(\mathbf A_N^*)$.
\end{enumerate}
\end{example}
The following lemma tells that $^*$-graph polynomials permute with the action by conjugation by permutation matrices, which is a crucial ingredient for our main Theorem \ref{MainTh}.

\begin{lemma}\label{Lem:PermInvDesTraces} Let $\mathbf A_N=(A_j)_{j\in J}$ be a family of $N$ by $N$ matrices. For any permutation matrix $V_N$, denote $V_N\mathbf A_NV_N^*:=(V_NA_jV_N^*)_{j\in J}$. Then, for any $^*$-graph operation $T$ and any permutation matrix $V_N$, one has $g( V_N\mathbf A_N V_N^*\big) = V_N g(\mathbf A_N) V_N^*.$ In particular, $\mathbf A_N$ and $V_N \mathbf A_N V_N^*$ have the same traffic distribution.
\end{lemma}
\begin{proof} If $\sigma$ denotes the permutation associated to $V_N$, then the entry $(n,m)$ of a matrix $V_N A_j V_N^*$ is $A_j\big( \sigma^{(-1)}(n), \sigma^{(-1)}(m)\big)$. Hence, we obtain the results thanks to the change of variable $\phi' = \sigma^{(-1)} \circ \phi$ below:
	\eq
		\Big( g( V_N \mathbf A_N V_N^*) \Big) (i,j) & = & \sum_{\substack{ \phi:V\to [N] \\ \phi(in)=j, \, \phi(out)=i}}  \prod_{e=(v,w)\in E}   A_{\gamma(e)}^{\varepsilon(e)}\big( \sigma^{(-1)} \circ \phi(w), \sigma^{(-1)}\circ\phi(v) \big) \\
			& = &  \sum_{\substack{ \phi':V\to [N] \\ \phi'(in)=\sigma^{(-1)}(j), \\ \phi'(out)=\sigma^{(-1)}(i)}}  \prod_{e=(v,w)\in E} A_{\gamma(e)}^{\varepsilon(e)}\big(  \phi'(w), \phi'(v) \big)\, \\
			& = &\Big( V_N g(\mathbf A_N) V_N^*\Big) (i,j).
	\qe
\end{proof}

Note that $g(U_N \mathbf A_N U_N^*)\neq U_N g(\mathbf A_N) U_N^*$ in general for an arbitrary unitary matrix $U_N$. Hence the traffic distributions of two matrices of a same linear map in two different basis can be different. For instance, the matrices $A=\left( \begin{array}{cc} 1 & 0 \\ 0 & -1 \end{array}\right)$ and $B=\left( \begin{array}{cc} 0 & 1 \\ 1 & 0 \end{array}\right)$ are both the matrices of a reflexion in $\mathbb C^2$, but $\Delta(A)=A$ and $\Delta(B)=0$, $deg(A)=A$ and $deg(B)=I_2$, etc.

\subsection{Link with other notions of convergence}

The traffic distribution of $\mathbf A_N$ encodes its $^*$-distribution, as well as many other statistics on the matrices.

\begin{example}\label{Ex:DistrFromTraffics}Let $\mathbf A_N$ be a family of random matrices. 

 \begin{enumerate}
	\item  The graph operations in matrices contain the $^*$-polynomials, and so the traffic distribution $\bar \Phi_{\mathbf A_N}$ of $\mathbf A_N$ encodes the $^*$-distribution 
	\eq
		\Phi_{\mathbf A_N}: P \mapsto \esp\big[\frac 1 N\Tr \, P(\mathbf A_N) \big]
	\qe
 by restriction  of $\bar \Phi$ on $^*$-graph polynomials consisting in linear combinations of simple lines $(\underset{out}\cdot \overset{x_{\gamma_1}^{\varepsilon_1}}\leftarrow \cdot \, \dots \, \cdot \overset{x_{\gamma_K}^{\varepsilon_K}}\leftarrow  \underset{in} \cdot)$.
		\item The traffic distribution of $\mathbf A_N$ induces the bilinear form 
	\eq
		(P,Q) \mapsto \esp\Big[ \frac 1 N \Tr \big[ P(\mathbf A_N) \circ Q(\mathbf A_N) \big]\Big] ,
	\qe
 where $\circ$ denotes the entry-wise product of matrices. 
	\item The diagonal matrix $deg(A_N)$ is a linear function of the matrix $A_N=(A_{i,j})_{i,j}$, so the map 
	\eq
		\Psi_{\mathbf A_N}:  P \mapsto \esp\Big[ \frac 1 N \Tr \,  deg\big(P(\mathbf A_N)\big)\Big]
	\qe
 is a linear form. It plays an important role later in the second part of the article, Chapter \ref{Sec:ThreeIndep}, in the context of asymptotic Boolean independence.
	\end{enumerate}
\end{example}

 Let us compare the convergence in traffic distribution, the convergence in $^*$-distribution and the convergence in spectral distribution.

\begin{lemma}\label{Lem:EquivCVTraf} Let $\mathbf A_N$ be a family of random matrices. The following are equivalent.
\begin{enumerate}
	\item The family $\mathbf A_N$ converges in traffic distribution, i.e.  $\bar \Phi_{\mathbf A_N}$ converges pointwise on $ \mathbb C \mathcal G\lara$.
	\item The family of matrices $\big( g(\mathbf A_N) \big)_{g}$ indexed by all $g\in \mathbb C \mathcal G\lara$ converges in $^*$-distribution.
	\item For any $^*$-graph polynomial $g\in \mathbb C \mathcal G\lara$ such that the matrix $g(\mathbf A_N)$ is Hermitian, the mean empirical spectral distribution of $g(\mathbf A_N)$ converges in moments.
\end{enumerate}
\end{lemma}

To prove the lemma we use the following fact, that will be important in all the paper. There is a composition law for graph operations consisting in replacing variables of a $^*$-graph monomial $g$ by $^*$-graph monomials.

\begin{definition}[Substitution of edges of graph monomials] \label{Def:SubsGraph} For a $^*$-graph monomial $g$ in the variables $(y_1 \etc y_n)$ and for $^*$-graph monomials $g_1 \etc g_n$ in the variables $\mathbf x$, the $^*$-graph monomial $g( g_1 \etc g_n)$ in the variables $\mathbf x$ is the graph obtained from $g$ by replacing each edge $e$ of $g$ labeled $y_i$ by the graph $g_i$, the input of $g_i$ identified with the source of $e$ and the output of $g_i$ with the target of $e$. \end{definition}

This operation of substitution is compatible with the evaluation of matrices: with $g, g_1 \etc g_n$ as above, for any families of matrices $\mathbf A_N$, we have 
	\eqa
		\label{Eq:SubstMatrices}
		g\big( g_1(\mathbf A_N), \dots g_n(\mathbf A_N) \big) = \big( g( g_1 \etc g_n) \big) \ (\mathbf A_N).
 	\qea

This property implies a lot of relations between the graph polynomials, that can be see easily by drawing associated graphs, without writing a formula in terms of the entries of the matrices. For instance, the relation 
	\eqa\label{Eq:ExampleAssoGraph}
		\Delta\big(A_N deg(B_N)\big) = deg\big( \Delta(A_N) B_N \big),
	\qea
 valid for any $A_N$ and $B_N$, is obtained in Figure \ref{fig:Subs} below.
 \begin{figure}
     \includegraphics{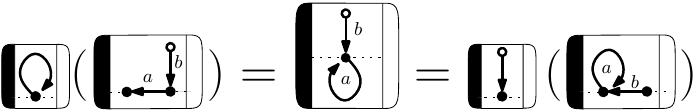}
    \caption{The substitution of edges of graph by $^*$-graph monomials that gives equality \eqref{Eq:ExampleAssoGraph}.}
    \label{fig:Subs}
  \end{figure}

\begin{proof}[Proof of Lemma \ref{Lem:EquivCVTraf}] Assume first (1) and let us prove (2). Let $g_1 \etc g_K$ be $^*$-graph monomials and let $M = y_{i_1}^{\varepsilon_1} \dots y_{i_L}^{\varepsilon_L}$ be a $^*$-monomial, where $i_\ell \in \{1\etc K\}$ and $\varepsilon_\ell\in \{1,*\}$, $\forall \ell=1\etc L$.  One can write
	$$M\big( g_1(\mathbf A_N), \dots , g_K(\mathbf A_N)\big) = g_{i_1}(\mathbf A_N)^{\varepsilon_1} \dots g_{i_L}(\mathbf A_N)^{\varepsilon_L}= \tilde g(\mathbf A_N)$$
 where $\tilde g =(\cdot \overset{y_1}\leftarrow \cdot \dots \cdot \overset{y_L}\leftarrow \cdot)
 (g_{i_1}^{\varepsilon_1} \etc g_{i_L}^{\varepsilon_L})$ with the notations of Definition \ref{Def:SubsGraph}. By assumption, $\mathbf A_N$ converges in traffic distribution so 
 	$$\esp\big[\frac 1 N \Tr M\big( g_1(\mathbf A_N), \dots , g_K(\mathbf A_N)\big)\big]=\esp\big[ \frac 1 N \Tr \, \tilde g(\mathbf A_N)\big]$$
 converges. Hence the convergence of $\big(g(\mathbf A_N)\big)_{g\in \mathbb C \mathcal G\lara}$ by multilinearity. 
\par We now prove that (2) implies (3). Consider a Hermitian matrix of the form $B_N = g(\mathbf A_N)$ where $g$ is a $^*$-graph polynomial. We can write $g = \sum_{i}\alpha_i g_i$ for a finite complex linear combination, where the $g_i$'s are $^*$-graph monomials.  For any integer $L\geq 1$, we can write $B_N^L=\big(g(\mathbf A_N)\big)^L = \sum_{i_1\etc i_L} \alpha_{i_1} \dots \alpha_{i_L} \tilde g_{i_1\etc i_L}(\mathbf A_N)$, where $\tilde g_{i_1\etc i_L} = (\cdot \overset{y_1}\leftarrow \cdot \dots \cdot \overset{y_L}\leftarrow \cdot)(g_{i_1} \etc g_{i_L})$. So the convergence of $\esp[\frac 1 N \Tr  \tilde g(\mathbf A_N)]$ for all $^*$-graph polynomials $\tilde g$ implies the convergence of $\esp[\frac 1 N \Tr B_N^L]$ for each $L$, and hence the convergence in moments of the empirical spectral distribution of $B_N$.

\par Let conclude by proving that (3) implies (1). For any $g\in \mathbb C \mathcal G\lara$, one can write $ g(\mathbf A_N) = g_1(\mathbf A_N) +  i g_2(\mathbf A_N)$ where  $g_1(\mathbf x) =\frac{   g(\mathbf x)  +    g^t(\mathbf x^*)}2 $ and $g_2 =  \frac{    g(\mathbf x) -  g^t(\mathbf x^*)}{2i} $. The matrices $g_1(\mathbf A_N)$ and $g_2(\mathbf A_N)$ are self-adjoint, so $\esp\big[\frac 1N \Tr \, g_i(\mathbf A_N)\big]$ converges for each $i=1,2$. Hence the convergence for $g(\mathbf A_N)$.
\end{proof}

\section{Statements of the results of Part 1}\label{Sec:MainTh}
\par We can now state the main result, omitting for the moment the explicit characterization of the limit.
\begin{theorem}[Asymptotic traffic independence]~\label{MainTh}~\\
Let $\mathbf A_N^{(1)} \etc \mathbf A_N^{(L)}$ be independent families of $N \times N$ random matrices. Denoting $\mathbf A_N^{(\ell)} = (A_j^{(\ell)})_{j\in J_\ell}$ for each $\ell=1\etc L$, we consider indeterminates $\mathbf x_\ell = (x_j^{(\ell)})_{j\in J_\ell}$.   Make the following hypotheses:
\begin{enumerate}
	\item  Each family (except possibly one) is permutation invariant (Definition \ref{Def:Intro}).
	\item  Each family converges in traffic distribution (Definition \ref{Def:StarGraphsMon}), namely: for any $\ell=1\etc L$ and any $^*$-graph polynomial $g \in \mathbb C \mathcal G \langle \mathbf x_\ell, \mathbf x_\ell^*\rangle$ one has
	\eqa
		 \bar \Phi_\ell (\, g\,  ) := \Nlim \esp\Big[ \frac 1 N \Tr\big[ g(\mathbf A_N^{(\ell)}) \big] \Big] \ \textrm{ exists.}
	\qea
	\item  Each family $\mathbf A_N^{(\ell)}$ satisfies the "factorization`` property: for all $\ell=1\etc L$, for any $^*$-graph polynomials $g_1,\dots,g_K \in \mathbb C \mathcal G \langle \mathbf x_\ell, \mathbf x_\ell^*\rangle$, $K\geq 2$,
	\eqa\label{eq:DecorrIntro}
		\esp\Big[ \prod_{k=1}^K \frac 1 N \Tr\big[ g_k(\mathbf A_N^{(\ell)}) \big] \Big] & \limN & \prod_{k=1}^K \bar \Phi( \,  g_k\, ) 
		  .
	\qea
	\end{enumerate}
Then the families $\mathbf A_N^{(1)} \etc \mathbf A_N^{(L)}$ are asymptotically traffic independent, that is:
	\begin{enumerate}
		\item They have a limiting joint traffic distribution: for any $^*$-graph polynomial $g\in \mathbb C \mathcal G \langle \mathbf x_\ell, \mathbf x_\ell^*\rangle_{\ell=1}^L$,
			\eqa
				\bar \Phi \big( \, g \, \big) := \Nlim \esp\Big[ \frac 1 N \Tr \big[ g(\mathbf A_N^{(1)} \etc \mathbf A_N^{(L)}) \big] \Big] \ \textrm{ exists.}
			\qea
		In particular, $\mathbf A_N=\mathbf A_N^{(1)} \cup \dots \cup \mathbf A_N^{(L)}$ has a limiting joint $^*$-distribution.
		\item The limiting traffic distribution $\bar \Phi$ of $\mathbf A_N$ depends only on the marginal limiting traffic distributions $\bar \Phi_\ell$ of the $\mathbf A_N^{(\ell)}$'s. It is called the product of the distributions $\bar \Phi_\ell$ and is given explicitly in Definition \ref{def:FreeProdGraphs}.
	\end{enumerate}
Moreover, the family of matrices $\mathbf A_N$ also satisfies the factorization property.
\end{theorem}
The proof of the theorem is given in Section \ref{Sec:ProofMainTh}
\\
\par Theorem \ref{MainTh} is applied for the matrix models of Definitions \ref{Def:Intro} and \ref{Def:WigMatrices}. In particular, it extends in a unified way known results of asymptotic $^*$-freeness in the setting of Theorem \ref{MainTh}.

\begin{corollary}\label{Cor:ApplClass} Let consider $\mathbf D_N, \mathbf V_N , \mathbf U_N, \mathbf X_N, \mathbf A_N$, independent families of random matrices  where
	\begin{itemize}
		\item $\mathbf D_N$ are independent diagonal matrices with independent and identically distributed diagonal entries whose moments of all orders exist.
		\item $\mathbf V_N$ are independent permutation matrices,
		\item $\mathbf U_N$ are independent unitary Haar matrices,
		\item $\mathbf X_N$ are independent Wigner matrices whose entries are invariant in law by complex conjugation,
		\item $\mathbf A_N$ is a family of random matrices satisfying the assumptions (2) and (3) of Theorem \ref{MainTh}.
	\end{itemize}
Then the matrices of $\mathbf D_N, \mathbf V_N , \mathbf U_N, \mathbf X_N$ and the family $\mathbf A_N$ are asymptotically traffic independent. Moreover, the matrices of $\mathbf U_N, \mathbf X_N$ and the family $\mathbf D_N\cup\mathbf V_N\cup \mathbf A_N$ are asymptotically free independent.
 \end{corollary}
 Corollary \ref{Cor:ApplClass} is the consequence of several results stated in Chapter 3. We summarize its proof in Section \ref{Sec:LinkTens}. 
\\
\par We prove two general criterions that are useful for a large class of matrices.

\begin{proposition}\label{Prop:RigAndCrit} Let $\mathbf A_N$ and $\mathbf B_N$ be two asymptotically traffic independent matrices.
\begin{enumerate}
	\item Denote $\Phi_N(A_N)=\esp\big[\frac 1 N \Tr \, A_N \big]$. If there exist two $^*$-polynomials $P,Q$ such that 
	\eq
		\lefteqn{\mathfrak K\big(P,Q,(\mathbf A_N)_{N\geq 1}\big)}\\
		& := &  \Nlim\Big(\Phi_N\big[ P( \mathbf A_N) \circ Q(\mathbf A_N)\big]-\Phi_N\big[ P( \mathbf A_N)\big] \Phi_N\big[ Q(\mathbf A_N)\big]\Big)
	\qe
	 is nonzero, and the same holds for $\mathbf B_N$, then $\mathbf A_N$ and $\mathbf B_N$ are not asymptotically free independent.
	 \item If $\mathbf A_N$ has the same limiting traffic distribution as a unitarily invariant families of matrices, then $\mathbf A_N$ and $\mathbf B_N$ are asymptotically free independent.
\end{enumerate}
\end{proposition}

\begin{remark} A partial reciprocal is true for this criterion: for two asymptotically traffic independent families $\mathbf A_N$ and $\mathbf B_N$, if for any $P,Q$ one has 
$$\mathfrak K\big(P,Q,(\mathbf A_N)_{N\geq 1}\big)=\mathfrak K\big(P,Q,(\mathbf B_N)_{N\geq 1}\big)=0,$$
then $\mathbf A_N$ and $\mathbf B_N$ they are asymptotically free independent. This is proved in \cite{CDM16} and relies on an equivalent formulation of traffic independence. If for any $P,Q$  one has $\mathfrak K\big(P,Q,(\mathbf A_N)_{N\geq 1}\big)=0$ but for some $P,Q$ one has$\mathfrak K\big(P,Q,(\mathbf B_N)_{N\geq 1}\big)\neq 0$, then different scenarios are possible.
\end{remark}

The first point is proved in Section \ref{Sec:CritLackInp} and the second point in Section \ref{Sec:LinkTens}.

 Theorem \ref{MainTh} and these two propositions are used in \cite{MP14,DLM15, MAL122} for several matrix models. We apply Proposition \ref{Prop:RigAndCrit} in the context of Corollary \ref{Cor:ApplClass} as follow. 

\begin{corollary}\label{Cor:ExIntro}~ \begin{itemize}
\item {\bf Asymptotic free independence with the transpose:} Independent complex Wigner matrices $X_N^{(j)},j\in J$, with parameter of the form  $(\alpha_j, 0)$, independent Haar matrices $U_N^{(j')}, j\in J'$, and the transposed matrices $X_N^{(j)t}, U_N^{(j)t},j\in J, j'\in J$ are asymptotically free independent.
	\item {\bf Entry wise product of matrices:}
 Let $W_N=(\omega_{i,j})_{i,j=1\etc N}$ be a random matrix with i.i.d. entries whose moments are finite and independent of $N$. For an independent unitary Haar matrix $U_N$ and an independent uniform permutation matrix $V_N$, consider the entry-wise products $M_N = W_N \circ U_N$ and $\tilde M_N = W_N \circ V_N$. If $\omega_{i,j}$ is centered then independent copies of $M_N$ are asymptotically free circular elements. If the modulus of $\omega_{i,j}$ is not constant then independent copies of $\tilde M_N$ are not asymptotically free. 
\end{itemize}
\end{corollary}

\par {\bf Illustration: the convolutions.} Let us point out some cases of applications for limiting distribution of the sum of two matrices. Consider two independent Hermitian random matrices $ A_N$ and $B_N$. Assume that $A_N$ is diagonal with independent entries identically distributed according to some distribution $\mu_a$. Note that $A_N$ is permutation invariant. Assume now that $B_N $ converges in traffic distribution. In particular its mean empirical spectral distributions converges to some measure $\mu_b$. It follows from Theorem \ref{MainTh}, Corollary \ref{Cor:ApplClass} and Remark \ref{Rk:Rk:MainThBack}, that $A_N$ and $B_N$ are asymptotically traffic independent. 

In particular, the empirical spectral distribution of the matrix $H_N = A_N + B_N$ converges to some measure $\mu_h$, which depends on $\mu_a$, $\mu_b$, and on the limiting traffic distribution of $B_N$.

\begin{enumerate}
	\item If $B_N$ is a diagonal matrix, then $\mu_h$ is the classical convolution $\mu_a * \mu_b$ of $\mu_a$ and $\mu_b$, that is the distribution of the sum of a random variable distributed according to $\mu_a$ and an independent (in the classical sense) random variable distributed according to $\mu_b$. Indeed it is straightforward to see that the empirical spectral distribution of $H_N$ is actually the convolution of the empirical spectral distribution of $A_N$ and $B_N$ for each $N\geq 1$. 
	\item If $B_N$ is unitary invariant, then $\mu_h$ is the so-called free convolution $\mu_a \boxplus\mu_b$ of Voiculescu by consequence of their asymptotically free independence. In the setting of free probability \cite{NS}, it is the distribution of the sum of a non random variable distributed according to $\mu_a$ and a free independent random variable distributed according to $\mu_b$. 
	\item For many examples of matrix models, it happens that $A_N$ and $B_N$ are asymptotically traffic independent but the limiting distribution of $A_N + B_N$ is neither the classical or the free convolution (see for example \cite{MAL122}. 
\end{enumerate}

\par In conclusion, at the level of the limiting $^*$-distribution of large matrices, traffic independence encodes both the free and the classical notions of independence, but also a large class of operations.

\chapter{Definition of Asymptotic Traffic Independence}\label{Sec:DefIndepTraffic}

We first introduce some combinatorial tools to manipulate the traffic distributions. Then we define asymptotic traffic independence, and prove natural properties we can expect from a notion of independence. We illustrate this notion with a computation, which leads to a criterion of lack of asymptotic free independence for matrices. We conclude this section with the proof of Theorem \ref{MainTh}.

\section[Combinatorial distributions and injective version]{Combinatorial form of traffic distributions and injective version}\label{Sec:InjTrace}

 Recall that a $^*$-test graph in variables $\mathbf x=(x_j)_{j\in J}$ is a collection $(V,E, \gamma, \varepsilon)$ where $(V,E)$ is a finite connected graph with at least one vertex (but possibly no edges), whose edges are labelled by symbols $x_j$ and $x_j^*$. An edge $e$ has label $x_{\gamma(e)}^{\varepsilon(e)}$, where $\gamma: E \to J$ and $\varepsilon: E\to \{1,*\}$. 
 
 \begin{definition} The set of $^*$-test graphs in variables $\mathbf x$ is denoted $\mathcal T\lara$ and the space of finite linear complex combinations of elements of $ \mathcal T \lara$ is denote by $\mathbb C \mathcal T\lara$.
 \end{definition}
Recall that a $^*$-graph monomial $g$ is the data of a $^*$-test graph $T_g$ and two vertices $in, out$ of $T_g$ and that the traffic distribution of $\mathbf A_N=(A_j)_{j\in J}$ the data of the map $\bar \Phi_{\mathbf A_N}:g\in \mathbb C \mathcal G\lara \mapsto \esp\big[ \frac 1 N \Tr \, g(\mathbf A_N) \big]\in \mathbb C.$ Given such a $^*$-graph monomial $g=(T_g, in , out)$, note that $\bar \Phi_{\mathbf A_N}(g) = \bar \Phi_{\mathbf A_N}\big(\Delta(g)\big) $ where $\Delta(g)$ is obtained from $g$ by identifying its input and output. Moreover, this quantity does not depend on the position of the input of $\Delta(g)$, but only on the $^*$-test graph $T$ such that $\Delta(g) = (T, in, in)$.

\begin{definition}[Combinatorial form of traffic distributions of matrices]~\label{Def:CombiDistrTraf} \begin{enumerate}
	\item For any family $\mathbf A_N = (A_j)_{j\in J}$ of matrices and any $^*$-graph monomial $g$ in variables $\mathbf x= (x_j)_{j\in J}$, we define $\tau_{\mathbf A_N}\big[ T] = \bar \Phi_{\mathbf A_N}(g)$ where $T$ is the $^*$-test graph obtained by identifying the input and output of $g$ and forgeting their position in the new $^*$-test graph. We use the notation $\tau_N\big[T(\mathbf A_N)\big]:=\tau_{\mathbf A_N}\big[ T]$.
	\item Equivalently, given $T=(V,E,\gamma,\varepsilon)$ a $^*$-test graph in variables $\mathbf x=(x_j)_{j\in J}$ and $\mathbf A_N = (A_j)_{j\in J}$ a family of random matrices, then $\tau_N\big[  T(\mathbf A_N)\big]:=  \esp\big[ \frac 1 N \Tr \, T(\mathbf A_N) \big]$ where
  	\eq
		\Tr \big[  T(\mathbf A_N)\big] :=\sum_{ \phi:V \to [N]} \prod_{e=(v,w)\in E} A_{\gamma(e)}^{\varepsilon(e)}\big( \phi(w), \phi(v) \big).
	\qe
Such a quantity is called the trace of the $^*$-test graph $T$ in the matrices $\mathbf A_N$. For any matrix $A_N$, denoting by $\circlearrowleft$ the $^*$-test graph with a single loop, we then have
	$$\Tr \, A_N = \Tr \big[ \circlearrowleft(A_N)\big].$$ 
 
\item The combinatorial distribution of $\mathbf A_N$ is the linear map
 	\eq
			\tau_{\mathbf A_N} : T \in  \mathbb C \mathcal T \lara \mapsto \tau_N\big[  T(\mathbf A_N)\big] \in \mathbb C.
	\qe
\end{enumerate}
 
 \end{definition}
 
\begin{remark}\label{Rk:CombiDistrTraf}
\begin{enumerate}
	\item Notation $\Tr \big[  T(\mathbf A_N)\big] $ and terminology are abusive since we do not define the object $T(\mathbf A_N)$ and $\Tr$ is not the trace of matrices. Hence the name \emph{combinatorial distribution} for $\tau_{\mathbf A_N}$. 
	\item As in Definition \ref{Def:SubsGraph}, for a $^*$-test graph in variables $y_1 \etc y_n$ and $^*$-graph monomials $g_1\etc g_n$, we define the $^*$-test graph $T(g_1 \etc g_n)$ by replacing edges labeled $y_i$ of $T$ by the graph $g_i$. Then we have the compatibility property 
	$$\tau_N\Big[ T\big( g_1(\mathbf A_N) \etc g_n(\mathbf A_N)\big)\Big] = \tau_N\Big[ T(g_1\etc g_n)(\mathbf A_N)\Big].$$
	In the l.h.s. the $^*$-test graph is $T$ and the matrices are the $g_i(\mathbf A_N)$, whereas in the r.h.s. the $^*$-test graph is $T(g_1\etc g_n)$ and the matrices are those of $\mathbf A_N$.
	\end{enumerate}
\end{remark}

\begin{example}\label{Ex:DistrFromTrafficsBis}Example \ref{Ex:DistrFromTraffics} continued. Let $\mathbf A_N$ be a family of random matrices.
\begin{enumerate}
	\item Let $M$ be a $^*$-monomial $M=x_{\gamma_1}^{\varepsilon_1} \dots x_{\gamma_n}^{\varepsilon_n}$. Then one has $\Tr \, M(\mathbf A_N) = \Tr \big[ T(\mathbf A_N) \big]$ where $T$ is the $^*$-test graph consisting in a simple oriented cycle, with vertex set $\{1 \etc n\}$ and edges $(i+1,i)$ labeled $x_{\gamma_i}^{\varepsilon_i}$, $i=1\etc n$ (with indices modulo $n$).
	\item Let $M'$ be a second $^*$-graph monomial and let $T'$ the $^*$-test graph with vertex set $\{1' \etc n'\}$ defined as $T$ is defined from $M$ in the previous item. Recalling that $\circ$ denotes the entry wise product of matrices, then $\Tr \big[ M(\mathbf A_N) \circ M'(\mathbf A_N) \big] = \Tr\big[ \tilde T(\mathbf A_N)\big]$ where $\tilde T$ is the $^*$-test graph consisting in the bunch of two simple oriented cycles obtained by identifying the vertices $1$ and $1'$ of $T$ and $T'$.
	\item Recall that $deg(A_N)$ is the diagonal matrix whose $i$-th diagonal matrix is the sum of the elements of $A_N$ over the $i$-th row. Then $\Tr \big[ deg(A_N)\big] = \Tr\big[ T(A_N)\big]$ where $T$ is the $^*$-test graph consisting in a simple edge.
\end{enumerate}
\end{example}

We can now introduce the following transformation of traffic distributions.
 
\begin{definition}[Injective trace of matrices]\label{Dej:InjDistr}Let $T=(V,E,\gamma,\varepsilon)$ be a $^*$-test graph in variables $\mathbf x=(x_j)_{j\in J}$ and let $\mathbf A_N = (A_j)_{j\in J}$ be a family of matrices, possibly random. We call injective trace of the $^*$-test graph $T$ in the matrices $\mathbf A_N$ the quantity
	\eqa \label{DefIntroGraphMon}
		 \Tr^0\big[ T(\mathbf A_N) \big] =  \sum_{
		\substack {
		\phi: V \to [N] \\ \textrm{injective}}
		} \prod_{e=(v,w)\in E} A_{\gamma(e)}^{\varepsilon(e)}\big( \phi(w), \phi(v) \big).
	\qea 
and we set $\tau^0_N\big[  T(\mathbf A_N)\big]:=  \esp\big[ \frac 1 N \Tr^0 \, T(\mathbf A_N) \big]$, where the expectation is relative to the matrices when they are random. The injective traffic distribution of $\mathbf A_N$ is the linear map 
	$$\tau_{\mathbf A_N}^0 : T \in \mathbb C \mathcal T \lara \mapsto  \tau_N^0\big[  T(\mathbf A_N) \big]\in \mathbb C.$$
\end{definition}

\begin{figure}
     \includegraphics{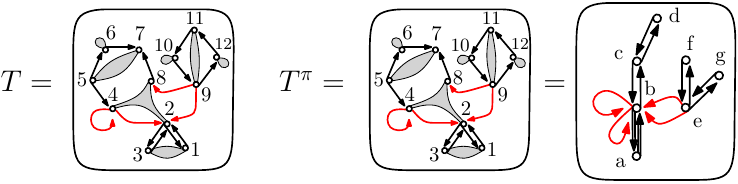}
    \caption{Left: a test graph with vertex set $\{1\etc 12\}$. Right: the quotient graph $T^\pi$ for the partition $\pi=\big\{ a=\{1,3\}, b=\{2,4,8\}, c=\{5,7\}, d=\{6\}, e=\{9,11\},f= \{10\}, g=\{12\} \big\}$.}
    \label{fig:ConstrTauPi}
  \end{figure}

Combinatorial distributions and their injective versions are related each other. Let $T$ be a $^*$-test graph with vertex set $V$ and let $\mathcal P(V)$ denote the set of partitions of $V$. For any $\pi\in \mathcal P(V)$, we denote by $T^\pi$ the $^*$-test graph obtained by identifying vertices in a same block of $\pi$ (the edges link the associated blocks). See an example Figure \ref{fig:ConstrTauPi}. Then for any $^*$-test graph $T$, one has

	\eqa\label{eq:TraffCum}
		 \Tr \big[ T(\mathbf A_N)\big] = \sum_{\pi \in \mathcal P(V)} \Tr^0 \big[ T^\pi(\mathbf A_N)\big].
	\qea

In this formula, the different cases of equality for the indices of the matrices in the definition of $\Tr \big[ T(\mathbf A_N)\big]$ are classified by choosing the partition of indices whose values are equal.
Reciprocally, it is possible to write $  \Tr^0$ in terms of $  \Tr$ using a generalized inclusion-exclusion principle (see \cite{NS,Stan12}).

\begin{lemma}\label{Lem:FormuleTraceTraceInj} Let $T$ be a $^*$-test graph with vertex set $V$. For any partition $\pi$ of $V$, denote $\mu_V(\pi) = \prod_{ B \in \pi}(-1)^{|B|-1}\big( |B|-1\big)!$. Then one has, for any family of matrices $\mathbf A_N$,
	\eqa\label{Eq:InversCum}
		\Tr^0 \big[ T(\mathbf A_N)\big] = \sum_{\pi \in \mathcal P(V)} \mu_V(\pi)  \Tr \big[ T^\pi(\mathbf A_N)\big].
	\qea
\end{lemma}

\begin{proof}Note that $\mathcal P(V)$ is a partially ordered finite set where $\pi \leq \pi'$ if and only if the blocks of $\pi$ contain the blocks of $\pi'$. Moreover, for any $T$ and any partition $\pi$ one has $\Tr \big[ T^\pi(\mathbf A_N)\big] = \sum_{\pi'\geq \pi} \Tr^0 \big[ T^{\pi'}(\mathbf A_N)\big]$. Hence by dual M\"obius inversion formula \cite[Lecture 10]{NS}, there exists a map $\mu_V:\mathcal P(V) \to \mathbb Z$ such that \eqref{Eq:InversCum} holds. One has $\mu_V(\pi) = \mu_{\mathcal P(V)}(0_{V}, \pi)$, where $ \mu_{\mathcal P(V)}$ is the M\"obius map associated to the poset $\mathcal P(V)$ and $0_{V}$ is the partition with only singletons. The values of $\mu_V$ are given in \cite[Example 3.10.4]{Stan12}.
\end{proof}

This relation between the trace and the injective trace of $^*$-test graph is used below as a definition in the general case.

\begin{definition}[Injective trace, general case]\label{Def:TrTrZERO} Let $\tau: \mathbb C \mathcal T\lara \to \mathbb C$ be a linear map. We call injective version of $\tau$ the linear map $\tau^0: \mathbb C \mathcal T\lara \to \mathbb C$ defined for any $^*$-test graph $T$ by
\eqa\label{Eq:InversCum2}
		\tau^0  [ T ] = \sum_{\pi \in \mathcal P(V)} \mu_V(\pi)  \tau  [ T^\pi],
	\qea
where $\mu_V$ and $T^\pi$ are as in \eqref{Eq:InversCum}, which implies that for any $^*$-test graph $T$
\eqa\label{Eq:InversCum3}
		\tau   [ T ] = \sum_{\pi \in \mathcal P(V)}  \tau^0  [ T^\pi].
	\qea
\end{definition}
\begin{remark} \begin{enumerate}
	\item There is no meaning for an expression like $\Tr^0 \, A_N$ or $\Tr^0 \, g(\mathbf A_N)$, the injective trace is always defined for $^*$-test graphs.
	\item A relation exists between the injective trace and the notions of free cumulants, but in the particular case of unitarily invariant matrices. It is the main motivation of \cite{CDM16} to state and use this relation.
	\item With notations as in the second item of Remark \ref{Rk:CombiDistrTraf}, in general we have
	$$\tau^0_N\Big[ T\big( g_1(\mathbf A_N) \etc g_n(\mathbf A_N)\big)\Big] \neq \tau^0_N\Big[ T(g_1\etc g_n)(\mathbf A_N)\Big],$$
	 see Lemmas \ref{Lem:SubsInjTrace0}.
	 
\end{enumerate}
\end{remark}

Let us formulate the assumptions of Theorem \ref{MainTh} in terms of the injective trace. 

\begin{lemma}\label{Lem:EqDeccorProp} Let $\mathbf A_N$ be a family of random matrices. There is equivalence between
\begin{enumerate}
	\item The convergence in traffic distribution of $\mathbf A_N$, namely the pointwise convergence of $\bar \Phi_{\mathbf A_N}: g \in \mathbb C \mathcal G \lara \mapsto \esp\big[\frac 1 N \Tr \, g(\mathbf A_N)\big]$,
	\item the pointwise convergence of $\tau_{\mathbf A_N}:  T \in \mathbb C \mathcal T \lara  \mapsto \esp\big[\frac 1 N \Tr \,T(\mathbf A_N)\big]$
	\item  the pointwise convergence of $\tau^0_{\mathbf A_N}:  T \in \mathbb C \mathcal T \lara  \mapsto \esp\big[\frac 1 N \Tr^0 \,T(\mathbf A_N)\big]$
\end{enumerate}

 Assuming this convergence, then $\mathbf A_N$ satisfies the factorization property, assumption (3) of Theorem \ref{MainTh}, if and only if for any $^*$-test graphs $T_1,\dots,T_K$, $K\geq 2$,
	\eqa\label{eq:DecorrInj}
		 \esp\Big[ \prod_{k=1}^K \frac 1 N \Tr\big[ T_k(\mathbf A_N ) \big] \Big] \limN \prod_{k=1}^K \tau[\, T_K\, ],
	\qea
which is also equivalent to the same property for the injective trace, namely for any $^*$-test graphs $T_1,\dots,T_K $, $K\geq 2$,
	\eqa
		  \esp\Big[ \prod_{k=1}^K \frac 1 N \Tr^0\big[ T_k(\mathbf A_N ) \big] \Big]  \limN \prod_{k=1}^K \tau^0[\, T_K\, ].
	\qea

\end{lemma}

\begin{proof} The equivalence of the formulations in terms of $\bar \Phi_{\mathbf A_N}$ and $\tau_{\mathbf A_N}$ is clear by definition. The equivalence between convergence in traffic distribution and pointwise convergence of $\tau^0_{\mathbf A_N}$ is a consequence of Formulas \eqref{DefIntroGraphMon} and \eqref{eq:TraffCum}. Let $T_1,\dots,T_K $ be $^*$-test graphs whose vertex sets are denoted $V_1\etc V_K$ respectively. The factorization property implies
	\eq
		  \esp\Big[ \prod_{k=1}^K \frac 1 N \Tr^0\big[ T_k(\mathbf A_N ) \big] \Big] &=&  \sum_{\substack{ \pi_k \in \mathcal P(V_k) \\ k=1\etc K}} \prod_{k=1}^K \mu_{V_k}(\pi_k)   \esp\Big[ \prod_{k=1}^K \frac 1 N \Tr\big[ T_k^\pi(\mathbf A_N ) \big] \Big] \\
		  & =&  \sum_{\substack{ \pi_k \in \mathcal P(V_k) \\ k=1\etc K}} \prod_{k=1}^K \mu_{V_k}(\pi_k) \bigg ( \prod_{k=1}^K \tau_N\big[ T_k^\pi(\mathbf A_N ) \big]\\
		  & & \ \  + \varepsilon(T^{\pi_k}_k)_{k=1\etc K} \bigg),
	\qe
where for each $k=1\etc K$ and $N$ large enough the term $\varepsilon(T^{\pi_k}_k)\leq 1$ tends to zero. Since $\tau_N\big[ T_k^\pi(\mathbf A_N ) \big] $ is bounded for each $\pi_1\etc \pi_k$, we get
	\eq
		  \esp\Big[ \prod_{k=1}^K \frac 1 N \Tr^0\big[ T_k(\mathbf A_N ) \big] \Big] & =&  \bigg(\sum_{\substack{ \pi_k \in \mathcal P(V_k)\\ k=1\etc K}} \prod_{k=1}^K \mu_{V_k}(\pi_k) \tau_N\big[ T_k^\pi(\mathbf A_N ) \big] \bigg) + \tilde \varepsilon\\
		  & = &  \prod_{k=1}^K  \tau_N^0\big[ T_k^\pi(\mathbf A_N ) \big]+ \tilde \varepsilon
	\qe
where $\tilde \varepsilon$ is bounded by a constant times the maximum of the $ \varepsilon(T^{\pi_k}_k)$ for any $\pi_k \in \mathcal P(V_k), k=1\etc K$. The reciprocal is true by the same computation where, starting with the trace instead of the injective one, we omit the terms $\mu_V(\pi_k)$ in the above computation.
\end{proof}

\section{Asymptotic traffic independence, properties}	\label{Sec:TrafficIndep}

\begin{definition}[Graph of colored components]~\label{Def:GCC}
\begin{enumerate}
	\item Let $T$ be a $^*$-test graph in variables $\mathbf x = \mathbf x_1\cup \dots \cup \mathbf x_p$, where the $\mathbf x_j$'s are families of pairwise disjoint variables (a variable appears at most in one family). A \emph{colored component} of $T$ with respect to $\mathbf x_1 \etc \mathbf x_p$ is a maximal connected subgraph of $T$, with at least one edge, whose edges are labeled by variables in only one family among $\mathbf x_1 \etc \mathbf x_p$. We denote by $\mathcal C\mathcal C(T)$ the set of colored components of $T$ with respect to $\mathbf x_1 \etc \mathbf x_p$.
	\item The graph of colored components of $T$ with respect to $\mathbf x_1 \etc \mathbf x_p$, denoted $\mathcal G\mathcal C\mathcal C(T)$, is the following bipartite undirected graph.
	\begin{itemize}
		\item The first kind of vertices are the colored components $T_1 \etc T_K$ of $T$.
		\item The second kind of vertices are the vertices $v_1 \etc v_L$ of $T$ that belong to at least two graphs among $T_1 \etc T_K$.
		\item There is an edge between $T_i$ and $v_j$ if $v_j$ is a vertex of $T_i$, $i=1\etc K$, $j=1\etc L$.
	\end{itemize}
\end{enumerate}
\end{definition}
\par
\begin{figure}
     \includegraphics{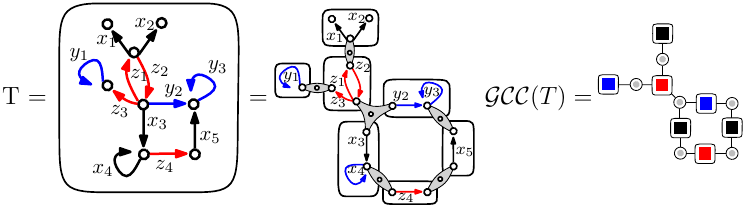}
    \caption{An example of construction of $\bar T$. The graph $\mathcal G \mathcal C \mathcal C(T)$ is not a tree since $\bar T$ has a cycle. Removing the edge labelled $z_4$ or $x_5$ in $T$, the graph of colored components becomes a tree.}\label{fig:DefFalseFreeness}
  \end{figure}

 In the leftmost picture of Figure \ref{fig:DefFalseFreeness}, we draw a $^*$-test graph in three families, represented by three colors $\mathbf x, \mathbf y$ and $\mathbf z$. In the intermediate figure, we have encircled the colored components of $T$, and drawn in the rightmost one its graph of colored components. 

\begin{definition}[Asymptotic traffic independence]~\label{def:FreeProdGraphs}
For each $\ell$ in $\{1\etc L\}$, let $\tau_\ell:\mathbb C\mathcal T \langle \mathbf x_\ell ,  \mathbf x_\ell^* \rangle \to \mathbb C$ be a linear map sending the graph with no edge to one. Up to a renaming of the symbols, we assume that the variables of the different families $\mathbf x_\ell$ for $\ell=1\etc L$ are pairwise distinct and we denote $\mathbf x = \mathbf x_1 \cup \dots \cup \mathbf x_L$.

The product\footnote{This product should be called ''free product`` as it satisfies a universal property, see \cite{CDM16}. Nevertheless we do not use this term to avoid confusion with free independence.} of the maps $\tau_1 \etc \tau_L$ is the linear map $\tau: \mathbb C \mathcal T \langle \mathbf x  ,  \mathbf x\rangle \to \mathbb C$ defined by the following:
\begin{enumerate}
	\item for each $\ell=1\etc L$ and any $^*$-test graph $T$ in the variables $\mathbf x_\ell$, $\tau [T ] = \tau_\ell[T]$ with the abuse of notation that $\mathbf x$ are considered as the variables of $T$, with no edges labeled $\mathbf x_{\ell'}, \ell'\neq \ell$.
	\item 
for any $^*$-test graph $T$ in the variables $\mathbf x_1 \etc \mathbf x_p$,
		\eqa\label{DefTrafInd}
			\tau^0[T] = \one\big( \mathcal G\mathcal C \mathcal C (T) \mathrm{ \ is \ a \ tree\,} \big)  \prod_{S \in \mathcal C \mathcal C(T)}\tau^0[S],
		\qea
	where $\mathcal G \mathcal C \mathcal C(T)$ (resp. $\mathcal C \mathcal C(T)$) is the graph (resp. the set) of colored components of $T$ with respect to $\mathbf x_1 \etc \mathbf x_p$.
\end{enumerate}
 We say that $\mathbf A_N^{(1)} \etc \mathbf A_N^{(L)}$ are asymptotically traffic independent whenever $\mathbf A_N:= \mathbf A_N^{(1)} \cup \dots \cup \mathbf A_N^{(L)}$ converges in traffic distribution and its combinatorial distribution converges to the product of the limiting distributions of the $\mathbf A_N^{(\ell)}$'s. 
\end{definition}

Hence the asymptotic traffic independence of $\mathbf A_N^{(1)} , \dots , \mathbf A_N^{(L)}$ is equivalent to the convergence for any $^*$-test graph $T$ with vertex set $V$,
	\eqa
		\lefteqn{\tau_N\Big[ T\big(\mathbf A_N^{(1)} , \dots , \mathbf A_N^{(L)}\big) \Big] }\nonumber \\ 
		\label{Eq:DefIndepTrace}   & \limN &
	\sum_{ \substack{\pi \in \mathcal P(V) \\ \mathrm{s.t.} \ \mathcal G \mathcal C \mathcal C(T^\pi) \ \\\mathrm{is \ a \ tree}}}    \prod_{S \in \mathcal C \mathcal C(T^\pi)}\Nlim \tau_N^0\Big[S\big(\mathbf A_N^{(\ell(S))}\big) \Big],
	\qea
where $\ell(S)$ is the index of the labels of the colored component $S$. Since the injective traces in the product can also be written in terms of (non-injective) traces, asymptotic traffic independence entirely characterizes the limiting traffic distribution of $\mathbf A_N^{(1)}\cup \dots \cup \mathbf A_N^{(L)}$ in terms of the traffic distributions of the $\mathbf A_N^{(\ell)}$'s.

Specified for graphs $T$ consisting in simple oriented cycles (Example \ref{Ex:DistrFromTrafficsBis}), Formula \eqref{Eq:DefIndepTrace} gives an expression for the limiting $^*$-distribution of $\mathbf A_N^{(1)} \cup \dots \cup \mathbf A_N^{(L)}$.

\begin{proposition}\label{Lem:Associativity} The product of Definition \ref{def:FreeProdGraphs} is symmetric and associative, in the following sense. Let $\tau_1, \tau_2,\tau_3$ denote three linear maps $\tau_i: \mathbb C \mathcal T \langle \mathbf x_i  ,  \mathbf x \rangle \to \mathbb C$ sending the graph with no edge to one. Then the product of $\tau_1$ and $\tau_2$ is the product of $\tau_2$ and $\tau_1$ (interchanging the variables). Moreover, the product $\tau$ of $\tau_1, \tau_2,\tau_3$ is the product $\tau_{3|1,2}$ of $\tau_3$ with the product $\tau_{1|2}$ of $\tau_1$ and $\tau_2$. In particular, three large families of matrices $\mathbf A_N, \mathbf B_N, \mathbf C_N$ are asymptotically traffic independent whenever  $\mathbf A_N$ and $ \mathbf B_N$ are asymptotically traffic independent and $\mathbf A_N\cup \mathbf B_N$ and $ \mathbf C_N$ are asymptotically traffic independent.
\end{proposition}

We shall need the following lemma, now and in the sequel of the article.

\begin{lemma}[Number of edges and vertices in a connected graph]~\label{lem:EdgesVertices}
\\Let $\mathcal G=(\mathcal V,\mathcal E)$ be a finite connected graph. Then, one has 
		$|\mathcal V|\leq |\mathcal E|+1, $
with equality if and only if $\mathcal G$ is a tree.
\end{lemma}

\begin{proof}[Proof of Lemma \ref{lem:EdgesVertices}]
If $G$ is a tree, we count the difference between the number of edges and vertices by removing its branches successively until it remains a single vertex, this yields the formula. If $G$ is not a tree, then we can remove edges until it becomes a tree which yields the inequality. 
\end{proof}

\begin{proof}[Proof of Proposition \ref{Lem:Associativity}] The symmetry is a immediate consequence of the definition. Let $T$ be a $^*$-test graph in the indeterminates $\mathbf x = \cup_{i=1}^3 \mathbf x_i$. Denote by $\mathcal G = (\mathcal V, \mathcal E)$ the graph of colored components ($\mathcal G\mathcal C \mathcal C$) of $T$ with respect to $\mathbf x_1, \mathbf x_2, \mathbf x_3$. Denote by $\mathcal G' = (\mathcal V', \mathcal E')$ the $\mathcal G\mathcal C \mathcal C$ of $T$ with respect to $\mathbf x_1\cup \mathbf x_2$ and $\mathbf x_3$. For each colored component $S$ of $\mathcal G'$ with label in $\mathbf x_1\cup \mathbf x_2$, we denote by $\mathcal G(S) = (\mathcal V(S), \mathcal E(S))$ the $\mathcal G\mathcal C \mathcal C$ of $S$ with respect to $\mathbf x_1$ and $\mathbf x_2$. Denote by $\mathcal C \mathcal C (\mathbf x_i)$, $i=1\etc 3$, the set of colored components of $T$ labeled $\mathbf x_i$.

We prove that $\mathcal G$ is a tree if and only if $\mathcal G'$ is a tree and the $\mathcal G(S)$'s are trees. Denote by
\begin{enumerate}
	\item $a$ the number of vertices of $T$ that belong both to an element of $\mathcal C \mathcal C (\mathbf x_3)$ and to a single other element of $\mathcal C \mathcal C (\mathbf x_1) \cup \mathcal C \mathcal C (\mathbf x_2)$,
	\item $b$ the number of vertices of $T$ that belong to an element of $\mathcal C \mathcal C (\mathbf x_3)$, of $\mathcal C \mathcal C (\mathbf x_1)$ and of $\cup \mathcal C \mathcal C (\mathbf x_2)$,
	\item $c$ the number of vertices of $T$ that belong both to an element of $\mathcal C \mathcal C (\mathbf x_1)$ and of $\mathcal C \mathcal C (\mathbf x_2)$, but not of $\mathcal C \mathcal C (\mathbf x_3)$,		\item $|\mathcal C \mathcal C(\mathbf x_1 \cup \mathbf x_2)|$ the number of colored components of $\mathcal G'$.
\end{enumerate}

One can enumerate easily the following quantities
	\eq
		|\mathcal V| = a + b + c+ \sum_{i=1}^3|\mathcal C \mathcal C(\mathbf x_i)|, \ \ |\mathcal E| = 2a+3b+2c,\\
	\sum_S|\mathcal V(S)| =  b + c+ \sum_{i=1}^2|\mathcal C \mathcal C(\mathbf x_i)|, \ \ \sum_S |\mathcal E(S)| = 2b+2c,\\
	|\mathcal V'| = a + b  + |\mathcal C \mathcal C(\mathbf x_3)| +|\mathcal C \mathcal C(\mathbf x_1 \cup \mathbf x_2)| , \ \ |\mathcal E'| = 2a+2b.
	\qe
Since $|\mathcal C \mathcal C(\mathbf x_1 \cup \mathbf x_2)|$ is the number of colored components $S$ of $\mathcal G'$ labeled in $\mathbf x_1\cup \mathbf x_2$, this results in the equality   
	$$|\mathcal V| - 1 -|\mathcal E| = \big( |\mathcal V'| - 1 -|\mathcal E'| \big) + \sum_S   \big( |\mathcal V(S)| - 1 -|\mathcal E(S)| \big),$$
telling by Lemma \ref{lem:EdgesVertices} that $\mathcal G$ is a tree if and only if $\mathcal G'$ is a tree and the $\mathcal G(S)$'s are trees. 

If we assume that $\mathcal G$, $\mathcal G'$ and the $\mathcal G(S)$ are trees, it is clear that 
	$$\prod_{i=1}^3\prod_{S \in \mathcal C \mathcal C( \mathbf x_i)} \tau^0[S] = \prod_{S \ \mathcal C\mathcal C( \mathbf x_3) }\tau^0[S] \times \prod_{S\ \mathcal C \mathcal C \textrm{ of } \mathcal G'} \prod_{\tilde S \in \mathcal C\mathcal C\textrm{ of } S} \tau^0[\tilde S].$$
	
\end{proof}

\begin{proposition}\label{Lem:IndPropAlg} Asymptotic traffic independence of a family of matrices is a property of the traffic space it generates, in the following sense.  Let $\mathbf A_N^{(1)} \etc \mathbf A_N^{(L)}$ be asymptotically traffic independent families of matrices. For each $\ell=1 \etc L$, let $\mathbf B_N^{(\ell)}$ be a family of matrices in the traffic space generated by $\mathbf A_N^{(\ell)}$, that is $\mathbf B_N^{(\ell)}$ is a collection of matrices of the form $g(\mathbf A_N^{(\ell)})$ for $g$ $^*$-graph polynomials in $\mathbb C\mathcal G\langle x^{(\ell)}, x^{(\ell)*}\rangle$. Then $\mathbf B_N^{(1)} \etc \mathbf B_N^{(L)}$ are asymptotically traffic independent.
\end{proposition}

In particular, if matrices $A_{N,1} \etc \tilde A_{N,L}$ are asymptotically traffic independent, then the families of matrices $(A_{N,\ell},A_{N,\ell}^2, A_{N,\ell}^t, Diag(A_{N,\ell}), Deg(A_{N,\ell}), A_{N,\ell}\circ A_{N,\ell}, \dots), \ell=1\etc L$ are asymptotically traffic independent.

 In order to prove the proposition, let us first state a property of the injective version $\tau^0$ of a map $\tau$. It is an analogue for traffics of the formula for the free cumulants whose entries are products of random variables \cite{KS00}

\begin{lemma}\label{Lem:SubsInjTrace0} Let $\tau: \mathbb C \mathcal T\lara \to \mathbb C$ be a linear form. For a $^*$-test graph $T$ in variables $\mathbf y=(y_1\etc y_n)$ and $^*$-graph monomials $g_1\etc g_n$ in variables $\mathbf x$, recall we define the $^*$-test graph $T(g_1\etc g_n)$ in variables $\mathbf x$ by replacing in $T$ the edges labeled $y_i$ by the graph $g_i$. We define $\tilde \tau : \mathbb C \mathcal T\langle \mathbf y, \mathbf y^* \rangle \to \mathbb C$ by $\tilde \tau\big[T] = \tau\big[ T(g_1\etc g_n)\big]$.

Let us fix $\tilde T = T(g_1 \etc g_n)$ and denote by $\tilde V$ the vertex set of $\tilde T$. To each vertex $v$ of $T$ is associated a vertex $f(v)$ of $\tilde T$ which gives the position of $v$ in the new graph. Note that this map $f:V \to \tilde V$ is not necessarily injective since some vertices of $T$ can be identified in $\tilde T$. For a partition $\tilde \pi \in \mathcal P(\tilde V)$, we define a partition $\tilde \pi_{|V} \in \mathcal P(V)$ via $ \tilde \pi_{|V}:= f^{-1}(\tilde \pi)$. We denote by $0_V$ the partition whose blocks are singletons. We then have
	\eq
		\tilde \tau^0\big[T\big] & =  & \sum_{ \substack{ \tilde \pi \in \mathcal P( \tilde V) \\ \mathrm{s.t.} \tilde \pi_{|V} = 0_V }} \tau^0\big[ \tilde T^{\tilde \pi} \big].
	\qe
	Note that if there is an edge
	 labeled $y_i$ which is not a self loop and such that $\Delta(g_i)=g_i$, then there is no $\tilde \pi$ in $\mathcal P(\tilde V)$ such that $\tilde \pi_{| V} =0_V$, so $\tilde \tau^0\big[ T] =0$.
\end{lemma}

This lemma is actually an observation for the trace of $^*$-test graphs in matrices. We will need this general version only later in Section \ref{Sec:TrafficSpaces}, using only the definition of the injective trace.

\begin{proof} 
 We have clearly
	\eq
		\tilde \tau\big[T\big] & :=  & \tau\big[ \tilde T \big] =  \sum_{ \tilde \pi \in \mathcal P(\tilde V)} \tau^0\big[ \tilde T^{\tilde \pi} \big]\\
			& = & \sum_{ \pi \in \mathcal P(V)}  \Bigg[ \sum_{ \substack{ \tilde \pi \in \mathcal P( \tilde V) \\ \mathrm{s.t.} \tilde \pi_{|V} = \pi }}  \tau^0\big[\tilde T^{\tilde \pi} \big] \Bigg].
	\qe
	Note that the term in the bracket depends only on $T^\pi$, and not on $T$: each $\tilde T^{\tilde \pi}$, where $\tilde \pi$ is as in the sum inside the bracket, is obtained from $T^\pi$ by constructing $(T^\pi)(g_1\etc g_n)$ (replacing each edge $y_i$ by the graph $g_i$) and then identifying some vertices of the latter graph according to a partition $\sigma$ of its vertex set, namely
		$$\big( T(g_1\etc g_n)\big)^{\tilde \pi} = \big( (T^\pi)(g_1\etc g_n)\big)^{\sigma}.$$	
	
	Hence the term in the bracket is indeed a linear form on elements of $\mathbb C \mathcal T\lara$. This implies by uniqueness of the injective trace that this term is $\tilde \tau^0[T^\pi]$, and we get the result since $\tilde \tau^0[T] = \tilde \tau^0[T^{0_V}]$. 
	
\end{proof}

\begin{proof}[Proof of Proposition \ref{Lem:IndPropAlg}] It is sufficient to prove the result for $L=2$, the general case follows by associativity and by a simple induction on the number of families of matrices. We denote $\mathbf A_N = \mathbf A_N^{(1)} \sqcup \mathbf A_N^{(2)}$ (the union of the families remembering their origin $\ell\in 1,2$). In the sequel we prove the following. Let $g$ be a $^*$-graph polynomial and denote $B = g(\mathbf A_N^{(1)})$. Then $(B, \mathbf A_N^{(1)})$ and $\mathbf A_N^{(2)}$ are asymptotically traffic independent. By induction and same reasoning for matrices of the form $B=g(\mathbf A_N^{(2)})$, we then obtain: for families $\mathbf B_N^{(1)}$ and $\mathbf B_N^{(2)}$ as in the statement of the lemma, the families $\mathbf A_N^{(1)}\cup \mathbf B_N^{(1)}$ and $\mathbf A_N^{(2)}\cup \mathbf B_N^{(2)}$ are traffic independent. Hence $\mathbf B_N^{(1)}$ and $\mathbf B_N^{(2)}$ are asymptotically traffic independent (heredity of traffic independence is immediate).

Hence in the following we introduce several times a matrix $B$ in the traffic space generated by $\mathbf A_N^{(1)}$.  For any $^*$-test graph in variables $ b, \mathbf a^{(1)}, \mathbf a^{(2)}$, we denote $  \tau\big[T\big]= \Nlim \tau_N\big[T(B, \mathbf A_N^{(1)},\mathbf A_N^{(2)})\big]$. We prove the asymptotic independence of $(B, \mathbf A_N^{(1)})$ and $\mathbf A_N^{(2)}$, namely that $  \tau^0$ satisfies Formula \ref{DefTrafInd},
		\eqa \label{Eq:ProofCompInd}
			  \tau^0\big[ T\big] = \one\big( \mathcal G \mathcal C \mathcal C(T) \ \mathrm{is \ a \ tree} \big)   \prod_{S \in \mathcal C \mathcal C(T)}  \tau^0\big[ S\big],
		\qea
where the notion of colored components is with respect to $(b,\mathbf a^{(1)})$ and $\mathbf a^{(2)}$.

{\bf 1.} Consider  a matrix $A$ of $\mathbf A_N$  and a complex number $\lambda$, and set the matrix $B=\lambda A$. The map $  \tau^0$ is anti-multilinear with respect to its edges in the following sense. For a $^*$-test graph $T$ in variables $b,  \mathbf a^{(1)}, \mathbf a^{(2)}$, denote $\tilde T$ the $^*$-test graph in the variables $\mathbf a^{(1)}, \mathbf a^{(2)}$ obtained from $T$ by replacing the label of the edges labeled $b$ (respectively $b^*$) by $a$ (respectively $a^*$). Let us set $p$ (respectively $q$) the number of edges of $T$ with label $b$ (respectively $b^*$). Then $  \tau^0\big[T\big] = \lambda^p \bar\lambda^q   \tau^0\big[\tilde T\big]$.  Moreover, by asymptotically traffic independence of $\mathbf A_N^{(1)}$ and $\mathbf A_N^{(2)}$, we have $\tau^0\big[\tilde T\big] = \one\big( \mathcal G \mathcal C \mathcal C(\tilde T) \mathrm{ \ is \ a \ tree} \big) \prod_{\tilde C\in \mathcal C \mathcal C(\tilde T)} \tau^0\big[ \tilde S \big]$, where the notion of colored components of $\tilde T$ is with respect to $ \mathbf a^{(1)}$ and $\mathbf a^{(2)}$. Furthermore, we have $\lambda^p \bar\lambda^q\prod_{\tilde S\in \mathcal C \mathcal C (\tilde T)} \tau^0\big[\tilde S\big]=\prod_{S\in \mathcal C \mathcal C (T)}   \tau^0\big[S\big]$ since the total number of edges with first label remains unchanged. Hence asymptotic traffic independence is stable by multiplication by scalars.

{\bf 2.} Let us consider now another matrix $ A' $ in the same family as $A $, and denote $B=A+A'$. Fix a $^*$-test graph $T$ in variables $b, \mathbf a^{(1)}, \mathbf a^{(2)}$ and denote by $E_1$ the set of edges of $T$ with label $b$. For any map $\gamma:E_1 \to \{a,a'\}$, we denote by $\tilde T_\gamma$ the $^*$-test graph in variables $\mathbf a^{(1)}, \mathbf a^{(2)}$ obtained from $T$ by replacing the label $b$ of all $e\in E_1$ by $\gamma(e)$. We then have $  \tau^0\big[T\big] = \sum_{\gamma:E_1\to \{a, \tilde a\}} \tau ^0\big[\tilde T_\gamma \big]$. Since both matrices $A$ and $A'$ are in a same family, then the graph of colored components of $\tilde T_\gamma$ is the same as $T$ for any $\gamma$. For a colored component $S$ of $T$ corresponds a colored component $\tilde S_\gamma$ of $\tilde T_\gamma$. With $E_{1,S}$ denoting the set of edges of $S$ with label $b$, we have  $\sum_{\gamma:E_{1,S}\to \{a,\tilde a\}}  \tau ^0\big[\tilde S_{\gamma}  \big] =   \tau ^0\big[S \big]$. Hence asymptotic traffic independence is stable by linear operations.

{\bf 3.} It remains to prove the stability under graph monomials. We first consider a matrix of the form $B=\Delta(A)$, i.e. $B$ consists in the diagonal elements of a matrix of $\mathbf A_N$. Let $T$ be a $^*$-test graph in variables $b, \mathbf a^{(1)}$ and $\mathbf a^{(2)}$. Assume there is at least one edge $(v,w)$ labelled $b$ which is not a self loop. Then, by the last remark of Lemma \ref{Lem:SubsInjTrace0}, $\tilde \tau^0 [T ]$ vanishes, and as well $ \tau^0 [S ]=0$ for the colored component of $T$ containing $(v,w)$, so Formula \eqref{Eq:ProofCompInd} is valid. On the other hand, if all the edges labeled $b$ in $T$ are loops, then with $\tilde T$ is obtained from $T$ by replacing $b$ by $a$, one has $\tilde \tau\big[T\big] = \tau\big[\tilde T\big]$, and so the same equality is true for non-injective trace. Hence the formula.

{\bf 4.} At last we consider a matrix $B=g(\mathbf A_N^{(1)})$ where $g$ is a $^*$-graph monomial such that the input and the output are distinct ($g\neq \Delta(g)$). By Lemma \ref{Lem:SubsInjTrace0}, with $\tilde T$ obtained from $T$ by replacing edges labeled $b$ by the graph monomial $g$, we have 
	$   \tau^0\big[ T\big] = \sum_{ \substack{ \tilde \pi \in \mathcal P(\tilde V) \\ s.t. \tilde \pi_{|V} = 0_{V}}} \tau^0 \big[ \tilde T^{\tilde \pi}\big].$
 We know by asymptotic traffic independence of $\mathbf A_N^{(1)}$ and $\mathbf A_N^{(2)}$ that
		$\tau^0 \big[ \tilde T^{\tilde \pi} \big] = \one\big( \mathcal G \mathcal C \mathcal C( \tilde T^{\tilde \pi}) \ \mathrm{is \ a \ tree} \big)   \prod_{S \in \mathcal C \mathcal C( \tilde T^{\tilde \pi})}  \tau^0\big[ S\big],$
where the notion of colored components is with respect to $\mathbf a^{(1)}$ and $\mathbf a^{(2)}$. Let $\tilde \pi \in \mathcal P(\tilde V)$ such that $\tilde \pi_{|V} = 0_V$. In other words, the vertices of $V$, which are not identified in $\tilde V$, are not identified in $\tilde V^{\tilde \pi}$ neither. We see from now $V$ as a subset of $\tilde V^{\tilde \pi}$.

We claim that $\mathcal G \mathcal C \mathcal C(\tilde T^{\tilde \pi})$ is a tree only if $\mathcal G \mathcal C \mathcal C(  T)$ is a tree. To prove this claim, in the graph $\tilde T$ (where the edges labeled $b$ are replaced by $g$) we chose an enumeration of the vertices $\{v_1 \etc v_Q\}$ of $\tilde V \setminus V$ (there are copies of the vertices of $g$ that are not the input nor the output). For each $q=0\etc Q$ we denote $V_q = V \cup \{v_1 \etc v_q\}$ ($V_0=V$). Moreover, denoting the partition $\tilde \pi = \{ \tilde B_1, \dots , \tilde B_{\tilde F}\}$, we set the partition $\tilde \pi(q)= \big \{ \tilde B_1\cap V_q \etc  \tilde B_{\tilde F}\cap V_q \big  \}$. We have $\tilde \pi(Q) = \tilde \pi$ by definition. Moreover, we have $\mathcal G \mathcal C \mathcal C( \tilde T^{\tilde \pi(0)}) = \mathcal G \mathcal C \mathcal C(T)$ since $\tilde \pi(0)$ is made of singletons (we can replace $g$ by every connected graph labeled in the same family without changing $\mathcal G\mathcal C\mathcal C(\tilde T)$). Denote by $|\mathcal V_q|$ and $|\mathcal E_q|$ the number of vertices and edges respectively of $\mathcal G \mathcal C \mathcal C( \tilde T^{\tilde \pi(q)})$. By Lemma \ref{lem:EdgesVertices}, we will get the claim if we prove that the sequence $|\mathcal V_q|-|\mathcal E_q|$ is non increasing. 

Passing from $\tilde T^{\tilde \pi(q)}$ to $\tilde T^{\tilde \pi(q+1)}$ means choosing a possible identification of the vertex $v_{q+1}$ with the vertices of $\tilde T^{\tilde \pi(q)}$ different from $v_{q+2} \etc v_Q$. If $v_{q+1}$ belongs to a colored component $S \in \mathcal C \mathcal C( \tilde T^{\tilde \pi(q)})$ and is identified with a vertex $w$ of $\tilde T^{\tilde \pi(q)}$ which is not in $S$, we see that we always have $|\mathcal V_{q+1}|-|\mathcal E_{q+1}| = |\mathcal V_q|-|\mathcal E_q| - 1$ by the following enumeration:
\begin{enumerate}
	\item if $w$ is in a colored component of $\tilde T^{\tilde \pi(q)}$ of the same color as $S$ then we merge the two colored components in $\tilde T^{\tilde \pi(q+1)}$, and so $|\mathcal V_{q+1}| = |\mathcal V_{q}| - 1$, $|\mathcal E_{q+1}|= |\mathcal E_{q}|$;
	\item otherwise, $w$ is not in a same colored component as $v_{q+1}$ in $\tilde T^{\tilde \pi(q+1)}$. Identified in $\tilde T^{\tilde \pi(q)}$, the fusion of the vertices creates a new connector vertex with degree two between $S$ and the colored component of $w$, and so $|\mathcal V_{q+1}| = |\mathcal V_{q}| +1$, $|\mathcal E_{q+1}|= |\mathcal E_{q}|+2$.
\end{enumerate}
Assume now that $v_{q+1}$ is identified with a vertex $w$ of $\tilde T^{\tilde \pi(q)}$ in $S$.
This does not modify the graph of colored components. We then get the claim that if $\mathcal G \mathcal C \mathcal C(\tilde T^{\tilde \pi})$ is a tree, then $\mathcal G \mathcal C \mathcal C(  T)$ is a tree. 

We also see that the $\tilde \pi$ (such that $\tilde \pi_V = 0_V$) for which $\mathcal G\mathcal C\mathcal C(\tilde T^{\tilde \pi})$ is a tree are given by the choice for each colored component $S$ of $T$ of a partition $\tilde \pi_S$ (with same restriction on $V$). Denote by $V_S$ the vertex set of $S$, by $\tilde V_S$ its image in $\tilde T$. Moreover, for a partition $\pi$ of the vertex set $V$ of $T$, denote by $\tilde \pi_{S|V_S}$ its restriction to $S$. 
We get the expected result as
	\eq
	\tau^0\big[ T(\mathbf b)\big] & = &   \one(\mathcal G \mathcal C \mathcal C(T) \mathrm{ \ is \ a \ tree}) \prod_{S\in \mathcal C \mathcal C(T^\pi)}  \sum_{ \substack{ \tilde \pi_S \in \mathcal P(\tilde V_S) \\ s.t. \tilde \pi_{S|V_S} = 0_{V_S}}} \tau^0\big[  \tilde S^{\tilde \pi_S}\big] \\
	& = &  \one(\mathcal G \mathcal C \mathcal C(T) \mathrm{ \ is \ a \ tree}) \prod_{S\in \mathcal C \mathcal C(T^\pi)} \tau^0\big[  \tilde S\big]. 
	\qe

\end{proof}

\section{A criterion of lack of asymptotic free independence}\label{Sec:CritLackInp}

\begin{figure}
     \includegraphics{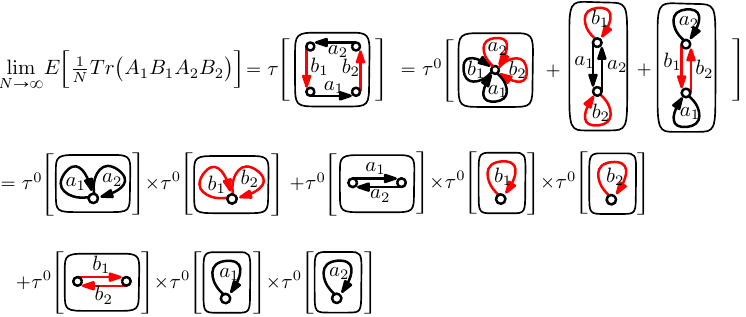}
\caption{Computation of \eqref{Eq:ComputationMoment4}}
    \label{fig:Proof1}
  \end{figure}

Let $\mathbf A_N =(A_1,A_2)$ and $\mathbf B_N =(B_1,B_2)$ be  two asymptotically traffic independent families of matrices. Let us first manipulate the definitions to understand how to compute 
		$$\Nlim \esp\big[ \frac 1 N \Tr \, A_1 B_1 A_2 B_2 \big].$$
	We denote for all $^*$-test graph $T$ and all $^*$-polynomials $P,Q$,
		\eq
			\tau\big[T(\mathbf a, \mathbf b)\big] & = & \Nlim\esp\big[ \frac 1 N \Tr  \, T(\mathbf A_N, \mathbf B_N)\big],\\
			 \Phi\big[P(\mathbf a, \mathbf b)\big] & = & \Nlim \esp\big[ \frac 1 N \Tr  \, P(\mathbf A_N, \mathbf B_N) \big], \\
			\Phi\big[P(\mathbf a, \mathbf b) \circ Q(\mathbf a, \mathbf b) \big]&  = &\Nlim \esp\Big[ \frac 1 N \Tr  \, \big[P(\mathbf A_N, \mathbf B_N) \circ Q(\mathbf A_N, \mathbf B_N)\big] \Big],
		\qe
where we recall that $\circ$ denotes the entry-wise product of matrices. Firstly, according to Definition \ref{Def:CombiDistrTraf}, we introduce the $^*$-test graph $T$ consisting in a simple cycle with edges labeled $b_2, a_2,b_1,a_1$ along it, so that we can write $\Nlim \esp\big[ \frac 1 N \Tr \, A_1 B_1 A_2 B_2 \big] = \tau\big[T(\mathbf a, \mathbf b)\big] = \sum_{\pi\in \mathcal P(V)} \tau^0\big[T(\mathbf a, \mathbf b)\big] $, where $V$ is the set of the four vertices of $T$. By Definition \ref{def:FreeProdGraphs} the asymptotic traffic independence of $\mathbf A_N$ and $\mathbf B_N$ implies that there are three partitions $\pi\in \mathcal P(V)$ that contribute on the limit (the computation is illustrated in Figure \eqref{fig:Proof1}). One has $\tau\big[T(\mathbf a, \mathbf b)\big]  = \sum_{i=1}^3 \tau^0\big[T_i(\mathbf a, \mathbf b)\big] $ for the following $^*$-test graph $T_1,T_2,T_3$.
\begin{itemize}
	\item The $^*$-test graph  $T_1$ with one vertex and four loops, labelled $a_1,a_2,b_1,b_2$ respectively. By the factorization property, one has 
	$$\tau^0 [T_1(\mathbf a, \mathbf b) ] =\tau^0  [T_1'(\mathbf a) ]  \times \tau^0   [T_1'(  \mathbf b) ],$$
 where $T'_1(x,y)$ is the $^*$-test graph with one vertex and two loops, labeled $x$ and $y$ respectively. Since $T'_1$ has a single vertex, one has $\tau^0 [T'_1 (\mathbf a)] = \tau [T'_1(\mathbf a) ]$. By the second case of Example \ref{Ex:DistrFromTrafficsBis}, we finally get $\tau [T'_1(\mathbf a)] = \Phi  [a_1 \circ a_2]$, as well as the similar expression for $\mathbf b$.
	\item The $^*$-test graph $T_2$ with two vertices, one loop on each vertex, and one edge between the two vertices in each direction. Labels are $a_1$ and $a_2$ on the loops, $b_1$ and $b_2$ on the simple edges, in such a way one can read the word $b_2,a_2, b_1, a_1$ by following a closed path. The factorization property yields 
		$$\tau^0 [T_2(\mathbf a, \mathbf b) ] = \tau^0 [T_0(  a_1)  ] \times  \tau^0 [T_0(  a_2) ] \times  \tau^0 [T'_2(\mathbf b) ],$$ where $T_0(a)$ is the $^*$-test graph with one vertex and one loop labeled $a$ and $T'_2(\mathbf b)$ the $^*$-test graph with two vertices and a pair twin edges labeled $b_1,b_2$ and with opposite orientation. Since $T_0$ has a single vertex, one has $\tau^0 [T_0(  a_i)  ]  =\tau [T_0(  a_i)  ] = \Phi(a_i)$. Moreover, the relation $\tau [T'_2(\mathbf b) ] = \tau^0 [T'_2(\mathbf b) ] + \tau^0 [T'_1(\mathbf b)  ]$ yields $ \tau^0 [T'_2(\mathbf b) ] = \Phi(b_1b_2) -  \Phi (b_1\circ b_2)$.
	\item The $^*$-test graph  $T_3$ obtained similarly with the roles of $\mathbf a$ and $\mathbf b$ interchanged.
\end{itemize}	
This gives by summing the three contributions
	\eqa\nonumber
		\lefteqn{ \Nlim \esp\big[\frac 1 N \Tr [ A_1B_1A_2B_2] \big] =  \Phi  ( a_1\circ a_2  )   \times \Phi ( b_1\circ b_2  )  } \\\label{Eq:ComputationMoment4}
	& & \ \ + \Big[ \Phi ( a_1^2 ) - \Phi ( a_1\circ a_2  )  \Big] \times \Phi ( b_1\big) \Phi (  b_2  ) +   \Big[ \Phi ( b_1^2 ) - \Phi ( b_1\circ b_2  )  \Big] \times \Phi ( a_1\big) \Phi (  a_2 \big).
	\qea

Using this formula for $ \overset \circ a :=  a_i -\Phi(a_i)$ and $ \overset \circ b_i = b_i -\Phi (b_i)$ yields $\Phi ( \overset \circ a_1 \overset \circ b_1 \overset \circ a_2 \overset \circ b_2) = \Phi ( \overset \circ a_1\circ \overset \circ a_2  )   \times \Phi  ( \overset \circ b_1\circ \overset \circ b_2  ).$ If this quantity is not zero then $\mathbf A_N$ and $\mathbf B_N$ are not asymptotically free independent. Writing $\Phi ( \overset \circ a_1\circ \overset \circ a_2  ) = \Phi(a_1 \circ a_2) - \Phi(a_1) \Phi(a_2)$ yields the following criterion which is useful for applications.

\begin{proposition}[Criterion of lack of freeness]\label{Prop:Criterion} Let $\mathbf A_N$ and $\mathbf B_N$ be two asymptotically traffic independent families. Denote $\Phi_N = \esp\big[   \frac 1 N \Tr \, \cdot \, \big]$. Define for any self-adjoint $^*$-polynomials $P,Q$,
	\eqa\label{Def:MFK}
		\lefteqn{\mathfrak K\big(P,Q,(\mathbf A_N)_{N\geq 1}\big) :=}\\
		 & \Nlim & \bigg(  \ \Phi_N\big(  P(\mathbf A_N)\circ   Q(\mathbf A_N) \big) - \Phi_N \big(P(\mathbf A_N)\big) \times \Phi_N \big(Q(\mathbf A_N)\big) \ \bigg)\nonumber
	\qea
where $\circ$ denotes the entry-wise product of matrices, and define $\mathfrak K\big(P,Q,(\mathbf B_N)_{N\geq 1}\big)$ similarly.
Assume that $\mathfrak K\big(P,Q,(\mathbf A_N)_{N\geq 1}\big)$ and $\mathfrak K\big(\tilde P,\tilde Q,(\mathbf B_N)_{N\geq 1}\big)$ are nonzero for some $P,Q, \tilde P, \tilde Q$. Then $\mathbf A_N$ and $\mathbf B_N$ are not asymptotically free independent.
\end{proposition}

\begin{example}~ \begin{itemize}
	\item {\bf Diagonal matrices:} Let $\mathbf D_N$ be a family of diagonal matrices. Then $\mathfrak K\big(P,Q,(\mathbf D_N)_{N\geq 1}\big)$ is the covariance  
		$$ \Nlim  \ \Phi_N\big(P(\mathbf D_N) \times Q(\mathbf D_N) \big) - \Phi_N\big(P(\mathbf D_N) \big) \times  \Phi_N \big(Q(\mathbf D_N) \big)$$
 of the limiting $^*$-distribution. So, by the Cauchy Schwarz inequality, it is zero for all $P$ and $Q$ if and only if for all matrices $D_N$ of $\mathbf D_N$ one has $\Nlim( \Phi_N(D_N^2)-\Phi_N(D_N)^2)=0$.
	\item {\bf Matrices with independent entries:} Let $d>0$ and let $A_N$ be a Hermitian matrix whose sub-diagonal entries are independent and distributed according to the Bernouilli distribution with parameter $\frac d N$. Then, by independence of the entries of $A_N$,
	\eq
		\lefteqn{\Phi_N\big( A_N^2 \circ A_N^2 \big) - \Phi_N\big( A_N^2\big) ^2}\\
		 & = & \esp\big[ \frac 1 N \sum_{i} \big( \sum_j A_N(i,j)^2\big)^2 \big] - 
		\esp\big[ \frac 1 N \sum_{i,j} A_N(i,j)^2\big]^2 \\
		& = & \esp\big[ \frac 1 N \sum_{i,j} A_N(i,j)^4\big] +o(1) = \frac 1 N \sum_{i,j} \frac d N = d >0.
	\qe
	This particular model is studied in a more general framework in \cite{MAL122} for which application of the convergence in traffic distribution is specified.
	\item {\bf Adjacency matrices of non regular graphs:} Let $G_N$ be a random graph on the set of vertices $[N]$. Assume that $G_N$ is a directed simple graph (with no loops nor multiple edges). The adjacency matrix $A_N$ of $G_N$ is the random matrix with entries in $\{0,1\}$, such that $A_N(i,j)=1$ if and only if there is an edge from $j$ to $i$, and $A_N(i,i)=0$ for all $i$.
	
 For a vertex $v$ uniformly chosen at random, let us denote $d_N = \sum_{j=1}^N A_N(v,j)$ (the number of whose output is $v$). Then 
 	$$\Phi_N(A_NA_N^* ) = \esp\Big[ \frac 1 N \sum_{i,j=1}^N A_N(i,j)^2  \Big] = \esp[ d_N]$$
 since $ A_N(i,j)\in \{0,1\}$ and
	$$\Phi_N\big( A_NA_N^* \circ A_NA_N^* \big) = \esp\Big[ \frac 1 N \sum_{i=1}^N \big( \sum_{j=1}^N A_N(i,j)^2 \big)^2 \Big] =  \esp\big[ d_N^2\big].$$
	Hence $\mathfrak K(xx^*,xx^*, (A_N)_{N\geq 1})$ is the limit of the variance of $d_N$. So if the degree of a vertex of $G_N$ uniformly chosen does not converge to a constant, then the adjacency matrix $A_N$ of $G_N$ satisfies Proposition \ref{Prop:Criterion}.
	\end{itemize}
\end{example}

\section{Proof of the main theorem}\label{Sec:ProofMainTh}

\subsection{Presentation of the method}

In order to compute the limiting traffic distribution of several random matrices (and to prove Theorem \ref{MainTh}) we often use the injective trace $\tau_N^0$ defined in Section \ref{Sec:InjTrace} since it is directly related to the moments of the entries of the matrices.

 \begin{lemma}\label{Lem:Tool} Let $\mathbf A_N = (A_j)_{j\in J}$ be a family of random matrices. Let $T=(V,E, \gamma, \varepsilon)$ be a $^*$-test graph and let $\phi_N$ is be a random injective map $V\to [N]$, uniformly distributed, independent of $A_N$. Denote 
 	$$ \delta^0_N\big[T( \mathbf A_N) \big] := \esp \bigg [    \prod_{ e=(v,w) \in E } A_{\gamma(e)}^{\varepsilon(e)}\big( \phi_N(w), \phi_N(v) \big) \bigg].$$ 
 	\begin{enumerate} 
	 \item If $\mathbf A_N$ is a permutation invariant family of random matrix, then for any $T=(V,E, \gamma, \varepsilon)$, one has
	 	$$ \delta^0_N\big[T( \mathbf A_N) \big] =\esp \bigg [    \prod_{ e=(v,w) \in E } A_{\gamma(e)}^{\varepsilon(e)}\big( \phi(w), \phi(v) \big) \bigg]$$
	 for any injective map $\phi:V\to [N]$. Hence the function $\delta^0_N$ computes joint moments in the entries of the matrices of $\mathbf A_N$.

	\item For any $^*$-test graph $T$ with vertex set $V$ and any family $\mathbf A_N$ of matrices, one has
	\eqa\label{Eq:LienTr0Delt0}	
		\tau_N^0\big[ T( \mathbf A_N)\big] = \frac 1 N \frac {N!} {(N-|V|)!} \delta_N^0\big[ T( \mathbf A_N)\big].
	\qea
\end{enumerate}
 \end{lemma}
 
 \begin{proof}1. Let $\phi: V\to [N]$ be an injective map. Let $\sigma_N$ a random permutation of $[N]$ uniformly distributed and independent of $\mathbf A_N$. By the proof of Lemma \ref{Lem:PermInvDesTraces}, since $\mathbf A_N$ is permutation invariant one has
 	$$ \esp\Big[  \prod_{(v,w) \in E} A_{\gamma(e)}^{\varepsilon(e)}\big( \phi(w) , \phi(v) \big) \Big] = \esp\Big[  \prod_{(v,w) \in E} A_{\gamma(e)}^{\varepsilon(e)}\big( \sigma_N \circ \phi(w) , \sigma_N \circ\phi(v) \big) \Big].$$
 But $\sigma_N \circ\phi$ is distributed as a uniform injective map $\phi_N:V\to [N]$ independent of $\mathbf A_N$, hence the result.
 
 2. We get the result by the following interpretation of the sum as an expectation
 \eq
\lefteqn{\tau_N^0\big[ T(\mathbf A_N) \big] }\\& = &
		\frac 1 N \times \frac {\frac {N!}{ (N-|V|)!} } {\textrm{Card } \Big\{ \substack{\phi:V \rightarrow \{1\etc N\} \\ \mathrm{injective}}  \Big\} }  \sum_{\substack{\phi:V \rightarrow \{1\etc N\} \\ \mathrm{injective}}} \prod_{e=(v,w) \in E} A_{\gamma(e)}^{\varepsilon(e)}  \big( \phi(w), \phi(v)  \big)   \\
		& = &\frac 1 N \frac {N!}{ (N-|V|)!} \esp\bigg[ \prod_{e =(v,w)\in E} A_{\gamma(e)}^{\varepsilon(e)}   \big( \phi_N(w), \phi_N(v)   \big)\bigg].
\qe
\end{proof}
 
Let see an application of this lemma, which is a sort of asymptotic independence theorem for the entry-wise product of matrices. We will use it later as illustrations of traffic independence. 
 
\begin{corollary}\label{Cor:VarProfile} Let $\mathbf A_N = (A_j)_{j\in J}$ and $\mathbf B_N = (B_j)_{j\in J}$ be two independent families of random matrices. Assume the following hypotheses. 
\begin{enumerate}
	\item The family $\mathbf A_N$ converges in traffic distribution, that is for any $^*$-test graph $T$ the limit $\tau^0[T]:=\Nlim\tau_N^0\big[T( \mathbf A_N) \big]= N^{|V|-1}\big( 1 + O(\frac 1 N) \big) \delta_N^0\big[T( \mathbf A_N) \big]$ exists.
	\item  For any $^*$-test graph $T$, the limit $\delta^0[T]:=\Nlim\delta_N^0\big[T( \mathbf B_N) \big]$ exists.

	\item One of the two families $\mathbf A_N$, $\mathbf B_N$ is permutation invariant.
\end{enumerate}
Denote $\mathbf C_N = (A_j \circ B_j)_{j\in J}$, where $\circ$ stands for the entry-wise product of matrices. Then $\mathbf C_N$ converges in traffic distribution. Moreover, for any $^*$-test graph $T$, one has
	$$\Nlim \tau_N^0\big[T(\mathbf C_N) \big] = \tau^0[T] \,  \delta^0[T] .$$
\end{corollary}

In the proof of the corollary we see an argument that appears also in the proof of Theorem \ref{MainTh}.

\begin{example} Let $A_N$ be a random matrix converging in traffic distribution. Let $B_N=(\omega_{i,j})_{i,j=1\etc N}$ be a matrix with i.i.d. entries, possibly up to the Hermitian symmetry. Assume that the distribution of the $\omega_{i,j}$'s is independent of $N$, invariant by complex conjugate, and that $\esp[|\omega_{i,j}|^K]<\infty$ for each $K\geq 1$. Then $A_N \circ B_N$ converges in traffic distribution. See Section \ref{Sec:ExamplesMatrices} for applications.
\end{example} 
\begin{proof} Let $T=(V,E,\gamma,\varepsilon)$ be a $^*$-test graph and $\phi_N$ a uniform injective map $V\to[N]$. Then $\tau_N^0\big[ T(\mathbf C_N) \big] = N^{|V|-1} \big( 1 + o\big( \frac 1 N \big) \big) \times \delta_N^0\big[T(\mathbf C_N)\big]$. Moreover, one has
	$$	\delta_N^0\big[ T(\mathbf C_N)\big] =  \esp \bigg [    \prod_{ e=(v,w) \in E } A_{\gamma(e)}^{\varepsilon(e)}\big( \phi_N(w), \phi_N(v) \big) \times  B_{\gamma(e)}^{\varepsilon(e)}\big( \phi_N(w), \phi_N(v) \big) \bigg].$$
	Assume that $\mathbf B_N$ is permutation invariant. Let $V$ be a random permutation matrix, associated to a uniform permutation $\sigma_N$ independent of $\mathbf A_N$. Then 
	\eq
		\lefteqn{	\delta_N^0\big[ T(\mathbf C_N)\big] }\\
		& = &  \esp \bigg [    \prod_{ e=(v,w) \in E } A_{\gamma(e)}^{\varepsilon(e)}\big( \phi_N(w), \phi_N(v) \big) \times  B_{\gamma(e)}^{\varepsilon(e)}\big( \sigma_N\circ \phi_N(w),  \sigma_N\circ \phi_N(v) \big) \bigg].
	\qe

 Since $\sigma_N$ is uniform and independent of $(\mathbf A_N, \mathbf B_N, \phi_N )$, then $(\mathbf A_N, \mathbf B_N, \phi_N , \sigma_N\circ \phi_N )$ has the same law as $(\mathbf A_N, \mathbf B_N, \phi_N ,  \phi_N' )$ where $\phi_N '$ is a uniform injection independent of $(\mathbf A_N, \mathbf B_N, \phi_N )$.
	Hence, by independence of $(\mathbf A_N, \phi_N)$ and $(\mathbf B_N, \phi_N')$ we obtain
	\eq
			\delta_N^0\big[ T(\mathbf C_N)\big] 
			&  = & \esp \bigg [    \prod_{ e=(v,w) \in E } A_{\gamma(e)}^{\varepsilon(e)}\big( \phi_N(w), \phi_N(v) \big) \bigg] \times \esp \bigg[  B_{\gamma(e)}^{\varepsilon(e)}\big(   \phi_N'(w),   \phi_N'(v) \big) \bigg] \\
			& = & \delta_N^0\big[ T(\mathbf A_N)\big] \times \delta_N^0\big[ T(\mathbf B_N)\big].
	\qe
	Note that the role of $\mathbf A_N$ and $\mathbf B_N$ can be interchanged so there is no limitation in assuming $\mathbf B_N$ permutation invariant. We obtain 
	\eq
		\delta_N^0\big[ T(\mathbf C_N)\big]  & = & N^{|V|-1}\Big( 1 + o\big( \frac 1 N \big) \Big)	\delta_N^0\big[ T(\mathbf A_N)\big] 	\delta_N^0\big[ T(\mathbf B_N)\big]\\
		 & = & \tau_N^0\big[T(\mathbf A_N)\big] 	\delta_N^0\big[ T(\mathbf B_N)\big] +o(1),
	\qe
	as expected.
\end{proof}
 
An important technical aspect of Lemma \ref{Lem:Tool} is that it tells how we have to renormalized the joint moments the entries of the matrices $\mathbf A_N$ in order to compute its limiting distribution, and relies directly this normalization to a topological constant of $T$. Indeed, $\tau_N^0\big[T(\, \cdot\, )\big]$ is multiple of $\delta_N^0\big[T(\, \cdot\, )\big]$ by a constant equivalent to $N^{V-1}$. In practice, since the map $\delta_N^0\big[T(\, \cdot\, )\big]$ is multilinear with respect to the edges of $T$, for several matrix models $\mathbf A_N$ it is possible to write $\tau_N^0\big[T(\mathbf A_N)\big] = N^{|\mathcal V| -1 + |\mathcal E|}\times \eta_N $ where $\mathcal G=(\mathcal V, \mathcal E)$ is a graph depending on $T$ well chosen and $\eta_N$ is bounded. Thank to Lemma \ref{lem:EdgesVertices}, we can classify the graphs $T$ that contribute in the limit. 

\subsection{Proof}
 We prove the theorem for two families, denoted in short $\mathbf A_1=\mathbf A_N^{(1)} $ and $\mathbf A_2 =\mathbf A_N^{(2)}$, assuming $\mathbf A_2$ permutation invariant. The general case is obtained by recurrence on the number of families using that the unions of the families also satisfies the assumptions of the theorem.

\subsubsection{Two lemmas}

In this proof, we call $^*$-graph the object formally defined as a $^*$-test graph with the connectedness condition on the graph omitted (a $^*$-graph is a finite disjoint union of $^*$-test graphs). The trace and invective trace of a $^*$-graph in matrices are defined as for $^*$-test graph. If the connected components of $T$ are the $^*$-test graphs $T^{(1)} \etc T^{(K)},$ one sees immediately that $\Tr \big[ T(\mathbf A_N)\big] = \prod_{k=1}^K  \Tr\big[ T^{(k)}(\mathbf A_N)\big]$.

\begin{lemma}[Splitting the contribution due to $\mathbf A_1$ and $\mathbf A_2$]~\label{Prop:PermutationInvariance}
Let $T$ be a $^*$-graph in the variables $\mathbf x_1$ and $\mathbf x_2$. For $i=1,2$, we denote by $T_i =( V_i, E_i,\gamma_i, \varepsilon_i)$ the $^*$-graph obtained from $T$ by considering only the edges with a label in $\mathbf x_{i}$ and by deleting the vertices that are not attached to any edge after this process. Hence $T_i$ is the union of the colored component labeled $\mathbf x_i$ of $T$. We have
\begin{eqnarray}
	 \lefteqn{\tau_N^0\big[ T(\mathbf A_1  , \mathbf A_2 )\big]}\nonumber\\
	 &  = & \label{Eq:LemSplit}
	  N \times  \frac{(N-|V_1| )!  \,  (N-|V_2| )!} { (N-|V|)! \, N!}  \tau_N^0 \big[ T_1 (\mathbf A_1 ) \big] \, \times \, \tau_N^0 \big[ T_2 (\mathbf A_2 ) \big].
\end{eqnarray}
\end{lemma}

\begin{proof}
Denote $\mathbf A_N = \mathbf A_1 \cup \mathbf A_2=(A_j)_{j\in J}$, as in Lemma \ref{Lem:Tool} we denote for any $^*$-graph $T=(V,E, \gamma, \varepsilon)$, 
	$$ \delta^0_N\big[T( \mathbf A_N) \big] = \esp \bigg [    \prod_{ e=(v,w) \in E } A_{\gamma(e)}^{\varepsilon(e)}\big( \phi_N(w), \phi_N(v) \big) \bigg],$$
where $\phi_N$ is a uniform injective map $V\to [N]$ independent of $A_N$, and we have
		\eqa	\label{Eq:TraceDensity}	 
			   \tau_N^0 \big[ T(\mathbf A_N) \big] & = & \frac 1 N \times \frac{ N!}{(N-|V|)!}  \times \delta^0_N\big[T(\mathbf A_N) \big],
		\end{eqnarray}
still valid when $T$ is not connected. Moreover, one has
		\eq
		 \lefteqn{\delta^0_N\big[T(\mathbf A_N) \big] }\\
		 &=&  \esp\bigg[ \prod_{e=(v,w)  \in E_{1}} A_{\gamma(e)}^{\varepsilon(e)} \big( \phi_N(w), \phi_N(v) \big)  \times \prod_{e=(v,w) \in E_{2}} A_{\gamma(e)}^{\varepsilon(e)} \big( \phi_N(w), \phi_N(v) \big)\bigg].
		\qe

Let $W$ be a random permutation matrix, associated to a uniform permutation $\sigma_N$ independent of $\mathbf A_N$. The family $\mathbf A_2$ has the same distribution as $W \mathbf A_2 W^*$, and denoting $\phi_N(e) = \big(\phi_N(w), \phi_N(v)\big)$ for any edge $e=(v,w)$, then for $e\in E_2$ one has $(WAW^*)_{\gamma(e)}^{\varepsilon(e)} \big( \phi_N(e) \big) = A_{\gamma(e)}^{\varepsilon(e)} \big( \sigma_N\circ\phi_N(e) \big)$. Since $\sigma_N$ is uniform and independent of $(\mathbf A_N, \phi_N)$, the triplet $(\mathbf A_N, \phi_N, \sigma_N\circ \phi_N)$ has the same law as $(\mathbf A_N, \phi_N,  \phi_N')$ where $\phi_N'$ is a uniform injection independent of $(\mathbf A_N, \phi_N)$. Hence we get
	\eq
		\lefteqn{ \tau_N^0\big[ T(\mathbf A_1 , \mathbf A_2)\big] }\\ 
	& = &
		 \frac 1 N  \times 
		 \frac{N!}{(N-|V|)!}  
		  \times  \esp\bigg[ \prod_{e \in E_{1}} A_{\gamma(e)}^{\varepsilon(e)} \big( \phi_N(e) \big)\bigg] \times
		   \esp\bigg[ \prod_{e \in E_{2}} A_{\gamma(e)}^{\varepsilon(e)} \big(\phi_N'(e) \big)\bigg]\\
		   & =& \frac 1 N \frac{ N!}{(N-|V|)!} \delta_N^0\big[T_1(\mathbf A_1)\big]\delta_N^0\big[T_2(\mathbf A_2)\big].
	\qe
	
	Using again \eqref{Eq:TraceDensity}	for the graphs $T_1$ and $T_2$ yields \eqref{Eq:LemSplit}. \end{proof}

\begin{lemma}[Decomposition of components]\label{Lem:IndComp}For any $j=1,2$ and any finite $^*$-graph $T=(V,E)$ in the variables $\mathbf x_j$ whose connected components are denoted by $T_1 \etc T_n$ one has
	$$ \esp\bigg[ \frac 1 {N^n} \Tr^0\big[T(\mathbf A_j)\big] \bigg] -   \esp\Big[ \prod_{i=1}^n \frac 1 N \Tr^0  \big[ T_{i}(\mathbf A_j) \big]  \Big] = O\big(\frac 1 N\big).$$
	
\end{lemma}

\begin{proof}
By the relation between the injective and the standard one,
\begin{equation} \label{eq:ProofAsympTh:InjNonInj}
	\Tr ^0\big[T(\mathbf A_j)\big] = \sum_{\pi \in \mathcal P(V)} \mu_{V}(\pi) \ \Tr\big[T^\pi (\mathbf A_j) \big].
\end{equation}

\noindent  Denote by $K_\pi$ the number of components of $T^\pi$. By the factorization property, we have that $\esp\big[ \frac 1 {N^{K_\pi}} \Tr \big[T^\pi (\mathbf A_j) \big]$ converges and so is bounded. Denote by $\Lambda$ the set of partitions $\pi$ of $\mathcal P(V)$ such that two vertices of different components never belong to a same block of $\pi$. For $\pi \in \Lambda$ one can decompose $\pi$ into partitions $\pi_i$ of $V_i$, $i=1\etc n$ and we denote $\pi = \pi_1 \sqcup \dots \sqcup \pi_n$. We then obtain, since $K_\pi=n$ if and only if $\pi \in \Lambda$,
	\eq
			\lefteqn{\esp\Big[ \frac 1 {N^n} \Tr^0 \big[T(\mathbf A_j)\big] \Big] }\\
			& = & \sum_{\substack{ \pi_i\in \mathcal P(V_i) \\ i =1\etc n}}  \mu_{V}( \pi_1 \sqcup \dots \sqcup \pi_n) \ \esp\Big[ \prod_{i=1}^n \frac 1 N \Tr\big[T^{\pi_i} (\mathbf A_j) \big] \Big] +O\big( \frac 1 N\big).
	\qe

By Lemma \ref{Lem:FormuleTraceTraceInj} and the formula for $\mu_{V}$, one has $\mu_{V}( \pi_1 \sqcup \dots \sqcup \pi_n) = \mu_{V}( \pi_1) \times \dots \times \mu_{V}( \pi_n)$. Hence the result.
\end{proof}
 
\subsubsection{Proof of the asymptotic traffic independence}

By Lemmas \ref{Prop:PermutationInvariance} and \ref{Lem:IndComp}, for any $^*$-test graph $T$ in the variables $\mathbf x_1$ and $\mathbf x_2$, one has
\begin{eqnarray*}
	 \tau_N^0\big[ T(\mathbf A_1  , \mathbf A_2 )\big]
	& = & \underbrace{N \times \frac{(N-|V_1| )! \, (N-|V_2| )!} { (N-|V|)! \ N!} \times N^{K_1-1} \times N^{K_2-1}}_{\Gamma_N} \\
	& & \ \   \times \bigg( \prod_{i=1}^2 \esp\Big[ \prod_{k=1}^{K_i} \frac 1 N \Tr^0  \big[ T_{i,k}(\mathbf A_i) \big]\Big]  + O\big(\frac1 N\big) \bigg),
\end{eqnarray*}

\noindent where the $T_{i,k}$'s are the connected components of $T$ that are labelled by $\mathbf x_i$, for $i=1,2$ and $k=1\etc K_i$,  $|V_i|$ is the number of vertices of $T$ attached to some edges labelled in $\mathbf x_i$, for $i=1,2$, and $|V|$ is the number of vertices of $T$. Remark that  $\Gamma_N = N^{\eta}\big( 1 + O(\frac 1 N) \big) $ where $\eta = K_1 + K_2 +|V| -|V_1| - |V_2| -1$. Note that by the factorization property and Lemma \ref{Lem:EqDeccorProp}, one has
	$$\esp\Big[ \prod_{k=1}^{K_i} \frac 1 N \Tr^0  \big[ T_{i,k}(\mathbf A_i) \big]\Big] = \prod_{k=1}^{K_i} \Nlim \tau_N^0\big[ T_{i,k}(\mathbf A_i) \big] +o(1).$$
So its remains to prove that $\Gamma_N\limN 1$ if the graph $\mathcal G \mathcal C \mathcal C(T)$ of colored components of $T$ with respect to $\mathbf x_1, \mathbf x_2$ is a tree, and $\Gamma_N\limN 0$ otherwise.

  Recall that the set of vertices of $\mathcal G \mathcal C \mathcal C(T)$ is the disjoint union of the set of colored components $T_{i,k}$, $i=1,2$, $k=1\etc K_i$, and the set $\delta V$ of vertices of $T$ that belong both to $T_1$ and $T_2$. For each $v$ in $\delta V$ there is an edge toward each connected component of $T_1$ and $T_2$. Hence the graph has $|\mathcal V|:= |\mathcal C \mathcal C(T)| +  | \delta V| $ vertices, it has $|\mathcal E|:= 2\delta V$ edges and since $K_1+K_2 = |\mathcal C \mathcal C(T)|$ we get $\eta =  |\mathcal V| -|\mathcal E| - 1$. By the relation between the number of vertices and the number of edges in a graph applied to $\mathcal G \mathcal C \mathcal C(T)$ (Lemma \ref{lem:EdgesVertices}), we get that $\eta\leq 0$ with equality if and only if $\mathcal G \mathcal C \mathcal C(T)$ is a tree. We then get the expected result: for any $^*$-test graph $T$,
\begin{eqnarray*}
		 \lefteqn{  \tau^0_N \big[ T(\mathbf A_1  , \mathbf A_2 ) \big]}\\   & = & \bigg(  \one\big( \mathcal G\mathcal C \mathcal C (T) \mathrm{ \ is \ a \ tree\,} \big) + O\big(\frac1 N\big) \bigg) \times  \bigg( \prod_{i=1}^2 \prod_{k=1}^{K_i} \tau_N^0  \big[ T_{i,k} (\mathbf A_i)\big]  + o(1) \bigg).
		   \end{eqnarray*}

\begin{remark}\label{Rk:Rk:MainThBack}
\begin{enumerate}
	\item Suppose that $\mathbf A_1$ and $\mathbf A_2$ do not satisfy the factorization property, namely assumption (3) of Theorem \ref{MainTh}. Suppose instead that $\esp\big[ \prod_{i=1}^K \frac 1 N \Tr[T_k(\mathbf A_j)\big]$ converges for any test graphs $T_1\etc T_K$, $j=1,2$. Then it remains true that $(\mathbf A_1, \mathbf A_2)$ converges in traffic distribution. The families of matrices are not asymptotically traffic independent. Indeed, we see from Lemma \ref{Lem:EqDeccorProp} that $\esp\big[ \prod_{i=1}^K \frac 1 N \Tr^0[T_k(\mathbf A_j)\big]$ converges also. It will be clear in the next step that we also get that $\esp\big[ \prod_{i=1}^K \frac 1 N \Tr[T_k(\mathbf A_1, \mathbf A_2)\big]$ converges.
	\item Assume that $\mathbf A_1$ consists in a family of diagonal matrices. Then it is clear that $\tau^0_N \big[ T(\mathbf A_1  , \mathbf A_2 ) \big]   = \tau^0_N \big[ \tilde T(\mathbf A_1  , \mathbf A_2 ) \big]  $ where $\tilde T$ is obtained by identifying source and target of each edge of $T$. Hence, there is a single colored component of $\tilde T$ whose labels correspond to $\mathbf A_2$. Hence we do not need to assume the factorization property for $\mathbf A_2$..
\end{enumerate}
\end{remark}

\subsubsection{Proof of the factorization property}

\begin{lemma}\label{lem:PropConvEsp} Let $\mathbf A_N$ be a family of matrices and $S$ be a $^*$-graph whose connected components are $T_1   \etc T_K $. Then,
\eq
	 \Tr^0 \big[ T_1(\mathbf A_N) \big] \dots  \Tr^0 \big[ T_K(\mathbf A_N) \big]  
	  & = &  \sum_{\pi\in \mathcal P^*(V)}     \Tr^0 \big[ S^\pi(\mathbf A_N) \big] ,
\qe

\noindent where $\mathcal P^*(V)$ is the set of partitions of $V$ that contain at most one vertex of each $S_k$, $k=1\etc K$. 
\end{lemma}

\begin{proof}[Proof of Lemma \ref{lem:PropConvEsp}] We write $S=(V,E,\gamma, \varepsilon)$ and denote by $V_k$ the set of vertices of $S_k$, $k=1\etc K$. Then, denoting $\phi( e ) = \big(\phi( w ), \phi( v )\big)$ for any edge $e=(v,w)$,
\eq
	 \Tr^0 \big[ T_1(\mathbf A_N) \big] \dots \Tr^0 \big[ T_K(\mathbf A_N) \big]  =  \sum_{\phi }\prod_{e \in E} A_{\gamma(e)}^{\varepsilon(e)} \big( \phi( e ) \big),
\qe
where the sum is over all maps $\phi: V \to \{1 \etc N\} $ such that $\phi_{|V_1}\etc \phi_{|V_K}$ are injective. The sum over $\pi$ in the lemma represents all the possible situations of overlapping of the images of $\phi_{|V_1}\etc \phi_{|V_K}$.
\end{proof}

\noindent Let $S, T_1\etc T_K$ be as in the lemma:
\eqa	\label{eq:ProofAsympFreeness:MultiCompo}
	\esp\Big[   \prod_{i=1}^K \frac 1 N \Tr^0 \big[ T_i(\mathbf A_1 , \mathbf A_2 ) \big]  \Big]  =   \sum_{\pi \in \mathcal P^*(V)}  \frac 1 {N^{K-1}}   \tau_N^0 \big[ S^\pi(\mathbf A_1 , \mathbf A_2 ) \big].
\qea

\noindent Let $\pi \in \mathcal P^*(V)$ and denote by $K_\pi\leq K$ the number of components of $S^\pi$. We apply Lemmas \ref{Prop:PermutationInvariance} and \ref{Lem:IndComp} as in the previous step with $T=S^\pi$, with the same notations:
		$$  \frac 1 {N^{K_\pi-1}} \tau_N^0\big[ T(\mathbf A_1  , \mathbf A_2 )\big] = \tilde \Gamma_N \times \bigg( \prod_{i=1}^2 \esp\Big[ \prod_{k=1}^{K_i} \frac 1 N \Tr^0  \big[ T_{i,k}(\mathbf A_i) \big]\Big]  + O\big(\frac1 N\big) \bigg)$$
		where the $T_{i,k}$ are the colored components of $T$ and now $\tilde \Gamma_N =  N^{|\mathcal V| -|\mathcal E| - K_\pi}\times \big( 1+O\big( \frac 1 N \big) \big)$. Note that $K_\pi$ is also the number of connected components of $\mathcal G\mathcal C \mathcal C(T)$, so by Lemma \ref{lem:EdgesVertices}, we get that $\tilde \Gamma_N \limN 1$ if $\mathcal G\mathcal C \mathcal C(T)$ is a forest and it tends to zero otherwise. Hence the only partition $\pi$ which contributes in (\ref{eq:ProofAsympFreeness:MultiCompo}) is the trivial partition and we get
		$$\esp\Big[  \prod_{i=1}^K \frac 1 N \Tr^0 \big[ T_i(\mathbf A_1 , \mathbf A_2 ) \big]\Big] 
		 \limN
	  \prod_{i=1}^K \Nlim \tau_N^0\big[T_i(\mathbf A_1, \mathbf A_2) \big] .$$
This yields, together with Lemma \ref{Lem:EqDeccorProp}, that $(\mathbf A_1, \mathbf A_2)$ satisfies the factorization property.

\chapter{Examples and Applications for Classical Large Matrices}\label{Sec:ApplWignerMatrices}

We prove the convergence in traffic distribution of a family of independent Wigner matrices and give several applications. We also prove the convergence in traffic distribution of a uniform permutation matrix and of a Haar unitary matrix, and compare these two models. The last section is devoted to the case of Diagonal matrices and to a remark on the factorization property.

 We denote as usual $\Phi_N: P \mapsto  \esp\big[\frac 1 N \Tr \, P \big]$ the expectation of the normalized trace of $^*$-polynomials, $\tau_N: T\mapsto   \esp\big[\frac 1 N \Tr \, T \big]$ the expectation of the normalized trace of $^*$-test graphs and $\tau_N^0 = \esp\big[ \frac 1N \Tr^0 \, \cdot \, \big]$ its injective version.

\section[Wigner matrices]{Wigner matrices}\label{Sec:ConvWigner}

\subsection{Limiting distribution of independent matrices}

\begin{figure}  
 \includegraphics{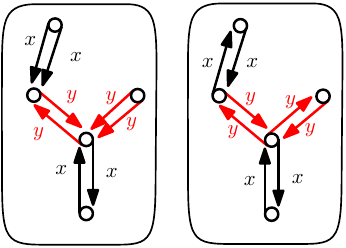} 
 \caption{Two test graphs. The leftmost contributes for the limiting injective traffic distribution of independent Wigner matrices. The rightmost contributes when the parameter of the matrices is of the form $(\alpha,0)$}
\label{fig:Semicircular}
\end{figure}

\begin{proposition}[The limiting distribution of Wigner matrices]\label{Prop:WigDistr} Let $\mathbf X_N=(X_j)_{j\in J}$ be a family of independent Wigner matrices (Definition \ref{Def:WigMatrices}), whose entries have the same law as their complex conjugates. Then $\mathbf X_N$ has a limiting traffic distribution. Denote the parameter of $X_j$ by $(\alpha_j , \beta_j)$ for any $j\in J$. Let $T=(V,E,\gamma)$ be a test graph in variables $\mathbf x=(x_j)_{j\in J}$, with no edge labeled $x^*$ (we deduce the general distribution as the matrices are Hermitian).

 Say that $T$ is a fat tree whenever it becomes a tree if the multiplicity of the edges and the orientation are forgotten. It is called a double tree if moreover the multiplicity of the edges is always two, see Figure \ref{fig:Semicircular}. Let call \emph{twin edges of a double tree} two edges that share the same vertices. We say that $T$ is colored if twin edges $e,e'$ of $T$ have the same label $\gamma(e)=\gamma(e')\in J$ so they correspond to the same matrix $X_{\gamma(e)}$. For a double tree $T$ we denote by $\ell_j(T)$ (respectively $k_j(T)$) the number of twin edges labeled $x_j$  with opposing (respectively similar) orientation.
 
 Then for any test graph $T$ one has
	\eqa\label{Eq:WignInjMoments3}
		  \tau^0_N\big[ T(\mathbf X_N) \big]  \limN   \one \big( T \textrm{ is a colored double tree}\big)
		  									\prod_{j\in J} \alpha_j^{\ell_j(T)}\beta_j^{k_j(T)}.
	\qea
\end{proposition}

 In particular, if the matrices are real Wigner matrices with parameter $(1,1)$, $\tau^0_N\big[ T(\mathbf X_N) \big]$ is asymptotically the indicator that $T$ is a colored double tree, and if they are complex matrices with parameter $(1,0)$, e.g. for a GUE matrix, the assumptions that twin edges of $T$ have different orientation is required.

\begin{proof}[Proof of Lemma \ref{Prop:WigDistr}]
 \par Without lose of generality we assume $\alpha_j=1$ for any $j\in J$. Denote $M_j=\sqrt N X_j$ and $\mathbf M_N = (M_j)_{j\in J}$. Consider a test graph $T=(V,E,\gamma)$ in variables $\mathbf x$. By multi-linearity of $\Tr^0\big[T(\, \cdot \, )\big]$ with respect to the edges of $T$ and since the family is permutation invariant, one has 
		$$  \tau_N^0\big[ T( \mathbf  X_N) \big] = N^{-\frac{|E|}2} \tau_N^0\big[ T(  \mathbf  M_N) \big]  = N^{|V|-1-\frac{|E|}2}\delta_N^0\big[T(  \mathbf M_N)]  \big(  1+O(N^{-1})\big),$$
 where $\delta_N^0\big[T( \mathbf  M_N)] = \esp \big[ \prod_{(v,w)\in E} M_\gamma(e) \big( \phi(w), \phi(v) \big)\big]$, as in Formula \eqref{Eq:LienTr0Delt0}, for any injection $\phi:V \to [N]$. The quantity $\delta_N^0\big[T( \mathbf  M_N)] $ is bounded and can be computed since the entries of $M_j$ are independent and independent of $N$. By centering of the entries it is zero whenever $T$ has an edge of multiplicity one. Let $T$ such that each edge has at least of multiplicity two.

 We apply Lemma \ref{lem:EdgesVertices} to the graph $\mathcal G = ( V, \bar E)$ obtained from $G$ by forgetting the multiplicity if its edges (and their orientation and labels). Since $N^{|V|-1-|E|/2} = N^{|V|-1-|\bar E|}\times N^{|\bar E|-|E|/2}$, we then get that $ \tau_N^0\big[ T( \mathbf X_N) \big] $ converges to zero if it $T$ is not a double tree. Moreover, by independence of the entries, when $T$ is a double tree then $ \tau_N^0\big[ T(  X_N) \big] =  \delta_N^0\big[ T(  M_N) \big] $ is the product of terms of the form $\esp\big [M_N(k,\ell )^2 \big]$ or $\esp\big [M_N(\ell,k)^2 \big]$ for $k\neq \ell$ along each twin edges. Hence if the double tree is not colored then  $ \tau_N^0\big[ T( \mathbf X_N) \big] $ converges to zero. Otherwise, the independence of the matrices and their entries give the expected formula.
\end{proof}

\subsection{Practical computations}
The limiting traffic distribution of independent Wigner matrices yields the following computations.
\begin{enumerate}
	\item Limit of $\Phi_N\big( P(\mathbf X_N)\big)$. Let $n\geq 1$ be an integer and let us first prove Wigner's law, namely for $\alpha=1$ the convergence of 
	$$m_n\toN:=\esp \big[ \frac 1 N \Tr X_N^n \big] \limN a_{\frac n 2},$$
 where $a_\ell$ is zero if $\ell$ is not an integer and is the $\ell$-Catalan number $\frac1 {\ell+1}\binom{2\ell}{\ell}$ otherwise. This proof is similar to the one in \cite{GUI}.
	
	We know from Example \ref{Ex:DistrFromTraffics} that we can write  $m_n\toN = \tau_N\big[ T_n(X_N) \big],$ where $T_n$ is the simple circle of length $n$ with edges oriented along the cycle (and all edges have the same label, associated to the matrix $X_N$). Now, denoting by $V_n$ the vertex set of $T_n$, by the relation \eqref{eq:TraffCum} relating the trace and the injective trace and by Lemma \ref{Prop:WigDistr} one gets
	$$m_n\toN= \sum_{\pi\in  V_n} \tau^0_N\big[ T^\pi_n(X_N) \big] = \sum_{\pi\in  V_n} \one(T^\pi_n \mathrm{ \ is \ a \ double \ tree}) +o(1).$$
We then see that the $n$-th moment $m_n\toN$ converges to the number of double trees $T^\pi_n$ we can obtain from the simple cycle $T_n$. This number is zero if $n$ is odd. Moreover, choosing a double tree is equivalent to chose the pairs that form double edges. These pair partitions are the non crossing partitions \cite{GUI}, for which it is known that the cardinals are the Catalan numbers \cite{NS}. 

For several independent Wigner matrices, the methods extends to get from Proposition \ref{Prop:WigDistr} the convergence of $\mathbf X_N$ to a free semicircular system, see \cite{CDM16}.
	\item {Limit of $\Phi_N\big( P(\mathbf X_N) \circ Q(\mathbf X_N) \big)$.} We can see that for any $^*$-polynomials $P,Q$, one has 
		$\Phi_N\big( P(\mathbf X_N) \circ Q(\mathbf X_N) \big) = \Phi_N\big( P(\mathbf X_N) \big) \Phi_N\big( Q(\mathbf X_N)\big) +o(1).$
	Indeed, according to the second application of Example \ref{Ex:DistrFromTraffics}, for any monomials $P, Q$, one has $\Phi_N\big( P(\mathbf X_N) \circ Q(\mathbf X_N)  \big) = \tau_N\big[T(\mathbf X_N)\big]$ where $T$ is the $^*$-test graph consisting in a simple cycle $T_1$ of length the degree of $P$ and a simple cycle $T_2$ of length the one of $Q$, both oriented and attached together by identifying a single vertex of each cycle. Edges are labeled accordingly to the monomials $P$ and $Q$ following the orientation of the respective cycle. 
	
	The partitions $\pi$ of $T$ such that $T^\pi$ is a double tree are those for which the induced subgraphs of $T_1$ and $T_2$ are double trees that have in common a single vertex, given by two partitions $\pi_1$ of $T_1$ and $\pi_2$ of $T_2$.  Hence we have 
	\eq
		 \tau_N\big[T(\mathbf X_N)\big] & = & \sum_{\substack{\pi_1 \in \mathcal P(V_1)\\\pi_2 \in \mathcal P(V_2)}} \tau^0_N\big[T_1(\mathbf X_N)\big]  \tau^0_N\big[T_2(\mathbf X_N)\big]  +o(1)\\
		 & = & \tau_N\big[T_1(\mathbf X_N)\big]  \tau_N\big[T_2(\mathbf X_N)\big]  +o(1)\\
		 & = &\Phi\big( P(\mathbf X_N) \circ Q(\mathbf X_N) \big) +o(1).
	\qe
 The result follows by linearity.
\end{enumerate}

\subsection{First example of asymptotic traffic independent matrices}
By a direct application of Proposition \ref{Prop:WigDistr} (without using Theorem \ref{MainTh}), we get:
\begin{lemma} Independent Wigner matrices whose entries have the same law as their complex conjugates are asymptotically traffic independent.
\end{lemma}

\begin{proof} Let $\mathbf X_N=(X_j)_{j\in J}$ be such a family. We have seen under the assumptions of Lemma \ref{Prop:WigDistr} that $\mathbf X_N$ converges in traffic distribution. Let $T=(V,E,\gamma,\varepsilon)$ be a $^*$-test graph in variables $\mathbf x$. We have to prove that $T$ is a colored double tree if and only if its graph of colored components is a tree and its colored components are double trees. This fact can be clearly seen in pictures, see Figure \ref{fig:SC01}. We prove it with the same method as in Lemma \ref{Prop:WigDistr}.

  \begin{figure}
     \includegraphics{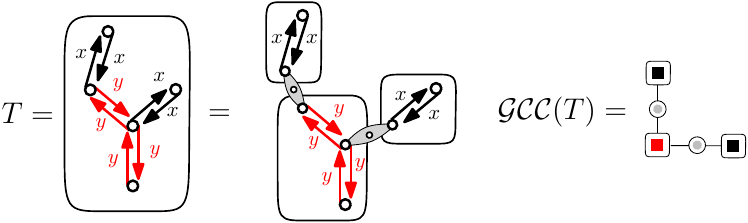}
\caption{ The decomposition in colored components of a colored double tree.  \label{fig:SC01}  }
  \end{figure}

\par If $T$ has an edge of multiplicity one then $T$ is not a double tree and it has a colored component that is not a double tree. So both the limiting distribution of Wigner matrices and the product of the marginal limits vanish on such graphs.
\par Assume from now $T$ has no edge of multiplicity one. For each $S\in \mathcal C \mathcal C(T)$ colored component of $T$, denote by $E_S$ the set of edges of $T$ and by $|\bar E_S|$ its number of edges without multiplicity. By Lemma \ref{lem:EdgesVertices} and the proof of Lemma \ref{Prop:WigDistr}, for a graph $T$ with no simple edge, the quantity
	\eq
		\eta(T) =  \Big(   |\bar E| -\sum_{S\in \mathcal C \mathcal C(T)} |\bar E_S| \Big) + \Big(  |\bar E| - \frac{ |E|}2 \Big) + \Big( |V| - 1 - |\bar E| \Big)   
	\qe
	is zero if and only if $T$ is a colored double tree. We can rewrite $\eta(T)$ in terms of quantities defined on the graph of colored components of $T$ as follow. First for each colored component $S$ of $T$ (with respect to the $x_j,j\in J$) denote by $V_S$ its set of vertices. Denote also $\delta V_S$ the set of vertices of $S$ that are also contained in another colored component, i.e $\delta V_S = V_S\cap \big( \cup_{S'\neq S} V_{S'} \big)$. Denote $\overset{\circ} V_{S} = V_S \setminus \delta V_S$ and $\delta V = \cup_S \delta V_S$. Note that we have the equality
	\eq
		\sum_{S}( |V_S| - 1 )   = \sum_{S}( |\delta V_S| + |\overset{\circ} V_{S}|- 1 )    = |V| - |\delta V| +  \sum_{S} |\delta V_S| -  \big| \mathcal C\mathcal C(T)\big| .
	\qe
Moreover, denoting $|\mathcal V|$ and $| \mathcal E|$ the number of vertices and edges of the graph of colored components of $T$, one has
	$$( |\mathcal V| - 1 -| \mathcal E|) =  \big| \mathcal C\mathcal C(T)\big| + |\delta V| - 1 - \sum_{S} |\delta V_S|.$$
Hence summing the two above identities yields an expression for $|V|-1$ that gives
	$$	\eta(T) = \Big(  |\bar E| - \frac{ |E|}2 \Big)+  ( |\mathcal V| - 1 -| \mathcal E|) + \sum_{S}( |V_S| - 1- |\bar E_S| ).$$
Since $T$ has no simple edges, the three terms in the r.h.s are nonpositive the above equation is zero if and only the multiplicity of the edges of $T$ is two, the graph of colored components of $T$ is a tree and the colored components are double trees. In particular, we get that for real Wigner matrices the limiting traffic distribution of $\mathbf X_N$ is indeed the product of the limiting distributions of the $X_j$'s. 

Let now consider the case of complex Wigner matrices. It is clear that the twin edges of $T$ are the twin edges of its colored components which yields again that the distribution of Lemma \ref{Prop:WigDistr} satisfies \eqref{DefTrafInd}.
\end{proof}

\subsection{Asymptotic free independence with the transpose}

\begin{lemma}\label{Lem:FreeTranspose} Let $\mathbf X_N=(X_j)_{j\in J}$ be a family of random matrices. Assume that it converges in traffic distribution and has the same limiting distribution as a family of independent complex Wigner matrices with parameter $(\alpha_j,0)_{j\in J}$. Then the matrices $X_j,X_j^t$, $j\in J$, are asymptotically free independent.
\end{lemma}

Asymptotic free independence with the transposed ensemble is studied for a larger class of examples in \cite{CDM16,MiPo14}.

\begin{proof} We show that $(X_j,X_j^t, j\in J)$ has the same limiting $^*$-distribution as $(Y_j,Y'_j, j\in J)$, where the $Y_j,Y'_j,j\in J$ are independent and $Y_j,Y_j'$ are distributed as $X_j$ for each $j\in J$. Since the matrices $Y_j,Y'_j, j\in J$, are
asymptotically $^*$-free this will yield the expected result.

By Lemma \ref{Lem:EquivCVTraf} the limiting $^*$-distribution of $\mathbf X_N \cup \mathbf X_N^t$, $j\in J$ is entirely determined by the traffic distribution of $\mathbf X_N$. More precisely, let $T$ be a test graph in the variables $(y,y')$. When there is a cycle visiting each edge once in the sense of their orientation, we say that $T$ is cyclic. By Example \ref{Ex:DistrFromTraffics} and since the quotient of cyclic graph is cyclic, the limits of $\tau_N^0\big[ T(\mathbf X_N,\mathbf X_N^t)\big]$ for any cyclic $^*$-test graph characterize the limiting joint $^*$-distribution of $(\mathbf X_N,\mathbf X_N^t)$.
\par Moreover, we have $\tau_N^0\big[ T(\mathbf X_N,\mathbf X_N^t)\big] = \tau_N^0\big[ \tilde T(\mathbf X_N) \big]$, where the orientation of the edges labelled $y'$ has been reversed and labels $y'$ where replaced by $y$. Hence, by the expression of the limiting distribution of $\mathbf X_N$, $\tau_N^0\big[ T(\mathbf X_N,\mathbf X_N^t)\big]$ converges to one whenever $T$ is a double tree whose twin edges have same label and opposite orientation or different labels and same orientation, and it converges to zero otherwise. But since $T$ is cyclic, it is not possible in a double tree that twin edges have same orientation. Hence for these graphs $T$ the quantities $\tau_N^0\big[ T(\mathbf X_N,\mathbf X_N^t)\big]$ and $\tau_N^0\big[ T(\mathbf Y_N,\mathbf Y'_N)\big] $ have the same limit, which concludes the proof.
\end{proof}

\subsection{Non-asymptotic traffic independent matrices}

\begin{lemma} A nonzero complex Wigner matrix $X_N$ is not asymptotically traffic independent from its transpose $X_N^t$.
\end{lemma}

Together with Lemma \ref{Lem:FreeTranspose}, this gives an example of asymptotically free independent matrices that are not asymptotically traffic independent.

\begin{proof} We have $\Phi_N(X_N^2) = \Phi_N\big(X_N(X_N^t)^t\big) =  \tau\big[ \cdot \overset{X_N}{\underset{X_N^t}\leftleftarrows} \cdot] = \tau^0\big[ \cdot \overset{X_N}{\underset{X_N^t}\leftleftarrows} \cdot] +o(1).$ Assume that $X_N$ and $X_N^t$ are asymptotically traffic independent. The graph of colored component of the graph $( \cdot \overset{x}{\underset{x^t}\leftleftarrows} \cdot)$ with respect to the variables $x$ and $x^t$ is not a tree, so we should have $\Phi_N(X_N^2) \limN 0$. But $\Phi_N(X_N^2)$ converges to the variance of $\sqrt N X_N(i,j)$ for $i\neq j$ which is nonzero since $X_N$ is a nonzero Wigner matrix.
\end{proof}

  \begin{figure}
     \includegraphics{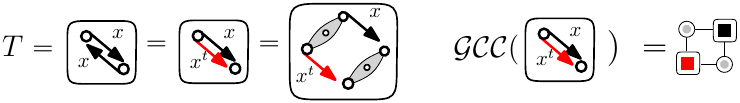}
\caption{   \label{fig:ExamplesFreeness2}  A graph of colored component of a double tree which is not a tree.  }
  \end{figure}

\subsection{Factorization property}

\begin{lemma}\label{lem:ConvEspWigner} A Wigner matrix $X_N$ satisfies the factorization property (\ref{eq:DecorrIntro}).
\end{lemma}

\begin{proof}
\noindent Let $T_1   \etc T_n $ be test graphs in one variable, and denote by $T$ the graph obtained as the disjoint union of $T_1 \etc T_n$. By Lemma \ref{lem:PropConvEsp}, 
\eq
	\esp\Big[ \prod_{i=1}^n \frac 1 N \Tr^0 \big[ T_i(X_N) \big] \Big] 
	  & = &  \sum_{\pi}  \frac 1 {N^{n}}   \esp\Big[  \Tr^0 \big[ T^\pi(X_N) \big]  \Big]\\
	  & = &  \sum_{\pi}    N^{V_\pi-\frac{E} 2-n}\Big( 1 +O\big( \frac 1 N \big) \Big) \delta_N^0\big[T^\pi(\sqrt N X_N) \big]  \Big],
	 \qe
where the sum is over all partitions whose blocks contain at most one vertex of each $T_i$ and $V_\pi$ denotes the vertex set of $T^\pi$. We have by Lemma \ref{lem:EdgesVertices} that $V_\pi-\frac E 2-n\leq 0$ with equality if and only if $\pi$ is the trivial partition with only singleton blocks and the graphs $T_1\etc T_n$ are double trees. By the independence of the entries of $X_N$, we obtain  the factorization property
\eq
	  \lefteqn{ \esp\Big[ \prod_{i=1}^n \frac 1 N \Tr^0 \big[ T_i(X_N) \big] \Big]  }\\
	  &  =  & \prod_{i=1}^{n} \one(T_i \textrm{ is a double tree}) \delta_N^0\big[T_i(\sqrt N X_N) \big] +o(1) \\
	  & \limN & \prod_{i=1}^{n} \Nlim \tau_N^0\big[T_i( X_N) \big].
	  \qe
\end{proof}

\section{Unitary Haar and uniform permutation matrices}\label{Sec:ExamplesMatrices}

\subsection{Limits in traffic distribution}
Let now compare the limiting traffic distributions of uniform permutation matrices and unitary Haar matrices.

\begin{figure}
 \includegraphics{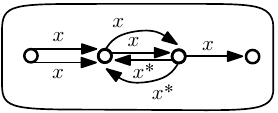}  \caption{     A $^*$-test graph that contributes for the injective trace of large uniform permutation matrices}\label{fig:ExamplesFreeness3}
\end{figure}

\begin{proposition}[The limiting distribution of a uniform permutation matrix]\label{Prop:PermDistr} A uniform permutation matrix $V_N$ has a limiting traffic distribution and satisfies the factorization property. Say that a $^*$-test graph is a directed line whenever there is an integer $K\geq 1$ such that the vertices of $T$ are $1\etc K$ and its directed edges are $(1,2) \etc (K-1,K)$ labelled $x$ and $(2,1) \etc (K,K-1)$ labelled $x^*$, with arbitrary multiplicity (see Figure \ref{fig:ExamplesFreeness3}). Then, for any $^*$-test graph $T$ in one variable, 
\eqa\label{eq:DefPerm}
		\tau_N^0\big[ T(V_N) \big] & \limN & \left\{ \begin{array}{cc} 1 & \textrm{ if } T \textrm{ is a directed line}\\
												0 & \textrm{ otherwise}. \end{array} \right.
	\qea

\end{proposition}

\begin{figure}
 \includegraphics{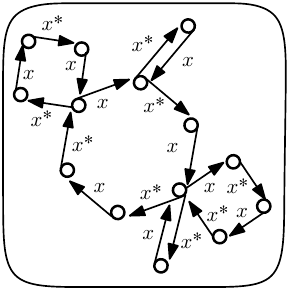} 
 \caption{     A $^*$-test graph that contributes for the injective trace of large Haar unitary matrices}
 \label{fig:ExamplesFreeness4}
 \end{figure}

\begin{proposition}[The limiting distribution of a unitary Haar matrix]\label{Prop:HaarDist} A unitary Haar matrix $U_N$ converges in traffic distribution and satisfies the factorization property. Let call simple oriented cycle of a $^*$-test graph $T$ a subgraph of $T$ consisting in a closed oriented path with no repetition of vertices, as in Example \ref{Ex:DistrFromTraffics}. The support of $\tau^0$ is the set of $^*$-test graphs $T$ such that each edge belongs to a unique simple cycle of $T$, and each simple cycle of $T$ is oriented and of even size with edges labelled $x$ and $x^*$ alternately (see Figure \ref{fig:ExamplesFreeness4}). For such a graph $T$, one has $\tau^0[T] = \prod_{i=1}^L c_{k_i/2-1}(-1)^{k_i/2-1}$ where $k_1\etc k_L$ are the length of the $L$ simple cycles of $T$ and $c_k=\frac{(2k)!}{(k+1)!k!}$ is the $k$-th Catalan number.
\end{proposition}

Before proving this statement (see Section \ref{Sec:ExMatricesProof}), we show some applications.

\subsection{Practical computations} 
\begin{enumerate}
	\item Limiting $^*$-distribution: Let $V_N$ be a uniform permutation matrix and $ U_N$ be a unitary Haar matrix. For any unitary matrix $W_N$, the $^*$-distribution of $W_N$ depends only the limits of $\esp \big[\frac 1 N \Tr W_N^k (W_N^*)^{\ell}\big]$ for any $k,\ell$. By \cite{NIC93} and \cite{VOI5}, for both $V_N$ and $U_N$ this limit is $\one(k=\ell)$. Let prove this result from the above propositions.

Let $W_N \in \{V_N, U_N\}$. Let $T=T_{k,\ell}$ be the $^*$-test graph in variables $w,w^*$ consisting in a simple cycle with $k$ edges labelled $w$ followed by $\ell$ edges labeled $w^*$ in such a way $\esp \big[\frac 1 N \Tr W_N^k (W_N^*)^{\ell}\big] = \tau_N\big[T(W_N)\big]$: denote by $1 \etc k+\ell$ its vertices, by $e_i=(i,i+1)$, $i=1\etc k$ its edges labelled $w$ and by $e_{k+i}=(k+i,k+i+1)$, $i=1\etc \ell$ its edges labelled $w^*$ (with the convention $i_{k+\ell+1}=i_1$).

There is a partition $\pi$ of the vertex set of $T$ such that $T^\pi$ is a line directed if and only if $k=\ell$ and in that case this partition is necessarily$\big\{ \{1\}, \{k\}, \{i, 2k-i\}, i=2\etc k-1\big\}$. Hence $\tau_N^0\big[T(V_N)\big] \limN \one(k=\ell)$. Moreover, there is a partition $\pi$ such that $T^\pi$ is a $^*$-test graph that contributes to the limiting traffic distribution of a unitary Haar matrix if and only if $k=\ell$, and then $\pi$ must be the same partition. Since all simple cycles of $T^\pi$ are of length two, when $k=\ell$ one has $\tau_N^0\big[T(U_N)\big] \limN c_0^k = 1$. We then obtained $\tau_N\big[T_{k,\ell}(W_N)\big] \limN \one(k=\ell)$ for $W_N = V_N,U_N$.
	\item Limits of entry-wise products: We prove that for $W_N$ a uniform permutation or unitary Haar matrix $\Phi_N\big( P(W_N) \circ Q(W_N)\big) -\Phi_N\big( P(W_N)\big) \times  \Phi_N\big( Q(W_N)\big) = o(1)$.
	
	As before we can assume $P(w)=w^{k_1}(w^*)^{\ell_1}$ and $Q=w^{k_2}(w^*)^{\ell_2}$. Denote by $T$ the $^*$-test graph in variables $w,w^*$ consisting in a bunch of two simple cycles, one with $k_1$ edges labelled $w$ followed by $\ell_1$ edges labeled $w^*$ and the second defined similarly with $(k_2,\ell_2)$ instead of $(k_1,\ell_2)$, in such a way $\esp \big[\frac 1 N \Tr (W_N^{k_1} (W_N^*)^{\ell_1}) \circ (W_N^{k_2} (W_N^*)^{\ell_2}) \big] = \tau_N\big[T(W_N)\big]$ (see Example \ref{Ex:DistrFromTraffics}).
	There is a partition $\pi$ of the vertex set of $T$ such that $T^\pi$ is a directed line if and only if $k_1=\ell_1$ and $k_2=\ell_2$. In that case the partition must be the one which makes the identification of each circle into a directed line, and then identifying the vertices of the shortest line into the one of the longest to create a directed line by starting for the point of identification of the vertices. Hence $\tau_N^0\big[T(V_N)\big] \limN  \one(k_1=\ell_1) \one(k_2=\ell_2)$. Moreover, $T^\pi$ contributes to the limiting traffic distribution of a unitary Haar matrix if and only the quotient subgraphs $T_1$ and $T_2$ induced by the circle also contribute and have in common a single vertex. In that case, since the limiting injective distribution of unitary Haar matrices is multiplicative with respect to the simple cycles of $T^\pi$, $\tau_N^0\big[T(U_N)\big] = \tau_N^0\big[T_1(U_N)\big] \times \tau_N^0\big[T_2(U_N)\big]+o(1) \limN  \one(k_1=\ell_1) \one(k_2=\ell_2)$.
	
\end{enumerate}

\subsection{Unitary Haar matrices and their transpose} Let $\mathbf U_N=(U_j)_{j\in J}$ be a family of independent unitary Haar matrices. By Voiculescu's theorem (Theorem \ref{Th:AsyFree}), the matrices $U_j, j\in J$ are asymptotically free independent (it will also be a consequence of Proposition \ref{RigFree}). We also have a similar statement as for complex Wigner matrices, Lemma \ref{Lem:FreeTranspose}, the matrices $U_j, U_j^*, j\in J$ are asymptotically free independent. The proof is similar to the one of Lemma \ref{Lem:FreeTranspose}, we get that $(\mathbf U_N, \mathbf U_N^t)$ has the same limiting $^*$-distribution as $(\mathbf U_N, \tilde {\mathbf U}_N)$, where $\tilde {\mathbf U}_N$ is an independent copy of $\mathbf U_N$, replacing twin edges by simple cycles of arbitrary length in the reasoning. 

\subsection{Multiplication of entries by independent weights} We now give a concrete example of Corollary \ref{Cor:VarProfile}.

\begin{proposition}\label{Prop:AppliMultEntry} Let $U_N$ be a unitary Haar matrix and let $V_N$ be a uniform permutation matrix. Let $W = (\omega_{i,j})_{i,j\geq 1}$ be an infinite array of independent identically distributed complex entries, independent of the parameter $N$, independent of $(U_N,V_N)$, and such that $\esp[|\omega_{i,j}|^K]$ is finite for any $K\geq 0$. Denote by $W_N$ the $N$ by $N$ matrix $(\omega_{i,j})_{i,j=1\etc N}$. Recall that the symbol $\circ$ denotes the entry wise product of matrices.
\begin{enumerate}
	\item Denote $M_N = U_N\circ W_N$ and let $\mathbf M_N$ be a family of independent copies of $M_N$. If $\omega_{i,j}$ is centered, then the matrices of $\mathbf M_N$ are asymptotically free independent and $\mathbf M_N$ has the same limiting traffic distribution as the matrix $\frac 1 {\sqrt N} \tilde W_N= \frac 1 {\sqrt N} (\tilde \omega_{i,j})_{i,j=1\etc N}$ with centered independent complex Gaussian entries with covariance $\esp[|\tilde \omega_{i,j}|^2] = \esp[| \omega_{i,j}|^2]$ and $\esp[\tilde \omega_{i,j}^2]=0$.
	\item Denote $M_N = V_N\circ W_N$ and let $\mathbf M_N$ be a family of independent copies of $M_N$. If the modulus of $\omega_{i,j}$ is not deterministic, then the matrices of $\mathbf M_N$ are not asymptotically free independent but are asymptotically traffic independent.
	\end{enumerate}
\end{proposition}

\begin{proof} The convergence in traffic distribution of $\mathbf M_N$ is consequence of Corollary \ref{Cor:VarProfile} and Theorem \ref{MainTh}.

{\bf 1.} By Corollary \ref{Cor:VarProfile}, we know that for any $^*$-test graph $T$, $\tau_N^0\big[ T(M_N)\big] = \tau_N^0\big[ T(U_N)\big] \times  \delta_N^0\big[ T(W_N)\big]$. Since the entries of $W_N$ are independent and centered, then $ \delta_N^0\big[ T(W_N)\big]=0$ if $T$ has an edge of multiplicity one. Moreover, if $\tau_N^0\big[ T(U_N)\big]$ does not converges to zero and has no edge of multiplicity one, then $T$ is a double tree whose twin edges have different label and orientation. For such $T$ with $2k$ edges, the quantity $ \delta_N^0\big[ T(W_N)\big]$ is $\esp[|\omega_{i,j}|^2]^{k}$. 
\par On the other hand let $\tilde W_N$ the matrix with independent Gaussian entries of the lemma and let us compute its limiting traffic distribution. The arguments are the same as for Wigner matrices. By Lemma \ref{Lem:Tool}, one has $\tau_N^0\big[ T(\frac 1 {\sqrt N}\tilde W_N )\big] = N^{|V|-1}\big( 1 + o\big( \frac 1 N\big) \big) N^{\frac{-| E|}2}\delta_N^0\big[T(\tilde W_N)\big]$ where $V$ and $E$ are the vertex and edge sets of $T$. This quantity is zero if $T$ has an edge of multiplicity one, and we assume from now that the edges are of multiplicity at least two. Denote by $|\bar E|$ its number of edges without multiplicity and write $|V|-1 -\frac{| E|}2 = (|V|-1 -|\bar E|) + (|\bar E| -\frac{| E|}2)$. Applying Lemma \ref{lem:EdgesVertices} to the graph obtained from $T$ by forgetting the multiplicity of its edges, we get that the only $^*$-test graphs $T$ for which $\tau_N^0\big[ T(\frac 1 {\sqrt N}\tilde W_N )\big] $ possibly does not vanish at infinity are the double trees. By independence of the entries, $\delta_N^0\big[T(\tilde W_N)\big]$ is multiplicative with respect to the twin edges of $T$. Since for the considered complex Gaussian random variable $\tilde \omega$ one has $\esp[\tilde \omega^2]=0$, then the only graphs that contribute are those for which the twin edges have opposite orientation and different labels. Hence $M_N$ and $\frac 1 {\sqrt N}\tilde W_N $ have the same limiting traffic distribution.

In particular $\mathbf M_N$ has the same limiting $^*$-distribution as independent copies of $\frac 1 {\sqrt N}\tilde W_N $. It is known that such matrices are asymptotically free independent (this is also a consequence of Theorem \ref{RigFree}). Hence the result.
\\
\par {\bf 2.} Let now denote $M_N = V_N \circ W_N$ and let prove that $\kappa_N:=\Phi( M_NM_N^* \circ M_NM_N^*)- \Phi( M_NM_N^*) \Phi(M_NM_N^*)$ does not converge to zero. By permutation invariance of $M_N$, for each $i$ in $[N]$ one has $\kappa_N = \mathbb V\mathrm{ar}\big( \sum_{j=1}^N |M_N(i,j)|^2\big)$. But $|M_N(i,j)|^2 = 0$ if $j$ is not the image of $i$ by the permutation associated to $\sigma$ and otherwise $|M_N(i,j)|^2 = |\omega_{i,j}|^2$. But independence of $\sigma$ and $(\omega_{i,j})_{i,j}$, one has $\kappa_N = \mathbb V\mathrm{ar} \big( |\omega_{i,j}|^2\big)$ for any $i,j$. We get the result thanks to Proposition \ref{Prop:Criterion}.
\end{proof}

\subsection{Proof of the convergence and factorization property}\label{Sec:ExMatricesProof}

We now go back to the proofs of Propositions \ref{Prop:PermDistr} and \ref{Prop:HaarDist}, the convergence in traffic distribution for uniform permutation and unitary Haar matrices and the factorization property  for these two models.

\begin{proof}[Proof of Proposition \ref{Prop:PermDistr}] Let $V_N$ be a uniform permutation matrix. First, remark that since the entries of $V_N$ are in $\{0,1\}$, then for any $^*$-test graph $T$ in one variable, $\tau_N^0\big[ T(V_N) \big]=\tau_N^0\big[ \tilde T(V_N) \big]$ where $\tilde T$ is obtained by
\begin{itemize}
	\item reversing the orientation of edges labelled $x^*$ and replacing this label by $x$,
	\item forgetting the multiplicity of each edges (assuming the multiplicity is one).
\end{itemize}
Hence, we can assume $T=\tilde T$ without loss of generality

Moreover, each row and column of $V_N$ has a single nonzero entry. Hence, $\frac 1 N \Tr^0_N\big[ T(V_N) \big]$ is zero as soon as two distinct edges leave (or start from) a same vertex. Hence, there are only two kinds of test graphs that possibly contribute: 
\begin{itemize}
	\item the simple cycles $T_K$ of length $K\geq 1$, with vertices $1\etc K$ and edges $(1,2) \etc (K-1,K),(K,1)$,
	\item the simple paths $S_K$ of length $K\geq 1$, with vertices $1\etc K$ and edges $(1,2) \etc (K-1,K)$.
\end{itemize}

\noindent Let $\sigma_N$ be the random permutation associated to $V_N$. For a simple cycle, the quantity $ \tau^0_N \big[T_K(V_N) \big]$ is the probability that a given integer $i$ in $\{1\etc N\}$ belongs to a cycle of $\sigma_N$ of length $K$. This probability is $\frac 1 N$. Indeed, there are $\frac{ (N-1)!}{(N-K)!} \times (N-K)!$ permutations such that $i$ is contained in a cycle of length $K$ (the first term counts the number of ways to chose the cycle of length $K$ containing $i$, while the second term accounts for the remaining freedom in the permutation). Then we get $ \tau^0_N \big[T_K(V_N)\big] \limN 0.$ At the contrary, for a simple path $ \tau^0_N \big[S_K(V_N)\big] $ is the probability that a given integer $i$ in $\{1\etc N\}$ belongs to a cycle of $\sigma_N$ of length greater than $K$. By the above, one has $ \tau^0_N \big[S_K(V_N)\big] \limN 1.$

Let now prove the factorization property. Let $\sigma_N$ be the permutation of $\{1\etc N\}$ associated to $V_N$. For any $K_1 \etc K_n, L_1 \etc L_m\geq 1$, the number 
		$$\esp\Big[ \prod_{i=1}^n \frac 1 N \Tr^0\big[ T_{K_i}  (V_N)\big]  \times \prod_{i=1}^m \frac 1 N \Tr^0\big[ S_{L_i} (V_N)\big] \Big]$$
is the probability that, choosing $i_1 \etc i_n,j_1\etc j_m$ uniformly and independently on $\{1\etc N\}$ one has
\begin{itemize}
	\item $i_k$ belongs to a cycle of length $K_k$ of $\sigma_N$ for any $k=1\etc n$,
	\item $j_k$ belongs to a cycle of length bigger than $L_k$ of $\sigma_N$ for any $k=1\etc m$.
\end{itemize}
By a straightforward computation, this probability tends to zero or one, depending if $n$ is positive or not respectively.
\end{proof}

\begin{proof}[Proof of Proposition \ref{Prop:HaarDist}] Let $U_N$ be a unitary Haar matrix. Let $T = (V,E,\varepsilon)$ be a $^*$-test graph in the variable $x$. By Formula \eqref{Eq:TraceDensity}, one has 
	$$\tau_N^0\big[ T(U_N) \big] = N^{|V|-1} \big( 1 + O(N^{-1}) \big) \delta_N^0\big[ T(U_N) \big].$$
 To compute the asymptotic of this quantity, we use Weingarten calculus (see \cite[Lecture 23]{NS} and \cite{CS06}). Denote by $(j_1,i_1) \etc (j_k,i_k)$ the edges of $T$ with label $x$ and $(i_1',j_1') \etc (i_{k'}',j_{k'}')$ the edges with label $x^*$, where the $i_n,j_n,i'_n,j'_n$ are integers. The invariance of $U_N$ by conjugation by permutation matrices and Weingarten formula \cite[Lemma 23.5]{NS} tell us that 
	\eq
		\delta_N^0\big[ T(U_N) \big] & = & \esp\big[ U_N(i_1,j_1) \dots U_N(i_k,j_k)  U_N^*(j_1',i_1') \dots U_N^*(j_{k'}',i_{k'})\big]
	\qe 
is zero if $k\neq k'$, and otherwise is equal to 
	\eqa\label{Weingarten}
		\sum_{\sigma,\tau \in \mathcal S_k} \delta_{i_\sigma(1),i'_1} \dots \delta_{i_{\sigma(k)},i'_k}  \delta_{j_{\tau(1)},j'_1} \dots \delta_{j_{\tau(k)},j'_k} \textrm{Wg}(\sigma \circ \tau^{-1},N),
	\qea
where $\mathcal S_k$ is the set of permutation of $\{1\etc k\}$ and $\textrm{Wg}$, known as the Weingarten function, has asymptotic $\textrm{Wg}( \sigma\circ \tau^{-1},N) = \phi(\sigma \circ \tau^{-1}) N^{-|E| + \sharp(\sigma \circ \tau^{-1})} \times \big(1 + O(N^{-2}) \big)$. The number $ \phi(\sigma \circ \tau^{-1})$ is $\prod_{i=1}^L c_{\ell_i-1}(-1)^{\ell_i-1}$ where $\ell_1\etc \ell_L$ are the sizes of the cycles of $\sigma \circ \tau^{-1}$ by \cite[Equation (14)]{CS06}, and $\sharp(\sigma \circ \tau^{-1})$ is the number of cycles of $\sigma \circ \tau^{-1}$ (counting cycles of size one).
\par Choosing $\sigma$ and $\tau$ as in the sum above for which all the indicator functions are nonzero amount to cover the edges of $T$ by disjoint cycles alternating between edges labeled $x$ and $x^*$ as follow. Firstly, we think $\sigma$ (resp. $\tau$) as the map sending the edge $(i'_n,j'_n)$ to $(j_{\sigma(n)},i_{\sigma(n)})$ (resp. $(j_{\tau(n)},i_{\tau(n)})$). Then with this convention we denote by $\pi$ the permutation of $E$ defined by $\pi(e) = \sigma(e)$ if $e$ has label $x^*$ and $\pi(e) = \tau^{-1}(e)$ if $e$ has label $x$. Note that $\sigma \circ \tau^{-1}$ is $\pi^2$ restricted to the edges labeled $x^*$. In particular, the number of cycles of $\pi$ is the one of $\sigma^{-1} \circ \tau$, the lengths of the cycles of $\pi$ is two times those of $\sigma^{-1} \circ \tau$.

\par We introduce the undirected graph (with possibly multiple edges) $\mathcal G = (\mathcal V, \mathcal E)$, where $\mathcal V$ is the union of $V$ and of the set of simple cycles of $\pi$, and $\mathcal E$ is the multi-set where for each edge $e$ of $E$ we add in $\mathcal E$ edges between the goal of $e$ and the cycles of $\pi$ that contain $e$. Note that assuming that their exist $\sigma$ and $\tau$ such that in the associated term of \eqref{Weingarten} the indicator function is nonzero and constructing the partition $\pi$, we get that each edge belong at least to a directed simple cycle of $T$. 

Let us prove that $\mathcal G$ is connected. Let $v$ and $v'$ be two vertices of $\mathcal G$ that correspond to vertices of $V$. In $T$ there is a path from $v$ to $v'$. While looking at consecutive steps of the path that stay in a same simple cycle $c_i$, the vertices that are visited are all linked to $c_i$ in $\mathcal G$. Now when from a step to another we move to a different cycle $c_{i+1}$, the vertex between these two steps belong both to the cycles $c_i$ and $c_{i+1}$. Since each vertex of $\mathcal G$ that corresponds to a simple cycle of $T$ are linked to some vertices of $T$ in $\mathcal G$, we get that $\mathcal G$ is connected.

 One has $|\mathcal V| = |V| + \sharp(\pi)$ and $|\mathcal E| = |E|$. We then get from the formulas for $\tau_N^0$, $\delta_N^0$, $\textrm{Wg}( \cdot , N)$ above and Lemma \ref{lem:EdgesVertices}, 
	\eqa
		\tau_N^0\big[T(U_N)\big] = \sum_\pi N^{|\mathcal V|-1-|\mathcal E|} \times \phi(\pi)= \sum_{\pi} \big( \one_{\mathcal G \textrm{ is a tree}} + O(N^{-1}) \big) \phi(\pi) \label{ProofDistrHaar},
	\qea
where $\phi(\pi)= \prod_{i=1}^L c_{k_i/2-1}(-1)^{k_i/2-1}$ for $k_1\etc k_L$ are the sizes of the cycles of $\pi$. 

When $\mathcal G$ is a tree we claim that each edge of $T$ belongs exactly to one simple cycle which is directed (we then call $T$ a \emph{cactus}). Indeed, assume that a vertex $(v,w)$ belongs to two distinct simple cycles $c_1$ and $c_2$. Then there is a link between $w$ and $c_1$ and a link between $w$ and $c_2$. But $v$ has the same property since it is both the target of an edge in $c_1$ and an edge in $c_2$, in contradiction with the fact that $\mathcal G$ is a tree.

Given that $T$ is a cactus, there is no choice in the partition $\pi$ in the sum of Equation \eqref{ProofDistrHaar}, its cycles must consist in the simple cycles of $T$. Indeed, otherwise there is a cycle of $\pi$ with two distinct edges with same target, and then $\mathcal G$ has an edge of multiplicity two and $\mathcal G$ is not a tree. We get the expected formula.

Let now prove the factorization property and let $T_1=(V_1,E_1) \etc T_n = (V_n,E_n)$ be $^*$-test graphs in one variable $x$. We first use Lemma \ref{lem:PropConvEsp}, 
\eq
	 \esp\Big[ \frac 1 N \Tr^0 \big[ T_1(U_N) \big] \dots  \frac 1 N \Tr^0 \big[ T_n(U_N) \big] \Big]  
	  & = &  \sum_{\sigma}  \frac 1 {N^{n}}   \esp\Big[  \Tr^0 \big[ T^\sigma(U_N) \big]  \Big],
\qe
where the sum is over the partition $\sigma$ that possibly identify vertices of different graphs $T_i$. Fix $\sigma$ as in the sum. Computing $\esp\Big[ \frac 1 N \Tr^0 \big[ T^\sigma(U_N) \big]  \Big]$, we reason as in the proof of Proposition \ref{Prop:HaarDist}. The difference is that now $T^\sigma$ is not connected in general. Hence, we now deduce as in \eqref{ProofDistrHaar} that
	\eq
		\frac 1 {N^{n}} \tau_N^0\big[T^\sigma(U_N) \big] = \sum_\pi\big( \one_{\mathcal G \textrm{ is a forest}} + O(N^{-1})\big) \phi(\pi) \times N^{n_\pi-n},
	\qe
where $n_\pi$ is the number of connected components of $T^\sigma$. Hence, this quantity goes to zero unless the partition $\sigma$ is trivial. Then, the sum over the choices for a covering of $T^\sigma$ by cycles splits into $n$ sums for the coverings of the $T_j$'s by cycles. As the map $\phi(\pi)$ is multiplicative with respect to its cycles, we get as expected
	\eq
		 \esp \Big[ \prod_{\ell=1}^n \frac 1 N \Tr^0\big[T_\ell(U_N) \big] \Big]  =  \prod_{\ell=1}^n \esp \Big[ \frac 1 N \Tr^0\big[T_\ell(U_N) \big] \Big] + O(N^{-1}).
	\qe
\end{proof}

\section{Diagonal matrices and the factorization property} 

\subsection{Diagonal matrices}
In the section we first consider diagonal matrices in the context of Theorem \ref{MainTh}. 

\begin{lemma}\label{Lem:DiagMatrices} Let $\mathbf D_N=(D_j)_{j\in J}$ be a family of random diagonal matrices. Then $\mathbf D_N$ converges in traffic distribution whenever it converges in $^*$-distribution: for any $^*$-test graph $T=(V,E,\gamma, \varepsilon)$, one has
	\eq
		\tau_N\big[T(\mathbf D_N)\big]  =   \Phi_N\big[  \prod_{e\in E} D_{\gamma(e)}^{\varepsilon(e)} \big], \ \ \tau_N^0\big[T(\mathbf D_N)\big]   =  \one(|V|=1) \times \Phi_N\big[  \prod_{e\in E} D_{\gamma(e)}^{\varepsilon(e)} \big] .
	\qe
It satisfies the factorization property (Assumption B3 of Theorem \ref{MainTh}) if and only if for any $K\geq 2$ and any $^*$-polynomials $P_1\etc P_K$ one has
	$$\prod_{k=1}^K \Phi_N\Big[  P_k(\mathbf D_N) \Big] - \esp\Big[\prod_{k=1}^K \frac 1 N \Tr P_k(\mathbf D_N) \Big] \limN 0.$$
\end{lemma}

The proof follows immediately from the fact that non diagonal entries of $\mathbf D_N$ are zero.

\begin{example} Let $D_N$ be a diagonal matrix with independent and identically distributed diagonal entries, distributed according to a probability distribution $\mu$ on $\mathbb C$, whose moments of all orders exist. Then, the traffic distribution of $D_N$ does not depend on $N$. For any $^*$-test graph $T=(V,E,\gamma,\varepsilon)$ in one variable $x$, one has $\tau_N\big[T(D_N)\big] = \int z^n\bar z^m\mathrm d \mu (z)$ where $n$ and $m$ are respectively the number of edges labeled $x$ and $x^*$. Hence $D_N$ converges in traffic distribution. It is immediate to see that it also satisfies the factorization property. 
\end{example}

We have now all the ingredients for the proof of Corollary \ref{Cor:ApplClass}.

\begin{proof}[Proof of Corollary \ref{Cor:ApplClass}] Let  $\mathbf X_N, \mathbf U_N , \mathbf V_N, \mathbf D_N$ be independent families of matrices as in the Corollary, namely:
\begin{enumerate} 
	\item $\mathbf X_N$ is a family of independent Wigner matrices whose entries are invariant in law by complex conjugation, for which the convergence in traffic distribution is proved in Proposition \ref{Prop:WigDistr} and the factorization property is proved at the end of Section \ref{Sec:ExamplesMatrices}.
	\item $\mathbf U_N $ and $ \mathbf V_N$ are respectively a family of independent unitary Haar and uniform permutation matrices. Their convergence in traffic distribution and the factorization property are proved in Section \ref{Sec:ExamplesMatrices}
	\item $\mathbf D_N=(D_j)_{j\in J}$ is a family of independent diagonal matrices with independent and identically distributed diagonal entries whose moments of all orders exist. By Lemma \ref{Lem:DiagMatrices} above and its example, the matrices of $\mathbf D_N$ also satisfy these assumptions.
	\end{enumerate}
The matrices of $\mathbf X_N, \mathbf U_N , \mathbf V_N, \mathbf D_N$ are permutation invariant. Hence by Theorem \ref{MainTh}, for any family of matrices $\mathbf A_N$ converging in traffic distribution and satisfying the factorization property, the matrices of $\mathbf X_N, \mathbf U_N , \mathbf V_N, \mathbf D_N$ and the family $\mathbf A_N$ are asymptotically traffic independent.

The limiting $^*$-distribution of $\mathbf X_N, \mathbf U_N , \mathbf V_N, \mathbf D_N, \mathbf A_N$ depends only of the limit for each family of normalized trace of $^*$-test graphs $T$ such that there exists a cycle visiting each edge of $T$ is the sense of their orientation. Hence the limiting $^*$-distribution is unchanged if the Wigner matrices are replaced by unitarily invariant Gaussian Wigner matrices with parameters of the form $(\alpha, 0)$. By Theorem \ref{Th:AsyFree}, the matrices of $\mathbf X_N \cup \mathbf U_N$ are asymptotic free independent and $\mathbf X_N \cup \mathbf U_N$ is asymptotically free independent from $\mathbf V_N \cup \mathbf D_N \cup \mathbf A_N$. Hence the corollary.
\end{proof}

\subsection{On the factorization property}

 The convergence of the traffic distribution of $\mathbf A_N=\mathbf A_N^{(1)} \cup \dots \cup \mathbf A_N^{(L)}$ remains true if we do not assume the factorization property but instead the convergence of $(g_1 \etc g_K) \mapsto \esp\big[ \prod_k \frac 1 N \Tr[ g_k(\mathbf A_N^{(\ell)})]$. See Remark \ref{Rk:Rk:MainThBack} in the proof of Theorem \ref{MainTh}. When the factorization property is not satisfied, the limiting traffic distribution of $\mathbf A_N$ is not the product of the limiting distribution distributions of the $\mathbf A_N^{(\ell)}$ in the sense of Definition \ref{def:FreeProdGraphs}, but this is the convention in \cite{GAB15}.

		 Let $\mathbf A_N$ and $ \mathbf B_N$ be two independent families of matrices, and assume that the matrices of $\mathbf A_N$ are diagonal. Assume that $\mathbf A_N $ and $\mathbf B_N $ converges in traffic distribution, that one of the family is permutation invariant, and that only the family of diagonal matrices $\mathbf A_N$ satisfies the factorization property. Then the conclusion of Theorem \ref{MainTh} is valid for $\mathbf A_N, \mathbf B_N$, even if $\mathbf B_N$ does not necessarily satisfies the factorization property. See Remark \ref{Rk:Rk:MainThBack} in the proof of Theorem \ref{MainTh}.
		
		 Let $A_N$ be a random matrix which is a uniform permutation matrix with probability one half and a unitary Haar matrix otherwise. Since independent uniform permutation and unitary Haar matrices have same limiting $^*$-distribution and are asymptotically free independent, independent copies of $A_N$ are asymptotically free independent. Nevertheless we know from Propositions \ref{Prop:PermDistr} and \ref{Prop:HaarDist} that a uniform permutation matrix and a unitary Haar matrix do not have the same limiting traffic distribution. It follows that $A_N$ does not satisfies the factorization property \eqref{eq:DecorrIntro} (this is a simple exercise left to the reader). So according to the first point of the remark, independent copies of $A_N$ are not asymptotically traffic independent. Let now consider a random matrix $B_N$ of the form of a block matrix
			$$B_N = \tilde V_N \left( \begin{array}{cc} U_{p_N} & 0 \\ 0 & V_{N-p_N} \end{array} \right)\tilde V_N^*,$$
where $V_{N-p_N}, \tilde V_N$ are independent uniform permutation matrices, independent from a unitary Haar matrix $U_{p_N}$, and $p_N$ a sequence of positive integers such that $\frac{p_N}N \limN \frac 1 2$. Then it is easy to see that $B_N$ has the same limiting traffic distribution as $A_N$ and it satisfies the factorization property \eqref{eq:DecorrIntro}, so independent copies of $B_N$ are asymptotically traffic independent.

\part{Traffics and their Independence}
\chapter{Algebraic Traffic Spaces}\label{Sec:TrafficSpaces}

In the previous sections we have considered the traffic distributions of matrices and their point-wise convergence. After a recall on free probability (see \cite{AGZ,NS} for detailed presentations), we define in this section the abstract \emph{traffics} which model the limits of large matrices for this mode of convergence. 

\section{Non commutative probability spaces}
\begin{definition}[Spaces of non commutative random variables]\label{Def:NCPS} A non commutative probability space (sometimes called \emph{algebraic} probability space) is a pair $(\mathcal A, \Phi)$, where 
\begin{itemize}
	\item $\mathcal A$ is a unital algebra over $\mathbb C$,
	\item $\Phi: \mathcal A \to \mathbb C$ is a unital linear functional. It is called a trace when it satisfies $\Phi(ab) = \Phi(ba)$ for any $a,b$ in $\mathcal A$.
\end{itemize}
A $^*$-probability space is a non commutative probability space $(\mathcal A, \Phi)$ such that
\begin{itemize}
	\item $\mathcal A$ is a $^*$-algebra, i.e. it is endowed with an anti linear involution $\cdot \, ^*$ such that $(ab)^*=b^*a^*$ for any $a,b\in \mathcal A$,
	\item $\Phi$ is a state, that is it satisfies the positivity condition $\Phi(a^*a)\geq 0$ for any $a\in \mathcal A$. It is called faithful if moreover $\Phi(a^*a)=0$ implies $a=0$ for any $a\in \mathcal A$.
	\end{itemize}
	
Elements of a non commutative probability space are called non commutative random variables. The non commutative distribution of a family $\mathbf a=(a_j)_{j\in J}$ of non commutative random variables is the linear map $\Phi_{\mathbf a}$, defined on polynomials in indeterminates $\mathbf x=(x_j)_{j\in J}$, by
	$$\Phi_{\mathbf a} : P \mapsto \Phi\big( P(\mathbf a)\big).$$
To emphasis the role of the linear form we sometimes say that $\Phi_{\mathbf a}$ is the distribution of $\mathbf a$ w.r.t. $\Phi$. For $\mathbf a$ in a $^*$-probability space, the $^*$-distribution of $\mathbf a$ is the non commutative distribution of $(\mathbf a, \mathbf a^*)$, or equivalently the map $\Phi_{\mathbf a}$ defined as above but for $^*$-polynomials.
The convergence of a family $\mathbf a_N$ of non commutative random variables is the point-wise convergence of their non commutative distribution.
\end{definition}

\begin{remark} In a $^*$-probability space $(\mathcal A, \Phi)$, the positivity of $\Phi$ yields the Cauchy-Schwarz inequality, that is $\Phi(ab)^2\leq \Phi(a^*a) \Phi(b^*b)$ for any $a,b\in \mathcal A$ and faithfulness of $\Phi$ implies that $(a,b)\mapsto \Phi(a^*b)$ is actually a scalar product.
\end{remark}

\begin{example}\label{Ex:Algebra}
\begin{enumerate}
	\item {\bf Classical random variables:} Let $(\Omega, \mathcal F, \mathbb P)$ be a probability space in the usual sense. The algebra $\mathcal L^{-\infty}(\Omega, \mathbb C)$ of measurable maps $\Omega \mapsto \mathbb C$ with finite moments of all orders is a $^*$-probability space, equipped with the complex conjugate and the expectation $\esp$. Its quotient $L^{-\infty}(\Omega, \mathbb C)$ by measurable maps null almost everywhere is a $^*$-probability space with faithful state. 
	\item {\bf Random matrices:} The space $\textrm M_N(\mathbf C)$ of deterministic $N$ by $N$ matrices is a $^*$-probability space with trace $ \frac 1 N \Tr $. Let consider an algebra $L^{-\infty}\big(\Omega, \textrm M_N(\mathbf C)\big)$ of random matrices whose entries, defined in a same probability space, have finite moments of all orders. It is a $^*$-probability space with faithful trace $\tau_N = \esp \big[ \frac 1 N \Tr \big]$.
\end{enumerate}
\end{example}

In the next section, we define the algebraic traffic spaces, that can be seen  as non commutative probability spaces with more structure. We do not define the associated notion of positivity in this article, which do not play an important role in the questions considered here (see \cite{CDM16}).

\section{Algebraic traffic spaces}\label{Sec:DefTrafficSpaces}

\subsection{Algebras over the operad of graph operations}

 The main point in defining algebraic traffic spaces is to formalize the good structure that replaces the notion of algebra. For that task, we use the notion of symmetric operads. An operad is a graded set $\mathcal G = \bigcup_{K\geq 0}\mathcal G_K$, where an element $g$ of $\mathcal G_K$ is seen as an operation which takes $K$ objects and gives a single one. Compatibility conditions are assumed in order to have a canonical notion of algebra over the operad, namely of vector space for which one can compose the elements according to the operations of $\mathcal G$.

\begin{definition} An operad is a set $\mathcal G=\bigcup_{K\geq 0}\mathcal G_K$ of sets endowed with operations of composition
	$$ (g,g_1 \etc g_K) \in \mathcal G_K \times \mathcal G_{L_1} \times \dots \times \mathcal G_{L_K} \mapsto g(g_1 \etc g_K)\in \mathcal G_{L_1 + \dots + L_K}$$ and a fixed element $id_\mathcal G$ in $\mathcal G_1$, called the identity of the operad, satisfying the following properties.
\begin{itemize}
	\item The element $id_\mathcal G$ is a \emph{unit for composition}, namely for any $g \in \mathcal G$ one has $g = id_\mathcal G( g) = g\big(id_\mathcal G \etc id_\mathcal G\big)$.
	\item The composition is \emph{associative}: for any $g \in \mathcal G_K$, any $g_k \in \mathcal G_{L_k}$ and any $g_{k,\ell} \in \mathcal G$, $k=1\etc K, \ell=1\etc L_k$, one has 
	\eq
		\lefteqn{g \big( g_1( g_{1,1} \etc g_{1,L_1}) \etc g_K( g_{K,1} \etc g_{K,L_K})\big)}\\
		& = & \big( g(g_1 \etc g_K) \big) ( g_{1,1} \etc g_{1,L_1} \etc g_{K,1} \etc g_{K,L_K}).
	\qe
\end{itemize}
A symmetric operad is an operad $\mathcal G =  \bigcup_{K\geq 0}\mathcal G_K$ endowed with an action $(\sigma, g)\mapsto g_\sigma$ of the symmetric group $S_K$ on $\mathcal G_K$ such that
	\eq
		g_\sigma (g_1 \etc g_K) = g(g_{\sigma(1)} \etc g_{\sigma(K)}),\\
		g\big( (g_1)_{\sigma_1} \etc (g_K)_{\sigma_K})=g(g_1 \etc g_K)_{ \sigma_1 \times \dots \times \sigma_K}.
	\qe
\end{definition}

To define the traffic spaces we use the following operad.

\begin{definition}[The operad of graph operations]~ \begin{itemize}
	\item A graph operation $g$ of $K$ elements is a finite connected oriented graph $(V,E)$ with $K$ edges, with the data of an ordering of the edges and of two vertices $in, out$, (possibly equal) called the input and the output. An ordering of the edges can be thought as a bijection $\gamma: E \to [K]$. We denote $g=(T, in, out)$ and $T=(V,E,\gamma)$. We set $\mathcal G_K$ the set of graph operations of $K$ elements and $\mathcal G=\bigcup_{K\geq 0} \mathcal G_K$.

	\item The graph operation of $\mathcal G_1$ with exactly two vertices $in\neq out$ and one single edge from $in $ to $out$ is denoted by its graph $(\cdot \leftarrow \cdot)$, where implicitly in this picture the input is on the right and the output on the left.
	\item As we have seen before for graph monomials, we define the operation consisting in the replacement of the edges of a graph operation $g$ of $\mathcal G_K$ by $K$ graphs $g_1\etc g_K$ of $\mathcal G$. The input (respectively the output) vertex of $g_k$ replaces the source (respectively the target) vertex of the $k$-th edge of $g$. The induced order for this new graph is the lexical order, namely the edges of $g_i$ come before edges of $g_{i+1}$ for $i=1\etc n-1$. This produces an new graph operation $g(g_1\etc g_K)$ in $\mathcal G$
	\item For any permutation $\sigma$ of $\{1\etc K\}$ and any $g=(T, in , out)\in \mathcal G_K$, $T=(V,E,\gamma)$, denote $g_\sigma$ the graph operation obtained from $g$ by considering the $i$-th edge of $g$ as the $\sigma(i)$-th one for $g_\sigma$, namely $g_\sigma = (\tilde T, in ,out)$, $\tilde T = (V,E,\sigma^{-1} \circ\gamma  )$.
\end{itemize}
\end{definition}
 \noindent We let the reader verify that $\mathcal G$ is a symmetric operad with identity $id_\mathcal G=( \cdot \leftarrow  \cdot)$.

\begin{definition}[Algebras over $\mathcal G$]~\label{Def:GAlgebra} \begin{enumerate}
	\item An algebra over the operad $\mathcal G$ of graph operations (in short a $\mathcal G$-algebra) is a vector space $\mathcal A$ over $\mathbb C$ endowed with an action of $\mathcal G$ as follow. 
	\begin{itemize}
		\item {\bf Linearity:} For any $K\geq 0$, $Z_g:\mathcal A^{\otimes K} \to \mathcal A$ is a linear map.
		\item {\bf Unity:} By convention, for the single graph operation $(\cdot)\in \mathcal G_0$ with one vertex and no edge, $Z_{(\cdot)}$ is a fixed element $\mathbb I$ of $\mathcal A$.
	 	\item {\bf Identity:} The identity of the operad $id_\mathcal G=(\cdot \leftarrow \cdot)$ is associated to the identity map, that is $Z_{id_\mathcal G} = id_\mathcal A$ or equivalently $Z_{(\cdot \leftarrow \cdot)}(a) = a$ for any $a\in \mathcal A$.
	 
		\item  {\bf Equivariance:} An action of a graph operation consists in replacing edges by elements of $\mathcal A$, so it depends on the locations of the edges not on their ordering: for any $a_1\etc a_K \in \mathcal A$, one has
		$$ Z_{g_\pi}(a_1 \otimes \dots \otimes a_K)=Z_g(a_{\pi(1)}  \otimes \dots \otimes a_{\pi(K)}).$$

		\item {\bf Substitution:} The action of graph operations on $\mathcal A$ are compatible with the substitution of graph operations: for any $g\in \mathcal G_K$, $g_1\etc g_K\in \mathcal G$,
		$$Z_g\big(Z_{g_1}  \otimes \dots \otimes Z_{g_K} \big) = Z_{g(g_1 \etc g_K)}.$$

\end{itemize}

	\item A $^*$-algebra over the operad $\mathcal G$ (in short, a $\mathcal G^*$-algebra) is a $\mathcal G$-algebra endowed with an antilinear involution $a \mapsto a^*$ with the following property. Given a graph operation $g$, we call transpose of $g$ and denote $g^t$ the graph operation obtained by reversing the orientation of the edges and interverting input and output. Then for any graph operation $g$ and any $a_1 \etc a_n$ in $\mathcal A$, one has $\big( Z_g(a_1 \otimes \dots \otimes  a_n) \big)^* = Z_{g^t}(a_1^* \otimes \dots \otimes  a_n^*)$.
	
		\item The $\mathcal G$-algebra spanned by a subset $A \subset \mathcal A$ of a $\mathcal G$-algebra is the linear space generated by $\{Z_g( a_1\otimes \dots \otimes a_K), g\in \mathcal G_K, a_1 \etc a_K \in A\}$. The $\mathcal G^*$-algebra spanned a subset $A \subset \mathcal A$ of a $\mathcal G^*$-algebra is defined similarly with $a_1\etc a_K\in A\cup A^*$.

\end{enumerate}
\end{definition}

\begin{example}\label{Ex:GAlgebra}
\begin{enumerate}
	\item {\bf Abelian algebras:} Let $\mathcal A$ be an abelian unital algebra with product $(a,b)\mapsto a\times b$ and unit $\mathbb I$. There is a trivial structure of $\mathcal G$-algebra on $\mathcal A$: for any family $\mathbf a= (a_1 \etc a_K)$ of elements of $\mathcal A$ and any graph operation $g\in \mathcal G_K$, we set $Z_g(a_1 \otimes \dots \otimes a_K) =a_1 \times \dots \times a_K$. This defines well a structure of $\mathcal G$-algebra on $\mathcal A$. If moreover $\mathcal A$ is a $^*$-algebra with involution $\cdot ^*$ then $\mathcal A$ is a $\mathcal G^*$-algebra.	
	 Note that if $\mathcal A$ is not abelian, then $Z_g(a_1\otimes \dots \otimes a_K)$ can be defined similarly but depends on the ordering of its edges of $g$. Hence the equivariance axiom is not verified.
	\item {\bf Matrix algebras:} The space of $N \times N$ matrices with complex entries is a $\mathcal G^*$-algebra for the operations given with similar formula as Definition \ref{Def:StarGraphsMon}: for any matrices $A_1\etc A_K$ and any graph operation $g=(T, in, out)$ with $T=(V,E,\gamma)$ and $|E|=K$, the entry $(i,j)$ of $Z_g(A_1\otimes \dots \otimes  A_K)$ is
		$$Z_g(A_1\otimes \dots \otimes  A_K)(i,j) = \sum_{\substack{ \phi:V \to [N] \\ \phi(i)=out, \phi(j)=in}} \prod_{k=1}^K A_{k}\big( \phi(w_k), \phi(v_k)\big),$$
		where the $k$-th edge is denoted $(v_k,w_k)$.
	\end{enumerate}
\end{example}

\subsection{Polynomials and graph polynomials}

In a $\mathcal G$-algebra, we define a product $\times$ as the bilinear map $(a,b) \mapsto Z_{(\cdot \overset{1} \leftarrow \cdot \overset{2} \leftarrow \cdot)}(a\otimes b)$, where the graph operation consists in two consecutive edges $e_1=(v,out)$ and $e_2=(in,v)$ with $out,v, in $ pairwise distinct. The product is associative since by the axiom of substitution one has $(a \times b)\times c = Z_{(\cdot \overset{1} \leftarrow \cdot \overset{2} \leftarrow \cdot\overset{3} \leftarrow \cdot)}(a,b,c) =a \times (b\times c)$. This defines a structure of algebra on $\mathcal G$-algebras for which $\mathbb I=Z_{(\cdot)}$ is the unit and we simply denote $ab=a\times b$. Moreover, in a $\mathcal G^*$-algebras, we have by definition $(a b)^*  = Z_{(\cdot \overset{1} \leftarrow \cdot \overset{2} \leftarrow \cdot)^t}(a^*\otimes b^*) =Z_{(\cdot \overset{2} \leftarrow \cdot \overset{1} \leftarrow \cdot)}(a^*\otimes b^*)  =  b^*a^*$. Hence a $\mathcal G^*$-algebra is in particular a $^*$-algebra. 
\\
\par Recall that a graph monomial is a finite connected directed graph whose edges are labeled by formal variables with an input and an output, and that $\mathbb C \mathcal G\langle \mathbf x \rangle$ denotes their linear space. Note that contrary to graph operations, in a graph monomial $g=(V,E,\gamma)$ we do not consider an ordering of the edges but a labeling $\gamma:E\to J$, and a same label can appear on several edges. The space of graph polynomials in some variables is a $\mathcal G$-algebra: for $g$ a graph operation with $K$-edges and $g_1\etc g_K$ graph monomials, the graph $Z_g(g_1\otimes \dots \otimes g_K)$ is given by substitution of edges as usual.

 Let $\mathbf a=(a_j)_{j\in J}$ be a family of elements of a $\mathcal G$-algebra $\mathcal A$ and let $g$ be a graph monomial in the variables $\mathbf x=(x_j)_{j\in J}$ with $K$ edges. Choosing an arbitrary ordering of its edges gives a graph operation that we denote $\tilde g$ and allows us to consider the element $Z_{\tilde g}(a_{j_1} \otimes \dots \otimes a_{j_K})$ of $\mathcal A$, where $j_k \in J$ is the label of the $k$-th edge of $\tilde g$. By the equivariance axiom for the action of operads, this element depends only on $g$ and $\mathbf a$, not on the ordering of the edges. It is denoted $g(\mathbf a)$ and this definition is extended by linearity for $g \in \mathbb C\mathcal G\langle \mathbf x\rangle$.

By definition of the product and thanks to the axiom of substitution on $\mathcal G$-algebras, the subspace of $\mathbb C\mathcal G\langle \mathbf x\rangle$ generated by simple directed lines $(\underset{out}\cdot \overset{x_{j_1}} \leftarrow \cdot \dots \cdot \overset{x_{j_n}} \leftarrow \underset{in}\cdot)$, $j_k\in J$, can be identified with the space of non commutative polynomials $\mathbb C\langle \mathbf x\rangle$.

Similarly, recall that a $^*$-graph monomial in variables $\mathbf x$ is a graph monomial with edges labeled by variables $\mathbf x$ and $\mathbf x^*$. With same notations we define $g(\mathbf a) = Z_{\tilde g}(a_{j_1}^{\varepsilon_1} \otimes \dots \otimes a_{j_K}^{\varepsilon_K})$ where $\varepsilon_k$ is $1$ or $*$ depending if the $k$-edge is labeled by a variable in $\mathbf x$ or $\mathbf x^*$ respectively.

The axioms for the action of graph operations implies similar properties for graph polynomials. 
\begin{enumerate}
	\item The {\bf substitution axiom} implies that for any $g, h_j, j\in J$, elements of $\mathcal G \langle \mathbf x \rangle$ and any $\mathbf a \in \mathcal A^J$, the element $g\big(h_j(\mathbf a), j \in J\big)$ is equal to $ \big(g(h_j)_{j\in J}\big)(\mathbf a)$ where $g(h_j)_{j\in J}$ is obtained from $g$ by replacing each edge labeled $x_j$ by the graph $h_j$. 
	\item The linearity of the $Z_g$'s implies that a graph monomial is {\bf multi-linear with respect to its edges}, in the following sense. Let $g$ be a graph monomial in variables $\mathbf x$ and $y$, let $\mathbf a \in \mathcal A^J$, $b_1,b_2 \in \mathcal A$ and $\lambda \in \mathbb C$. Denote by $E_y$ the set of edges of $g$ labeled $y$. Then one has $g(\mathbf a, \lambda b) = \lambda^{|E_y|}g(\mathbf a,b)$ and 
		$$g(\mathbf a,b_1+b_2) = \sum_{\gamma:E_y\to\{1,2\}} g_\gamma(\mathbf a, b_1,b_2),$$
 where $g_\gamma$ is obtained from $g$ by declaring that an edge $e\in E_y$ has label associated to $b_{\gamma(e)}$.
\end{enumerate}

For elements $a,b$ of a $\mathcal G$-algebra, we define $a^t$, $\Delta(a)$, $deg(a)$ and $a\circ b$ with the same graph monomials as for matrices in Example \ref{Ex:GraphPolyMatrices}. In any $\mathcal G$-algebras, there are infinitely many relations between the operations. For instance, the relation 
	\eqa \label{Eq:ExRelation}
		\Delta\big( a \times deg(b)\big) = deg\big( \Delta(a) \times b\big),
	\qea
 valid for any $a,b\in \mathcal A$, was obtained in Figure \ref{fig:Subs} in Chapter 1.

\subsection{Algebraic traffic spaces}

Defining the traffic distribution of a family of matrices $\mathbf A_N$, we have firstly considered the map $\bar\Phi_{\mathbf A_N} : \mathbb C \mathcal G\lara \to \mathbb C$ in Definition \ref{Def:StarGraphsMon} which generalizes the notion of $^*$-distribution. Then we have introduced the equivalent data given by the combinatorial distribution $\tau_{\mathbf A_N} : \mathbb C \mathcal T\lara \to \mathbb C$ in Definition \ref{Def:CombiDistrTraf}, from which is defined the injective distribution $\tau_{\mathbf A_N}^0$.  In the definition of abstract traffics below (Definition \ref{Def:TrafficSpace}), we consider the data of the combinatorial distribution to be the intrinsic one. For this reason, the $\mathcal G$-algebras considered to introduce the traffic spaces come together with the following space of observables.

\begin{definition}
 Let $\mathcal A$ be an arbitrary set. A test graph labeled in $\mathcal A$ is a triplet $T=(V,E,\gamma)$ where $(V,E)$ is a finite connected graph, endowed with a map $\gamma:E\to \mathcal A$. We define by $\mathcal T\langle \mathcal A \rangle$ the set of test graphs labeled in $\mathcal A$ and by $\mathbb C \mathcal T \langle \mathcal A \rangle$ its linear space. 
\end{definition}

A test graph $ T$ in variables $\mathbf x=(x_j)_{j\in J}$ and a family $\mathbf a=(a_j)_{j\in J}$ of elements of $\mathcal A$ defines obviously an element $T(\mathbf a) \in \mathcal T\langle \mathcal A \rangle$.

\begin{definition}
 \label{Def:TrafficSpace}[Algebraic traffic spaces] An algebraic traffic space is a pair $(\mathcal A, \tau)$, where $\mathcal A$ is a $\mathcal G$-algebra and $\tau : \mathbb C \mathcal T\langle \mathcal A \rangle \to \mathbb C$ is a linear form satisfying the following properties.
\begin{enumerate}
	\item {\bf Unity:} $\tau$ sends the graph with no edges to one.
	\item {\bf Substitution:} For any $T \in \mathcal T\langle \mathcal A \rangle$ having an edge $e_0$ labeled $h(\mathbf a)$ for a graph monomial $h$, then $\tau\big[T\big] = \tau\big[ T_h\big]$ where $T_h$ is obtained from $T$ by replacing the edge $e_0$ by the graph $h$, with labels as in the evaluation $h(\mathbf a)$. 
	\item {\bf Multi-linearity w.r.t. the edges:} For any $T \in \mathcal T\langle \mathcal A \rangle$ having an edge $e_0$ labeled $a+\lambda b$, we have $\tau[T] = \tau[T_a]+ \lambda \tau[T_b]$ where $T_a$ and $T_b$ are obtained from $T$ by declaring that the label of $e_0$ is now $a$ and $b$ respectively.
\end{enumerate}

Elements of $\mathcal A$ are called \emph{traffics} and $\tau$ is called the combinatorial trace on $\mathcal A$. The \emph{traffic distribution} of a family $\mathbf a=(a_j)_{j\in J}$ of elements of $\mathcal A$ is the data of the map 
	$$\tau_{\mathbf a} : T \in \mathbb C\mathcal T\langle \mathbf x \rangle \mapsto \tau\big[T(\mathbf a)\big] \in \mathbb C.$$
	As we will see equivalent formulation of the traffic distribution, we name specifically $\tau_{\mathbf a} $ the \emph{combinatorial distribution} of the family $\mathbf a$. Let $\mathbf a_N=(a_{N,j})_{j\in J}$, for each $N\geq1$, and $\mathbf a=(a_j)_{j\in J}$ be families of traffics, possibly in different spaces for different $N$. We say that $\mathbf a_N$ converges to $\mathbf a$ in traffic distribution as $N$ goes to infinity whenever $ \tau_{\mathbf a_N}$ converges point-wise to $ \tau_{\mathbf a}$.

\end{definition}

A motivation for using the term \emph{traffics} is because of the algebraic structure of $\mathcal G$-algebra which allows to compose them not only by multiplication, but thanks to schemes given by graph monomials. The term \emph{traffic space} means space of traffics as for the term vector space.

\begin{example}\label{Ex:TrafficSpaces} Example \ref{Ex:Algebra} continued.
\begin{itemize}
	\item Let $(\mathcal A, \Phi)$ be an abelian non commutative probability space. We endow $\mathcal A$ with its trivial structure of $\mathcal G$-algebra structure. Moreover we define $\tau:\mathbb C \mathcal T\langle \mathcal A \rangle \to \mathbb C$ by $\tau[ T ]\mapsto  \Phi\big[ \prod_e a_e \big]$ where the product is over the edges $e$ of $T$ and $a_e$ denotes the label of $e$. Hence $\tau$ is simply the expectation of the product of the labels of $T$. Then $(\mathcal A, \tau)$ is an algebraic traffic space.
	\item The $\mathcal G$-algebra of $N$ by $N$ matrices is an algebraic traffic space when endowed with the combinatorial trace of test graphs in matrices defined in the first part of the article. Note yet that the traffic distribution of a family of matrices $\mathbf A_N$ as defined in the first part (defined on $^*$-test graphs) corresponds to the traffic distribution of $(\mathbf A_N, \mathbf A_N^*)$ in the sense of algebraic traffic spaces.
\end{itemize}
\end{example}

\begin{lemma} For each $N\geq 1$, let $\mathbf a_N=(a_{N,j})_{j \in J}$ be a family of traffics in some space $(\mathcal A_N, \tau_N)$. Assume that $\tau_{\mathbf a_N}$ converges point-wise. Then there exists a family $\mathbf a = (a_j)_{j\in J}$ in some space $(\mathcal A, \tau)$ such that $\mathbf a_N$ converges to $\mathbf a$ in traffic distribution.
\end{lemma}

In particular, for a family of matrices $(A_{N,j})_{j\in J}$ converging in traffic distribution in the sense of the first part of the article, the family $(A_{N,j}, A_{N,j}^*)_{j\in J}$ converges to a family $(a_j, a_j^*)_{j\in J}$ in an algebraic traffic space.

\begin{proof} Let $\tau:\mathbb C \mathcal G\langle \mathbf x \rangle \to \mathbb C$ be the limiting distribution of $\mathbf a_N$, where $\mathbf x = (x_j)_{j\in J}$. The $\mathcal G$-algebra $\mathcal A:=\mathbb C\mathcal G\langle \mathbf x \rangle$ of graph polynomials in variables $\mathbf x$ endowed with $\tau$ is a traffic space and $\mathbf a:=\big( (\cdot \overset {x_j}\leftarrow \cdot)\big)_{j \in J}$ is limit in traffic distribution of $\mathbf a_N$.
\end{proof}

Similarly if $\mathbf a_N$ is a family of traffics in a $\mathcal G^*$-algebra and $\tau_{\mathbf a_N, \mathbf a_N^*}$ converges point-wise, then $(\mathbf a_N, \mathbf a_N^*)$ has a limit in traffic distribution $(\mathbf a, \mathbf a^*)$ in some algebraic traffic space which is a $\mathcal G^*$-algebra.

\subsection{Trace, anti-trace and injective trace}

A traffic can be seen as a non commutative random variable, but with two different points of view to do so, as the following proposition shows.

\begin{proposition}\label{Def:TrafficSpaceOtherDef} Let $(\mathcal A, \tau)$ be an algebraic traffic space. We define two functions $\Phi$ and $\Psi$ from $\mathcal A$ to $\mathbb C$ by
	$$\Phi(a) =  \tau[^a\circlearrowleft], \ \ \Psi(a) = \tau\big[ \cdot \overset{ a } \leftarrow \cdot\big], \ \forall a \in \mathcal A,$$
namely the combinatorial trace of a self-loop and a simple edge respectively.
Then $\Phi$ and $\Psi$ are unital linear forms on $\mathcal A$, defining two structures of non commutative probability space. The form $\Phi$ is tracial and is called the \emph{trace} associated to $\tau$. The linear form $\Psi$ is called the \emph{anti-trace} associated to $\tau$ (it is not a trace in general).
\end{proposition}

\begin{proof} We have $\Phi(\mathbb I) =  \tau[^{\mathbb I}\circlearrowleft]$ and by convention, $\mathbb I=Z_{(\cdot)}$ is associated to the graph monomial with no edges, so the substitution axiom implies $\tau[^{\mathbb I}\circlearrowleft] = \tau[\, \cdot \,]$. By the unity axiom for $\tau$, $\Phi$ is unital. Moreover, $\Phi$ is linear since $\Phi(a + \lambda b) = \tau[^{a + \lambda b}\circlearrowleft] = \tau[^{a}\circlearrowleft] + \lambda  \tau[^{b}\circlearrowleft]$ by the linearity of $\tau$ with respect to the edges of graphs. The same reasoning implies that $\Psi$ is linear. It remains to prove that $\Phi$ is tracial. But we have $\Phi(ab) = \tau[^{ab}\circlearrowleft] = \tau\big[ \cdot \overset a {\underset b\leftrightarrows} \cdot\big]$ by the substitution axiom. The expression is symmetric in $a$ and $b$, which yields $\Phi(ab) = \Phi(ba)$ as expected. Concrete examples where $\Psi$ is not a trace appear for instance in Section \ref{Sec:ExBool} below.\end{proof}

Let us first consider the trace $\Phi$ and find a characterization of this map in $(\mathcal A, \tau)$.

\begin{lemma} Let $(\mathcal A, \tau)$ be an algebraic traffic space. The trace $\Phi$ associated to $\tau$ satisfies the following properties.
\begin{enumerate}
		\item {\bf Diagonality:} Recalling that $\Delta$ is the graph monomial with a single vertex and a single self-loop, one has $ \Phi = \Phi \big( \Delta(\, \cdot \,) \big).$
		\item {\bf Input-independence:} For any graph monomial $g=(T,in,in)$ with same input and output and for any family $\mathbf a$ of elements of $\mathcal A$, $\Phi\big(g(\mathbf a)\big)$ does not depend on the place of input in $g(\mathbf a)$ but only on the test graph $T(\mathbf a)$ labeled in $\mathcal A$.
\end{enumerate}
Reciprocally, in a $\mathcal G$-algebra $\mathcal A$, assume that $\Phi$ is a unital linear map satisfying the two above properties. We define the linear form $\tau : \mathbb C \mathcal T\langle \mathcal A \rangle \to \mathbb C$ by $\tau\big[ T(\mathbf a)\big] = \Phi\big(g(\mathbf a)\big)$ for any $g=(T,in,in)$. Then $(\mathcal A, \tau)$ is an algebraic traffic space for which $\Phi$ is the associated trace.
\end{lemma}

\begin{proof}
For any $a\in \mathcal A$, we have $\Phi\big( \Delta(a) \big) =  \tau[^{\Delta(a)}\circlearrowleft]$. Using the substitution axiom for $\tau$ means replacing a self-loop by a self-loop, which yields $\Phi\big( \Delta(a) \big)  =\tau[^{ a}\circlearrowleft] = \Phi(a)$. Hence $\Phi$ is diagonal. Let $g=(T,in, in)$ be a graph monomial with same input and output. Then $\Phi\big( g(\mathbf a) \big) = \tau[^{g(a)}\circlearrowleft] = \tau\big[ T(\mathbf a)\big]$ by substitution axiom, which then depends only on the labeled graph $T(\mathbf a)$.

 Let now $\Phi$ be a unital linear form on a $\mathcal G$-algebra satisfying the two properties of the lemma, namely the diagonality and input-independence, and let $\tau$ as defined therein. Since $\Phi$ is unital and by definition of the action of the graph $g_0$ with no edges, we have $1=\Phi(\mathbb I) = \Phi(g_0) = \tau[ \, \cdot \,]$. Hence $\tau$ satisfies the unity axiom. Let now $T$ be a test graph labeled in $\mathcal A$ and write $\tau [ T  ] =  \Phi\big( g(\mathbf a) \big)$ for a graph polynomial obtained from $T$ by choosing an arbitrary vertex as common value of input and output. Assume that $T$ has an edge labeled $h(\tilde {\mathbf a})$ for a graph monomial $h$. By the substitution property for graph monomials stated in the previous section, $g(\mathbf a) =  g_h(\tilde {\mathbf a}, \mathbf a)$, where $g_h$ is obtained from $g$ by replacing the edge evaluated in $h(\tilde{\mathbf a})$ by the graph $h$. And then we get $\tau[T] = \tau[T_h]$ as expected in the Substitution axiom. The property of multi-linearity w.r.t. the edges of $\tau$ follows similarly from the same property for graph monomials stated in Definition \ref{Def:TrafficSpace}.\end{proof}

We now consider the anti-trace. We omit the proof of the lemma below that follows with same arguments as in proof of the analogue statement for $\Phi$.

\begin{lemma} Let $(\mathcal A, \tau)$ be an algebraic traffic space. The anti-trace $\Psi$ associated to $\tau$ satisfies the following {\bf input/output-independence} property. For any graph monomial $g=(T,in,out)$ (with arbitrary input and output) and for any family $\mathbf a$ of elements of $\mathcal A$, $\Psi\big(g(\mathbf a)\big)$ does not depend on the place of the input and the output but only on the test graph $T(\mathbf a)$ labeled in $\mathcal A$.

Reciprocally, in a $\mathcal G$-algebra $\mathcal A$, assume that $\Psi$ is a unital linear map satisfying the above input/output-independence property. Then the linear form $\tau : \mathbb C \mathcal T\langle \mathcal A \rangle \to \mathbb C$ defined by $\tau\big[ T(\mathbf a)\big] = \Phi\big(g(\mathbf a)\big)$ for any $g=(T,in, out)$ defines on $\mathcal A$ a structure of algebraic traffic space, for which $\Psi(a) = \tau[ \cdot \overset a \leftarrow \cdot ]$.
\end{lemma}
 We can also relate directly $\Phi$ and $\Psi$ as follow. Recall that $deg$ denotes the graph monomial with two distinct vertices $in=out$ and $v$ and an edge from $v$ to $in=out$. 

\begin{lemma} Let $(\mathcal A, \tau)$ be an algebraic traffic space with $\Phi$ and $\Psi$ defined as above. Then we have $\Psi(a) = \Phi\big( deg(a) \big)$ and $\Phi(a) =\Psi\big( \Delta(a) \big)$. 
\end{lemma}

\begin{proof} We have $\Phi\big( deg(a) \big) =\tau[^{ deg(a)}\circlearrowleft] $. Using the substitution axiom for $\tau$ means replacing the self-loop labeled $deg(a)$ by a simple edge, identifying the vertex of the loop with the vertex $in=out$ of the graph of $deg$. Hence $\Phi\big( deg(a) \big) =\tau[\cdot \overset a \leftarrow \cdot]  = \Psi(a).$ On the other hand, $\Psi\big(\Delta(a)\big) = \tau[\cdot \overset {\Delta(a)} \leftarrow \cdot]  = \tau[^{a}\circlearrowleft]  = \Phi(a)$.
\end{proof}

\begin{lemma}\label{Lem:TraceAdj} Let $(\mathcal A, \tau)$ be an algebraic traffic space with trace $\Phi$. Assume moreover that $\mathcal A$ is a $\mathcal G^*$-algebra and that $\Phi$ is a state, that is $\Phi(a^*a)\geq 0$ for any $a\in \mathcal A$. Then $\Psi$ is a state. Moreover, for a test graph $T\in \mathcal T\langle \mathcal A \rangle$, denote $T^*$ the test graph in same variables obtained by exchanging the orientations of each edge and replacing each label $a$ by its adjoint $a^*$. Then for any $T\in \mathcal T\langle \mathcal A \rangle$, one has $\tau[T^*] = \overline{ \tau[T]}$.
\end{lemma}

\begin{proof} If $\Phi$ is a state, then $\Psi$ is a state since 
	$$\Psi(a^*a) = \tau[\cdot \overset{a^*a}\leftarrow \cdot ] =  \tau[\cdot \overset{a^*}\leftarrow \cdot \overset{a^*}\leftarrow \cdot] = \Phi\big[ deg(a)^* deg(a) \big] \geq 0.$$
Moreover, if $\Phi$ is a state recall that $\Phi(a^*) = \overline{ \Phi(a)}$ \cite[Remarks 1.2]{NS}. Hence, given $\tau\big[ T(\mathbf a)\big] = \Phi\big( g(\mathbf a)\big)$, one verifies that $\tau\big[ T^*(\mathbf a)\big] = \Phi\big( g(\mathbf a)^*\big) $ which is equal to $\overline{ \Phi\big( g(\mathbf a)\big)} = \overline{ \tau\big[ T(\mathbf a)\big]}$.
\end{proof}

We can then define two equivalent formulations for the distribution $\tau_{\mathbf a}: T \in \mathbb C\mathcal T\langle \mathbf x \rangle \mapsto \tau\big[T(\mathbf a)\big] $ of a family $\mathbf a$ of traffics: they are given by the linear forms $\Phi$ and $\Psi$ on the space $\mathbb C \mathcal G \langle \mathbf x \rangle$ of graphs polynomials, namely 
\begin{itemize}
	\item the map $\bar \Phi_{\mathbf a}: g\in \mathbb C \mathcal G \langle \mathbf x \rangle   \mapsto   \Phi\big( g(\mathbf a) \big)  $, 
	\item and the map $\bar \Psi_{\mathbf a}: g\in \mathbb C \mathcal G \langle \mathbf x \rangle   \mapsto   \Psi\big( g(\mathbf a) \big)$. 
\end{itemize}
The \emph{traffic distribution} is generic term for the three maps $\tau_{\mathbf a}, \bar \Phi_{\mathbf a}, \bar \Psi_{\mathbf a}$. The restrictions $\Phi_{\mathbf a}$ of $\bar \Phi_{\mathbf a}$ and $\Psi_{\mathbf a}$ of $\bar \Psi_{\mathbf a}$ on $\mathbb C  \langle \mathbf x \rangle$ are the non commutative distributions of $\mathbf a$ with respect to $\Phi$ and $\Psi$ respectively.

Moreover, from Chapter \ref{Sec:DefIndepTraffic}, we can also characterize the traffic distribution of $\mathbf a$ thanks to the injective version $\tau_{\mathbf a}^0$ of $\tau_{\mathbf a}$. The same definition is valid to define from $\tau: \mathcal T\langle \mathcal A \rangle \to \mathbb C$ the linear map $\tau^0: \mathcal T\langle \mathcal A \rangle \to \mathbb C$, namely for any test graph  $T\in \mathcal T\langle \mathcal A \rangle$ with vertex set $V$, one has $\tau[T]= \sum_{\pi \in \mathcal P(V)} \tau^0\big[T^\pi\big]$, where $T^\pi$ is the graph obtained by identifying vertices of $T$ in a same block of $\pi$.

\begin{lemma}\label{PropTraceInj} Let $(\mathcal A, \tau)$ be an algebraic traffic space. The injective version $\tau^0$ of $\tau$ satisfies the following property.
\begin{enumerate}
	\item It sends the graph with no edges to one.
	\item It is multi-linear with respect to the edges of the graphs in the sense of Definition \ref{Def:TrafficSpace}.
	\item It satisfies the following property. For any $T \in \mathcal T\langle \mathcal A \rangle$ having an edge $e_0$ labeled $h(\mathbf a)$ for a graph monomial $h$, denote $T_h\in \mathcal T\langle \mathcal A \rangle$ the graph obtained from $T$ by replacing the edge $e_0$ by the graph $h(a)$ and denote by $V_h$ its vertex set. Then one has 
	$$\tau^0\big[T\big] =\sum_{ \substack{ \tilde \pi \in \mathcal P(  V_h) \\ \mathrm{s.t.} \tilde \pi_{|V} = 0_V }}  \tau^0\big[ T_h^{\tilde \pi}\big],$$ 
where the notation $\tilde \pi_{|V} = 0_V $ (introduced in Lemma \ref{Lem:SubsInjTrace0}) means that the sum is over the partitions $\tilde \pi$ of $V_h$ such that any couple of vertices $v$ and $w$ in $T$, when seen in $T_h$ after insertion of $h$, belong to two distinct blocks of $\tilde \pi$.
\end{enumerate}

Reciprocally, in a $\mathcal G$-algebra $\mathcal A$, assume that $\tau^0: \mathbb C \mathcal T \langle \mathcal A\rangle \to \mathbb C$ is a linear map satisfying the three above properties. Then $\tau^0$ is the injective version of a map $\tau$ which defines a structure of algebraic traffic space on $\mathcal A$.
\end{lemma}

\begin{proof} We have $\tau[ \cdot ] = \tau^0[\cdot]$ so the unity axioms coincides for both maps. Likewise,  since $\tau$ and $\tau^0$ are related each other by linear combinations, the multi-linearity of the maps are also equivalent.

Let $\tau^0$ be the injective version of a linear form $\tau$ on $\mathbb C \mathcal T\langle \mathcal A \rangle$. Let $T, e_0, h(\mathbf a)$ and $T_h$ be as in the statement. Assuming that $\tau$ satisfies the substitution axiom, we already proved that $\tau^0$ satisfies the expected formula in Proposition \ref{Lem:IndPropAlg}.

So we assume now that $\tau^0$ satisfies the formula of the lemma and prove that $\tau[T] = \tau[T_h]$. On the one hand, we have $\tau[T] = \sum_{\pi \in \mathcal P(V)} \tau^0[T^\pi]$. Let us denote by $(T^\pi)_h$ the graph obtained from $T^\pi$ by replacing the edge $e_0$ by $h$ and by $V_{\pi, h}$ its vertex set. Then the formula for $\tau^0$ tells that 
	$$\tau[T]= \sum_{\pi \in \mathcal P(V)} \sum_{\substack{ \sigma\in \mathcal P(V_{\pi, h}) \\ \mathrm{s.t. \ } \sigma_{|V_\pi}=0_{V_\pi}}} \tau^0\Big[ \big( (T^\pi)_h\big)^\sigma\Big] = \sum_{\tilde \pi \in \mathcal P(V_h)} \tau^0[T_h^{\tilde \pi}].$$  

Using again the definition of $\tau^0$ yields $\tau[T] =  \tau\big[T_h\big]$.\end{proof}

\begin{lemma}\label{Lem:PremilIndp}Let $(\mathcal A, \tau)$ be an algebraic traffic space. For any test graph $T$, one has $\tau\big[T(\mathbb I , \, \cdot \,)\big] = \tau\big[\tilde T( \, \cdot \,)\big]$ where $\tilde T$ is obtained from $T$ by identifying source and target of each edge labelled $\mathbb I$, and removing these edges. Moreover, with same notations, $\tau^0\big[T(\mathbb I, \, \cdot \,)\big] = \tau^0\big[\tilde T( \, \cdot \,)\big]$ if all the edges of $T$ labelled $\mathbb I$ are self loops and it vanishes otherwise.

\end{lemma}

\begin{proof}  We have $\Delta(\mathbb I) = \mathbb I$, where $\Delta$ denotes as usual the graph monomial consisting in a single loop. So by definition of $\mathbb I$ and by the substitution axiom, $\tau\big[T(\mathbb I, \, \cdot \,)\big]$ is not modified if we identify source and target of edges labeled $\mathbb I$ and removing them. Similarly, let $T$ be a test graph with an edge $e_0=(v,w)$ labeled $\mathbb I$ which is not a self loop. Replacing the edge by its label as in Lemma \ref{PropTraceInj}, the resulting graph $T_h$ is obtained by identifying $v$ and $w$. Hence $\tau^0\big[T\big] = 0$, since in the sum of the formula for $\tau^0$ in the lemma,  there is no partition that separates $v$ and $w$. Let now $T$ be a test graph such that all edges labeled $\mathbb I$ are self loops. By the explicit definition of $\tau^0$, one has $\tau^0\big[T(\mathbb I, \cdot)\big] = \sum_{\pi \in \mathcal P(V)} \mu_V(\pi) \tau\big[ T^\pi(\mathbb I, \cdot)\big]$ where $V$ is the vertex set of $T$. Denoting by $\tilde T^\pi$ the graph obtained from $T^\pi$ by erasing self loops labeled $\mathbb I$, by the above we have $\tau\big[ T^\pi(\mathbb I, \cdot)\big] = \tau\big[ \tilde T^\pi(\cdot)\big]$, and then $\tau^0\big[T(\mathbb I, \cdot)\big]  = \sum_{\pi \in \mathcal P(V)} \mu_V(\pi) \tau\big[ \tilde T^\pi(\, \cdot\, )\big] = \tau^0\big[\tilde T(\, \cdot \, )\big]$
\end{proof}

\begin{remark}\label{Rk:DeltaLoop}The reasoning shows actually that $\tau^0\big[T( \Delta(a), \, \cdot \,)\big] = 0$ for any traffic $a$ and any $T$ such that there is at least one edge labeled $\Delta(a)$ which is not a self-loop.
\end{remark}

\chapter{Traffic Independence and the Three Classical Notions}\label{Sec:ThreeIndep}
The product of traffic distributions, discovered in the first part for large matrices, is considered in the context of algebraic traffic spaces, for which it defines the notion of traffic independence. It is shown to unify the three classical notions of non commutative independence.
\section{Definition and statement}

For convenience we recall here the classical notions of non commutative independence.
\begin{definition}[Non commutative notions of independence]
\label{Def:ClassicalIndp}~
	\begin{enumerate}
		\item The unital subalgebras $\mathcal A_1 \etc \mathcal A_L$ of a non commutative probability space $(\mathcal A, \Phi)$ are said to be freely independent if and only if for any $n\geq1$, any $a_j\in \mathcal A_{\ell_j}, j =1\etc n$, such that $\Phi(a_j)=0$ and $\ell_{j}\neq \ell_{j+1}$, $\ell_j\in \{1\etc L\}$, one has $\Phi(a_1 \dots a_n)=0$. 
		\item  The subalgebras $\mathcal A_1\etc \mathcal A_L$ (non necessarily unital) of an algebra $\mathcal A$ are said to be Boolean independent with respect to a linear form $\Psi$ on $\mathcal A$ if and only if for any $n\geq1$, any $a_j\in \mathcal A_{\ell_j}, j =1\etc n$, such that $\ell_{j}\neq \ell_{j+1}$, one has $\Psi(a_1 \dots a_n)=\Psi(a_1) \times \dots \times \Psi(a_n)$. 
		\item The unital subalgebras $\mathcal A_1 \etc \mathcal A_L$ of a non commutative probability space $(\mathcal A, \Phi)$ are said to be tensor independent if and only if they commute (i.e. $ab = ba, \, \forall a\in \mathcal A_\ell, b\in \mathcal A_{m}, \ell\neq m $) and for any $a_\ell \in \mathcal A_\ell$, $\ell=1\etc L$, one has $\Phi( a_1 \dots a_L) = \Phi(a_1) \times \dots \times  \Phi(a_L) $.
	\end{enumerate}
\end{definition}

\begin{definition}[Traffic independence] Let $(\mathcal A, \tau)$ be an algebraic traffic space and let $\mathcal A_1 \etc \mathcal A_L$ be unital $\mathcal G$-subalgebras of $\mathcal A$. We say that $\mathcal A_1 \etc \mathcal A_L$ are independent (or traffic independent if necessary to avoid confusion) whenever the restriction of $\tau$ on the $\mathcal G$-algebra spanned by $\mathcal A_1\etc \mathcal A_L$ is the product of $\tau_1 \etc \tau_L$ in the sense of Definition \ref{def:FreeProdGraphs}: for every family $\mathbf a_\ell$ of $\mathcal A_\ell$, $\ell=1\etc L$, for any test graph $T$ in variables $\mathbf x_1 \etc \mathbf x_L$, 
	\eqa\label{Eq:DefIndTraffBis}
		\tau^0\big[T(\mathbf a_1\etc \mathbf a_L)\big] = \one\big( \mathcal G\mathcal C\mathcal C(T)\mathrm{ \ is \ a \ tree} \big) \prod_{S \in \mathcal C \mathcal C(T)}\tau^0\big[S\big],
	\qea
where $\mathcal G\mathcal C\mathcal C(T)$ (respectively $\mathcal C\mathcal C(T)$) is the graph (respectively the set) of colored components of $T$ with respect to $\mathbf x_1 \etc \mathbf x_L$ (see Definition \ref{Def:GCC}).

Subsets of $\mathcal A$ or families of elements of $\mathcal A$ are said to be (traffic) independent whenever the $\mathcal G$-subalgebra they spanned are independent.  Let $(\mathcal A_N, \tau_N)_{N\geq 1}$ be a sequence of algebraic traffic spaces. A sequence of families $\mathbf a\toN_1 \etc \mathbf a\toN_L$ of element of $\mathcal A_N$ is said to be asymptotically independent whenever it converges toward independent families. \end{definition}

The map $\tau$ on the $\mathcal G$-subalgebra spanned by independent $\mathcal G$-subalgebras is completely determined by the restriction of $\tau$ on each $\mathcal G$-subalgebras: for any test graph $T$ with vertex set $V$ and any $\mathbf a_1 \etc \mathbf a_L$
	$$\tau\big[T(\mathbf a_1 \etc \mathbf a_L)\big] = \sum_{\substack{ \pi \in \mathcal P(V)   \\ \mathrm{s.t.} \ \mathcal G\mathcal C\mathcal C(T^\pi)\mathrm{ \ is \ a \ tree} }}  \prod_{S \in \mathcal C \mathcal C(T^\pi)}   \tau^0\big[S(\mathbf a_1 \etc \mathbf a_L)\big].$$

\begin{lemma}\label{Lem:TraffiInpSpTr} Independence of traffics subspaces is symmetric and associative, in the sense stated in Proposition \ref{Lem:Associativity}. Moreover, if the traffic distribution of $\mathbf a_1 \etc \mathbf a_L$ is the product of the traffic distributions of the $\mathbf a_\ell$'s (Formula \eqref{Eq:DefIndTraffBis} is true for $\mathbf a_1 \etc \mathbf a_L$ given and for all $T$), then $\mathbf a_1 \etc \mathbf a_L$ are independent (Formula \eqref{Eq:DefIndTraffBis} is true for any graph polynomials $g(\mathbf a_\ell)$). 
\end{lemma}

The first part of the Lemma is a direct consequence of Proposition \ref{Lem:Associativity}. The second part follows thanks to Proposition \ref{Lem:IndPropAlg}, valid with no modification of the proof. Hence families of matrices $\mathbf A_N^{(1)} \etc \mathbf A_N^{(L)}$ are asymptotically traffic independent in the sense of the first part of the article if and only if the families  $\mathbf A_N^{(1)}\cup\mathbf A_N^{(1)*}  \etc \mathbf A_N^{(L)}\cup \mathbf A_N^{(L)*}$ are asymptotically traffic independent in the sense of the above definition.

\begin{lemma}\label{Lem:Elementary} If $a$ and $b$ are independent traffics then $\Phi(a b ) = \Phi(a) \Phi(b)$ and $\Psi(ab) = \Psi(a) \Psi(b)$.
\end{lemma}

The proof is left as an exercice for the reader (see the computations in Section \ref{Sec:CritLackInp}).

We can distinguish three particular kinds of traffics, for which traffic independence can be interpreted in terms of the different notions of independence.

\begin{theorem}[Unification of the notions]\label{Th:ThreeIndepTraf} Let $(\mathcal A, \tau)$ be an algebraic traffic space, with trace $\Phi:a \mapsto \tau[\circlearrowleft_a]$ and anti-trace $\Psi:a \mapsto \tau[ \cdot \overset a \rightarrow \cdot]$. Let $\mathbf a_1 \etc \mathbf a_L$ be independent families of elements of $\mathcal A$.
\begin{enumerate}
	\item Say that $\mathbf a =(a_j)_{j\in J}$ is "unitarily invariant" in $(\mathcal A, \tau)$ if and only if it has the same traffic distribution as $u \mathbf a u^* = (ua_ju^*)_{j\in J}$ where $(u,u^*)$ is independent from $\mathbf a$, satisfies $uu^*=u^*u=\mathbb I$, and is the limit in traffic distribution of $(U_N,U_N^*)$ for $U_N$ a large unitary Haar matrix. If $\mathbf a_\ell$ is unitarily invariant for each $\ell=1\etc L$ (except possibly one), then $\mathbf a_1 \etc \mathbf a_L$ are free independent in the non commutative probability space $(\mathcal A, \Phi)$.
	\item Say that $\mathbf a =(a_j)_{j\in J}$ is ''of Boolean type`` in $(\mathcal A, \tau)$ if and only if its combinatorial distribution is supported on trees. If $\mathbf a_\ell$ is of Boolean type for each $\ell=1\etc L$, then $\mathbf a_1 \etc \mathbf a_L$ are Boolean independent with respect to $\Psi$. 
	\item Say that an element of a $\mathcal G$-algebra $\mathcal A$ is "diagonal" whenever $a=\Delta(a)$ for $\Delta$ the graph operation with a single vertex and a single self loop. If the traffics of $\mathbf a_\ell$ are diagonal for each $\ell=1\etc L$, then $\mathbf a_1 \etc \mathbf a_L$ are tensor independent, both in the non commutative probability spaces $(\mathcal A, \Phi)$ and $(\mathcal A, \Psi)$. Reciprocally the tensor independence of diagonal families of traffics in the non commutative probability space $(\mathcal A, \Phi)$ (resp. $(\mathcal A, \Psi)$) characterizes their traffic independence.
\end{enumerate}

\end{theorem}
 The theorem is proved in the three following sections. In the case of Boolean independence, we exhibit examples of random matrices whose limits are of Boolean type.

Given traffic spaces $\mathcal A_1\etc \mathcal A_L$, it is natural to ask if there exists a traffic space $\mathcal A$ containing the $\mathcal A_\ell$ as independent spaces. This fact is true and is proved in \cite{CDM16}. This implies the existence, for any traffic $a$, of a space containing a sequence $(a_n)_{n\geq 1}$ of independent traffics distributed as $a$. In this article, we always assume for granted the existence of such sequences.

\section{Link with free independence}\label{Sec:LinkFree}

In order to prove the first point of Theorem \ref{Th:ThreeIndepTraf}, it suffices to prove the following.

\begin{lemma}\label{Lem:RigFree} Let $\mathbf a, \mathbf b$ be two arbitrary families of traffics. 
\begin{enumerate}
	\item Assume that $(\mathbf a, \mathbf b)$ is traffic independent from $(u,u^*)$, satisfying $uu^*=u^*u=\mathbb I$, and limit in traffic distribution of $(U_N,U_N^*)$ for a Haar unitary matrix $U_N$. Then $ \mathbf a$ and $  u \ \mathbf b \, u^*$ are free independent.
	\item Assume that $\mathbf a, \mathbf b$ are traffic independent and unitary invariant. Then the joint family $\mathbf a \cup \mathbf b$ is unitary invariant.
\end{enumerate}
\end{lemma}

Indeed, let $\mathbf a^{(1)} \etc \mathbf a^{(L)}$ be as in Proposition \ref{RigFree}. Assume that $\mathbf a^{(2)}\etc \mathbf a^{(L)}$ are unitarily invariant and denote $\mathbf b = \mathbf a^{(2)} \cup \dots \cup \mathbf a^{(L)}$. Let $(u,u^*)$ be traffic independent from $(\mathbf a^{(1)} \etc \mathbf a^{(L)})$ and limit of $(U_N,U_N^*)$ as in the lemma. Then the first point of the above lemma implies that $ \mathbf a^{(1)} $ and $u \mathbf b u^*$ are free independent. 

On the other hand, the associativity of traffic independence implies that $  \mathbf a^{(1)} $ and $\mathbf b$ are traffic independent, then that $\mathbf a^{(1)} $ and $u \mathbf b u^*$ are traffic independent. Their joint traffic distribution depends only on the marginal distributions. But thanks to the second part of the lemma, $u \mathbf b u^*$ has the same traffic distribution as $\mathbf b$ which is unitarily invariant. Hence $(\mathbf a^{(1)}, u \mathbf b u^*)$ has the same traffic distribution as $(\mathbf a^{(1)}, \mathbf b)$. Hence they have the same $^*$-distribution and so $\mathbf a^{(1)}$ and $\mathbf b$ are free independent. We get the proposition by induction on $L$ and thanks to the associativity of free independence \cite{NS}.

This implies a useful criterion of asymptotic free independence in the context of the asymptotic traffic independence theorem.

\begin{corollary}~\label{RigFree}
Let $\mathbf A_N^{(1)} \etc \mathbf A_N^{(L)}$ be as in Theorem \ref{MainTh}. Assume moreover that each family $\mathbf A_N^{(\ell)}$ has the same limiting traffic distribution as $U_N \mathbf A_N^{(\ell)} U_N^*$ for any unitary matrix $U_N$, except possibly for one index $\ell \in \{1\etc L\}$. Then $\mathbf A_N^{(1)} \etc \mathbf A_N^{(L)}$ are asymptotically free independent.
\end{corollary}
Note that the additional assumption in Proposition \ref{RigFree} is much weaker than assuming the unitary invariance of the $\mathbf A_N^{(\ell)}$ (that is $\mathbf A_N^{(\ell)} \overset{\mathcal L}= U_N\mathbf A_N U_N^*$ for each $N$ and any unitary matrix $U_N$). For instance, independent Wigner matrices with parameter of the form $(\alpha, 0)$ satisfy this proposition. This result is applied for a class of large random graphs with large degree in \cite{MP14}. 
A much detailed analysis of the relation between traffic independence, free independence and unitarily invariance is made in \cite{CDM16}, with an explicit description of the limiting traffic distribution of unitary invariant traffics. 

\begin{proof}[Proof of Lemma \ref{Lem:RigFree}] We first prove the first part of the lemma. Let $n\geq 1$ be an integer and $P_1 \etc P_n, Q_1 \etc Q_n$ be non commutative polynomials and denote $y_i = P_i( \mathbf a )$, $z_i = Q_i(  \mathbf b)$, $\mathbf y = (y_i)_{i=1\etc n}$ and $\mathbf z=  (z_i)_{i=1\etc n}$. Note that $u y_iu^* = P_i( u\mathbf a u^* )$ for any $i=1\etc n$ since $uu^*=u^*u=\mathbb I$. We assume  $\Phi(y_i)=0$, $\Phi( z_i ) =0$.  Proving for any $i=1\etc n$ that  $ \Phi( y_1 z_1  \dots  y_nz_n  ) $ vanishes we will get the free independence of $u \, \mathbf a\, u^*$ and $\mathbf b$.

One has  $\Phi( y_1 z_1  \dots  y_nz_n  ) =  \Phi (u y_1u^* z_1  \dots  uy_nu^*z_n )$. Let $T=(V,E,\gamma,\varepsilon)$ be the test graph in variables $u, u^*,y_1 \etc y_n, z_1 \etc z_n$ such that 
	$$ \Phi (u y_1u^* z_1  \dots  uy_nu^*z_n ) = \tau[T],$$
 namely
\begin{itemize}
	\item the set of vertices is $V = \{1,2 \etc 4n\}$,
	\item the edges are $(1,4n), (4n,4n-1), (4n-1,4n-2) \etc  (3,2), (2,1)$,
	\item with notation of indices modulo $2n$, the edges $(4i+2,4i+1)$ are labelled $u$, the edges $(4i+3, 4i+2)$ are labelled $y_i$, the edges $(4i+4, 4i+3)$ are labelled $u^*$, and the edge $(4i+5, 4i+4,)$ are labelled $z_i$.
\end{itemize}

By the relation between the combinatorial trace and its injective version (Formula \eqref{eq:TraffCum}), one has $\tau[T]= \sum_{\pi\in \mathcal P(V)} \tau^0[T]$. Moreover, since $(u,u^*)$ and $(\mathbf a, \mathbf b)$ are traffic independent, then $(u,u^*)$ and $(y_i, z_i)_{i=1\etc n}$ are traffic independent. Hence by definition of traffic independence, we get 
	\eq
		\tau^0[T^\pi] =  \one\big( \mathcal G \mathcal C \mathcal C(T^\pi) \textrm{ is a tree} \big)   \prod_{ S\in \mathcal C \mathcal C(T^\pi)} \tau^0[S],
	\qe
where $ \mathcal G \mathcal C \mathcal C(T^\pi)$ (respectively $\mathcal C \mathcal C(T^\pi)$) is the graph (the set) of colored components of $T^\pi$ with respect to the variables $(u,u^*)$ and $(\mathbf y,\mathbf z)$. 
\par Recall that the injective distribution of $u$ is supported on cacti (Proposition \ref{Prop:HaarDist}), namely graphs such that each edge belongs exactly to one simple cycle. Moreover, the edges of a same cycle must be oriented in a same direction and labels must alternate between $u$ and $u^*$.

 Given $\pi$ as in the sum, denote by $S(T^\pi)$ the $^*$-test graph obtained from $T^\pi$ by identifying the vertices attached to a same connected component labelled in $\mathbf y$ or $\mathbf z$, and forgetting the edges labelled in $\mathbf y, \mathbf z$. To ensure that $\pi$ contributes, each connected component of $T^\pi$ labeled in $\{x,x^*\}$ must be a cactus. When $\mathcal G \mathcal C \mathcal C (T^\pi)$ is a tree, the graph $S(T^\pi)$ itself must be a cactus. 

Denote by $\mathcal S(n)$ the set of cacti with $2n$ edges. Then one has
	\eqa
		\tau[T] & = & \sum_{ S_0 \in   \mathcal S(n)}  \sum_{ \substack{ \pi \in \mathcal P(V) \\ \mathrm{  s.t. } \ S(T^\pi) = S_0 }}   \one\big( \mathcal G \mathcal C \mathcal C(T^\pi) \textrm{ is a tree} \big) \prod_{ S \in \mathcal C \mathcal C(T^\pi)} \tau_N^0[S].\label{Eq:LinkFreeness}
	\qea

For any $S_0 \in \mathcal S(n)$ and any $\pi$ in $\mathcal P(V)$ such that $S(T^\pi)=S_0$, there is a colored component $S\in \mathcal C \mathcal C(T^\pi)$ labeled in $\mathbf y$ or $\mathbf z$ consisting in a loop. Indeed, there is a vertex of $S(T^\pi)$ which belong to a single cycle (each cycle has length at least two). So since the variables $u,y_i,u^*,z_i$ in $T$ alternates there is such a loop in $T^\pi$ corresponding to this vertex.  This implies that $\tau^0[T^\pi] $ is zero since the combinatorial trace of a loop is the trace of the variable.

It remains to prove the second part of the lemma. Let $\mathbf a, \mathbf b$ be traffic independent and unitary invariant. Let $u,v,w$ be traffics, limits of Haar unitary matrices (as long with their adjoint), such that $\mathbf a, \mathbf b, u, v, w$ are independent. We prove that $(u \mathbf a u^*, u \mathbf b u^*)$ has the same traffic distribution as $( \mathbf a ,\mathbf b)$. By unitary invariance and associativity of independence, $(u \mathbf a u^*, u \mathbf b u^*)$ has the same distribution as $(uv \mathbf a v^*u^*, uw \mathbf b w^*u^*)$. But $(uv, uw)$ and $(a,b)$ are independent and the matrix approximation (Theorem \ref{MainTh} for independent Haar unitary matrices) shows that $(uv,uw)$ has the same traffic distribution as $(v,w)$. Hence $(u \mathbf a u^*, u \mathbf b u^*)$ has the same distribution as $(v\mathbf a v^*, w \mathbf b w^*)$, and so the same distribution as $(\mathbf a, \mathbf b)$.
\end{proof}

\section{Link with tensor independence}\label{Sec:LinkTens}

\begin{lemma}\label{Lem:DiagTrafff} Let $(\mathcal A, \tau)$ be an algebraic traffic space with trace $\Phi$ and anti-trace $\Psi$. The space $\Delta(\mathcal A)= \{ \Delta(a), a\in \mathcal A\}$ of diagonal elements is a commutative $\mathcal G$-subalgebra of $\mathcal A$. Moreover, for any family $\mathbf a$ of diagonal traffics and for any test graph $T=(V,E,\gamma)$, one has
	\eq
		\tau\big[T(\mathbf a)\big]  =   \Phi\big[  \prod_{e\in E} a_{\gamma(e)} \big] =  \Psi\big[  \prod_{e\in E} a_{\gamma(e)} \big], \ \ \tau^0\big[T(\mathbf a)\big]   =  \one(|V|=1) \times \tau\big[T(\mathbf a)\big]  .
	\qe
\end{lemma}

\begin{proof} Let $\mathbf a$ be a family of elements of $\Delta(\mathcal A)$ and $g=(T,in,out)$, $T=(V,E,\gamma)$, a graph monomial. Since $\Delta(\mathbf a) = \mathbf a$ and by the associativity of the composition of graph monomials, $g(\mathbf a) = \tilde g(\mathbf a)$ where $\tilde g$ is obtained from $g$ by identifying source and goal of each edge. So $\tilde g$ is a bunch of self-loops, independently of the geometry of the initial graph $g$. Hence, we have $\tilde g(\mathbf a) = \prod_{e\in E} a_{\gamma(e)}$ and $g(a)$ is diagonal. Since $\Delta$ is linear, we get that $\Delta(\mathcal A)$ is a $\mathcal G$-sublagebra. Moreover $\Delta(a) \times \Delta(b) = Z_g(a,b) = \Delta(b) \times \Delta(a)$, where $g$ is the graph operation with two self-loops attached to a single edge. Hence the chain of equality comes from the associativity of composition and the equivariance axiom. 

The formula for $\tau\big[T(\mathbf a)\big]$ follows from the previous paragraph since we can write $\tau\big[T(\mathbf a)\big] = \Phi\big( g(\mathbf a)\big) = \Psi \big( h(\mathbf a)\big)$ for some graph monomials $g,h$. It remains to prove the formula for $\tau^0$. Let $T=(V,E,\gamma)$ be a test graph. Denote by $1_V= \{V\}$ the partition of $V$ with a single block. By the previous point we have $\tau\big[T] = \tau\big[T^{1_V}\big]$. One the other hand,
	$$\tau\big[ T\big] = \sum_{\pi \in \mathcal P(V)} \tau^0\big[T^\pi\big], \ \tau\big[ T^{1_{V}}\big] = \tau^0\big[ T^{1_{V}}\big].$$
Hence we get that $\sum_{\pi \in \mathcal P(V) \setminus \{1_{V}\} } \tau^0\big[T^\pi\big], \ \tau\big[ T^{1_{\mathcal P(v)}}\big] =0$, which gives that $\tau^0\big[T]$ for $|V|\geq 2$ by a induction on $|V|$.
\end{proof}

\begin{proof}[Proof of the third item of Theorem \ref{Th:ThreeIndepTraf}]

Assume that the families of traffics $\mathbf a^{(1)} \etc \mathbf a^{(L)}$ are independent and diagonal. For each $\ell=1\etc L$, let $P_\ell$ be a commutative monomial. We can write $\Phi  \big[ \prod_{\ell=1}^L P_\ell(\mathbf a^{(\ell)}) \big]  = \tau\big[T(\mathbf a)\big]$ where $T$ is the graph with a single vertex and one self-loop for each variable appearing in the monomials.

Since $T$ has a single vertex, $ \tau\big[T(\mathbf a)\big] =  \tau^0\big[T(\mathbf a)\big]$, and by traffic independence, $ \tau^0\big[T(\mathbf a)\big] = \prod_{\ell=1}^L  \tau\big[T_\ell(\mathbf a)\big],$
where $T_\ell$ is the subgraph of $T$ whose edges labels correspond to $\mathbf a^{(\ell)}$. But for each $\ell$, by the same reasoning, one has 
	$$\tau ^0\big[T_\ell(\mathbf a)\big] = \tau\big[T_\ell(\mathbf a)\big] =\Phi \big[ P_\ell(\mathbf a^{(\ell)})\big] .$$
Hence $\Phi  \big[ \prod_{\ell=1}^L P_\ell(\mathbf a^{(\ell)}) \big] = \prod_{\ell=1}^L \Phi \big[ P_\ell(\mathbf a^{(\ell)})\big]$ so the families of matrices are tensor independent with respect to $\Phi$. Moreover, $\Phi$ and $\Psi$ are equal for diagonal traffics so the families are also tensor independent with respect to $\Psi$. 

Reciprocally, we now assume the tensor independence of diagonal elements $\mathbf a^{(1)} \etc \mathbf a^{(L)}$. Let $T$ be a test graph. If $T$ has more than a single vertex, then there is a colored component $T'$ of $T$ which has the same property and so $\tau^0\big[T(\mathbf a)\big]=0$ and $\tau^0\big[T'(\mathbf a)\big]=0$. Hence the rule of traffic independence is satisfied for these graphs. If now $T$ has a single vertex, we have $ \tau^0\big[T(\mathbf a)\big] = \prod_{\ell=1}^L  \tau^0\big[T_\ell(\mathbf a)\big]$ by the above reasoning. Hence we get the traffic independence of the families.
\end{proof}

\section{Link with Boolean independence}\label{Sec:LinkBool}

Note that from the illustration of Theorem \ref{MainTh} given at the very end of Section \ref{Sec:MainTh}, the link between traffic independence and tensor and free independence is not a surprise. Nevertheless, these two notions are not sufficient if we want to explain the central limit theorem for the sum of independent traffics presented in the last section. This was the original motivation to find the second item of Theorem \ref{Th:ThreeIndepTraf}.

\subsection{Generalities and proof}

\begin{lemma} Let $\mathbf a$ be a family of elements in an algebraic traffic space $(\mathcal A, \tau)$ with associated trace $\Phi$. 
	\begin{enumerate}
		\item If $\mathbf a$ is of Boolean type, then the $^*$-distribution of $\mathbf a$ with respect to $\Phi$ is the distribution of the null element.
		\item $\mathbf a$ is of Boolean type if and only if the injective combinatorial distribution of $\mathbf a$ is supported on trees, in which case for any $^*$-test graph $T$ one has $\tau\big[T(\mathbf a)\big] = \tau^0\big[T(\mathbf a)\big]$.
		\item If $\mathbf a_1\etc \mathbf a_L$ are traffic independent and of Boolean type, then $\mathbf a_1\cup \dots \cup \mathbf a_L$ is of Boolean type.
	\end{enumerate}
\end{lemma}

\begin{proof}(1) Let $\mathbf a$ be of Boolean type. For any $^*$-monomials $M\neq \mathbb I$,  $\Phi\big( M(\mathbf a) \big)$ is equal to 
	$$\tau\big[ \circlearrowleft\big(M(\mathbf a)\big) \big] = \tau\big[ \circlearrowleft(M)(\mathbf a) \big] = 0,$$
since the graph $\circlearrowleft(M)$ is not a tree (but a simple cycle).

(2) Recall that for any test graph $T$, we have $\tau [ T ] = \sum_{\pi \in \mathcal P(V)}  \tau^0  [ T^\pi]$ and $\tau^0  [ T ] = \sum_{\pi \in \mathcal P(V)} \mu_V(\pi)  \tau  [ T^\pi]$. But $T^\pi$ is a tree if and only if $T$ is a tree and $\pi$ is the partition $0_V$ consisting in singletons. Hence if $\mathbf a$ is of Boolean type then necessarily $\tau^0  [ T ] = \tau[T^{0_V}] = \tau[T]$ for any $T$. If $\tau^0$ is supported on trees, similarly $\tau[T] = \mu_V(0_V)\tau^0[T] = \tau^0[T]$.

(3) Let now $\mathbf a_1\etc \mathbf a_L$ be of Boolean type and traffic independent. For any test graph $T$, one has $\tau^0[T] = \one\big( \mathcal G \mathcal C \mathcal C(T) \mathrm{ \ is \ a \ tree}\big) \prod_{S\in \mathcal C \mathcal C(T)} \tau^0[S]$. If $T$ is a not a tree then either $\mathcal G \mathcal C \mathcal C(T) $ is not a tree or a colored component of $T$ is not a tree. Hence $\tau^0[T] =0$. 
\end{proof}

\begin{proof}[Proof of the second item of Theorem \ref{Th:ThreeIndepTraf}]
Let $\mathbf a_1 \etc \mathbf a_L$ be traffic independent and such that the combinatorial distribution of each $\mathbf a_\ell$ is supported on trees. Let $M_1\etc M_n$ be non constant monomials. Let us prove that 
	$$\Psi\big( M_1(\mathbf a_{i_1}) \dots M_n(\mathbf a_{j_n})\big) = \prod_{j=1}^n \Psi\big( M_j(\mathbf a_{i_j})\big)$$
 for any $i_1\neq i_2\neq \dots \neq i_n$. 

Let us denote $M=M_1 \dots M_n$. Then, by the substitution axiom, we have $\Psi\big( M(\mathbf a_1\etc \mathbf a_L)\big) = \tau[T_M]$ where $T_M$ consists in a simple line, namely $T_M=(\cdot \overset{x_1}\leftarrow \cdot \dots \overset{x_Q}\leftarrow \cdot)$ whenever $M=x_1 \dots x_Q$. By the above lemma, we have $\tau[T_M] = \tau^0[T_M] = \prod_{S \in \mathcal C \mathcal C(T)} \tau^0 [S ]$. But the colored components of $T$ are the graphs $T_{M_1} \etc T_{M_n}$ constructed as $T_M$ for $M$ replaced by $M_1 \etc M_n$. Hence $\tau[T_M]  = \prod_{j=1}^n \tau^0[T_{M_j}] = \prod_{j=1}^n \tau[T_{M_j}] =  \prod_{j=1}^n \Psi\big( M_j(\mathbf a_{i_j})\big)$. 
\end{proof}

\subsection{Application to random matrices}\label{Sec:ExBool}

Recall that $\deg(A_N)$ denotes the diagonal matrix $diag\big( \sum_{j=1}^N A_N(i,j)\big)_{i=1\etc N}$. 

\begin{corollary}\label{Prop:LinkBool} The anti-trace $\Psi_N$ of the algebraic traffic space of matrices is given by
	$$\Psi_N[A_N] =   \esp\Big[ \frac 1 N \Tr \deg A_N \Big] = \esp\Big[ \frac 1 N \sum_{i,j=1}^N A_N(i,j)\Big]= \langle  e_N, A_N e_N \rangle,$$
where $\langle \, \cdot \, , \, \cdot \, \rangle$ denotes the usual scalar product in $\mathbb C^N$ and $e_N$ the column vector whose all entries are $\frac 1 {\sqrt N}$. If $\mathbf A_N$ converges in traffic distribution, then it has a limiting distribution with respect to $\Psi_N$, that is $\Psi_N\big[P(\mathbf A_N)\big]$ converges for any non commutative polynomial $P$. If $\mathbf A_N^{(1)} \etc \mathbf A_N^{(L)}$ are asymptotically traffic independent families of matrices whose limiting combinatorial distributions are supported on simple trees, then the families $\mathbf A_N^{(1)} \etc \mathbf A_N^{(L)}$ are asymptotically Boolean independent with respect to $\Psi_N$.
\end{corollary}

Let us now give examples of such matrices. The matrix $\mathbb J_N$, whose all entries are $\frac 1 N$, converges to a traffic of Boolean type: for any $^*$-test graph $T$, by Lemmas \ref{Lem:Tool} and \ref{lem:EdgesVertices},
		$$\tau_N\big[T(\mathbb  J_N)\big] = N^{|V|-1-|E|}\Big( 1 + O\big(\frac 1 N\big) \Big) = \one(T \textrm{ is a tree}) +o(1).$$
	Note that one has $\Psi_N(A_N) = \esp\big[  \Tr ( A_N \mathbb J_N)\big]$ and $\mathbb J_N = e_N e_N^*$ where $e_N$ is as in Proposition \ref{Prop:LinkBool}. Furthermore, $\mathbb J_N$ is a deterministic permutation invariant matrix, and any permutation invariant deterministic matrix $A_N$ is of the form $\Phi_N(  A_N)\mathbb  I_N + \big( \Psi_N(  A_N) - \Phi_N(  A_N) \big) \mathbb J_N$, where $\Phi_N(  A_N) = \esp\big[ \frac 1 N \Tr   A_N \big]$ and $\mathbb I_N$ is the identity matrix. It is also a projection matrix, namely $\mathbb J_N^2=\mathbb J_N$. 
	
	The limiting distribution of $\mathbb J_N$ with respect to $\Psi_N$ is the distribution of a variable constant to one. Indeed, one has $deg(\mathbb J_N)=\mathbb I_N$ so for any monomial $K\geq 1$, one has $\Psi_N\big[\mathbb J_N^K\big]=\esp\big[ \frac 1N \Tr \,\mathbb  I_N\big] =1$ and so for any polynomial $P$ one has $\Psi_N\big[P(\mathbb J_N)\big]=P(1)$.
	
	Let see now a non constant example.
	
\begin{lemma}\label{Lem:BoolGaus}Let $X_1\etc X_N$ be independent and identically distributed complex random variables. Assume all the moments of the $X_i$'s are finite and do not depend on $N$. We denote the permutation invariant matrix $M_N = \big(\frac{X_i+\overline{ X_j}}N\big)_{i,j=1\etc N}$.  Then $M_N$ has a limiting traffic distribution supported on trees and it satisfies the factorization property (Assumption B3 of Theorem \ref{MainTh}): denoting $\alpha_{\ell,k} = \esp[X_i^\ell \overline{ X_i}^k]$, for any test graph $T=(V,E)$, one has 
	\eqa\label{Eq:DistrGaussBoolMatrices}
			  \tau_N^0\big[ T(M_N) \big]  \limN   \one \big( T \textrm{ is a tree}\big)
		  								\underbrace{\sum_{\pi\in \mathcal P(E|V)} \prod_{B_v\in \pi} \alpha_{\ell_{B_v} , k_{B_v}}}_{\alpha(T)},
	\qea
where $\mathcal P(E|V)$ is the set of partitions $\pi = \{ B_v, v \in V\}$ of edges of $T$ into blocks $B_v$ having in common a same vertex $v\in V$, and $\ell_{B_v}$ (respectively $k_{B_v}$) is the number of edges of $B_v$ for which $v$ is the source (respectively the target). Moreover, if the $X_j$ are centered then the limiting distribution of $M_N$ with respect to $\Psi_N$ is the Rademacher distribution  $\frac{\delta_{-\alpha_{1,1}} + \delta_{\alpha_{1,1}}}2$, ($\alpha_{1,1} = \esp[|X_j|^2]$), namely for any $K\geq 1$, $\Psi_N[ M_N^K]\limN  \alpha_{1,1}^{\frac K 2}\one(K \textrm{ even})$.

In particular, if the variables $X_j$ are real Gaussian random variables centered with variance one, then $\alpha(T)$ is the number of partitions of edges of $T$ whose blocks consist in two edges having a vertex in common. If the $X_j$ are complex Gaussian random variables such that $\esp[X_j] = \esp[X_j^2] =0$ and $\esp[|X_j|^2]=1$, then $\alpha(T)$ is the number of partitions of edges of $T$ whose blocks consist in two edges having a vertex in common, which is the source for one edge and the target of the other one.
\end{lemma}

The former lemma, Theorem \ref{MainTh} and Proposition \ref{Prop:LinkBool} yield the following example of asymptotic Boolean independent matrices.

\begin{corollary}\label{Cor:ExBooleanIndpMatrices}
Let $\mathbf M_N$ be a family of independent matrices as in Lemma \ref{Lem:BoolGaus} and let $\mathbb J_N$ be the matrix whose all entries are $\frac 1 N$. Then the matrices of $\mathbf M_N$ and $\mathbb J_N$ are asymptotically traffic independent and are asymptotically Boolean independent with respect to $\Psi_N$. Moreover, if the variables $X_i$ defining $\mathbf M_N$ are complex Gaussian variables such that $\esp[X_i] = \esp[X_i^2] = 0$ and $\esp[|X_i|^2]=1$, then the matrices of $\mathbf M_N$ and their transpose are asymptotically Boolean independent. 
\end{corollary}

\begin{proof}[Proof of Lemma \ref{Lem:BoolGaus}] For any $^*$-test graph $T = (V,E,\gamma)$ in one variable, one has (Lemma \ref{Lem:Tool})
	\eq
		 \tau_N^0\big[T(M_N)\big] & = &  N^{V-1}\delta_N^0\big[T(M_N)\big]\times \Big( 1 +O\big( \frac 1 N \big) \Big) \\
		 & =& N^{|V|-1-|E|} \delta_N^0\big[T(N M_N)\big]
 \times \Big( 1 +O\big( \frac 1 N \big) \Big),
	\qe
where $\delta_N^0\big[T(N M_N)\big]= \esp\big[\prod_{\{v,w\}\in E}(X_{\phi(v)} + \overline{X_{\phi(w)}})\big]$ for any injection $\phi:V\to [N]$. It does not depend on $\phi$ and $N$.  Hence by Lemma \ref{lem:EdgesVertices}, $\tau_N^0\big[T(M_N)\big]$ converges and its limit is zero if $T$ is not a tree. Let $T$ be a tree and denote $\alpha_N(T):=\delta_N^0\big[T(NM_N) \big] $. To compute the limit of $\alpha_N(T)$ we expand the product over $E$ and the sums in its definition, which amounts for each edge $e=(v,w)$ to keep either the variable attached to its source $X_{\phi(v)}$ or to its target $X_{\phi{(w)}}$. Since the variables $X_1\etc X_N$ are independent, this yields Formula \eqref{Eq:DistrGaussBoolMatrices}.

The proof of the factorization property is the same as for Wigner matrices. Let $T_1   \etc T_n $ be test graphs in one variable, and denote by $T$ the graph obtained as the disjoint union of $T_1 \etc T_n$. By Lemma \ref{lem:PropConvEsp}, 
\eq
	 \esp\Big[ \prod_{i=1}^n \frac 1 N \Tr^0 \big[ T_i(M_N) \big] \Big]  
	  & = &  \sum_{\pi}  \frac 1 {N^{n}}   \esp\Big[   \Tr^0 \big[ T^\pi(M_N) \big]  \Big],\\
	  & = &  \sum_{\pi}   N^{V_\pi-E_\pi-n}\Big( 1 +O\big( \frac 1 N \big) \Big) \delta_N^0\big[T^\pi(NM_N) \big]   ,
\qe
where the sum is over partitions that do not identify two vertices of a same $T_i$ and $V_\pi, E_\pi$ denote the vertex and edge sets of $T^\pi$. We have by Lemma \ref{lem:EdgesVertices} that $V_\pi-E_\pi-n\leq 0$ with equality if and only if $\pi$ is the trivial partition with only singleton blocks and the graphs $T_1\etc T_n$ are trees. Moreover, the matrix entries in $\delta_N^0\big[T^\pi(NM_N) \big] $ associated to edges of different components of $T^\pi$ are independent. Hence we obtain  the factorization property
\eq
	\esp\Big[ \prod_{i=1}^n \frac 1 N \Tr^0 \big[ T_i(M_N) \big] \Big]    &  =  &     \prod_{i=1}^{n} \one(T_i \textrm{ is a tree}) \delta_N^0\big[T_i(NM_N) \big] +o(1).
	 \qe

Let compute now the limiting distribution of $M_N$ with respect to $\Psi_N$, assuming the $X_i$ defining $M_N$ centered. If $T$ is a directed line of odd length $2n+1$, then it is not possible to find a term which is not zero (the partitions $\pi$ of $\mathcal P(E|V)$ possess a block of size one), and so $\Psi_N[M_N^{2n+1}]\limN 0$. If $T$ has an even length $2n$, then there is a unique way to get a non zero term in the expansion of $\alpha_N(T)$, which gives $\alpha_N(T) = \esp[|X_i|^2]^n$. This proves the convergence of $M_N$ with respect to $\Psi_N$ to the expected limit.

The formulas for $\alpha(T)$ when the $X_i$ are Gaussian follows from Wick formula: for a Gaussian random variable $X$, in the real case $\esp[ X^{k}]$ is equal to the number of pair partitions of $\{1\etc k\}$ and in the Gaussian case (with $\esp[X^2=0]$) $\esp[ X^{k} \overline{X}^\ell]$ is the number of bijections $\{1\etc k\} \to \{1\etc \ell\}$. This is a direct consequence of the stability of the Gaussian distribution.  
\end{proof}
\begin{proof}[Proof of Corollary \ref{Cor:ExBooleanIndpMatrices}]
The first part of Corollary \ref{Cor:ExBooleanIndpMatrices} is consequence of Theorem \ref{MainTh} and Corollary \ref{Prop:LinkBool}. Let us prove that when the variables defining $\mathbf M_N$ are complex Gaussian random variables $X_{j}$ such that $\esp[X_j^2]=0$, the matrices of $\mathbf M_N$ are asymptotically Boolean independent along with the transposed matrices. The proof is the same as for the analogue statement (Lemma \ref{Lem:FreeTranspose}) for Wigner matrices. Let $P$ be a monomial in $\mathbf M_N$ and $\mathbf M_N^t$. We can write $\Psi_N\big[ P(\mathbf M_N, \mathbf M_N^t)\big] = \tau_N\big[ T(\mathbf M_N)\big]$ where the test graph $T$ is a simple line, with edges corresponding to $\mathbf M_N$ in one direction and those corresponding to $\mathbf M_N^t$ in the other direction. But in \eqref{Eq:DistrGaussBoolMatrices}, the partitions must pair edges with same orientation. Hence the couple $(\mathbf M_N, \mathbf M_N^t)$ has the same limiting distribution with respect to $\Psi_N$ as $(\mathbf M_N, \tilde {\mathbf M_N})$ where $\tilde {\mathbf M_N}$ is an independent copy of $\mathbf M_N$. Hence their asymptotic Boolean independence.
\end{proof}

\begin{corollary} Let $M_N =\big( \frac{ X_j +\overline{ X_j} }N \big)_{i,j}$ as in Lemma \ref{Lem:BoolGaus} where $\esp[X_j]=\alpha\neq 0$ and $\esp[|X_j -\alpha|^2]=1$. Then the limiting distribution of $M_N$ with respect to $\Psi_N$ is 
	$$\frac 1 2 \left( 1+ \frac {\Re(\alpha)}{\sqrt{ \Re(\alpha)^2+1}}  \right) \delta_{\Re(\alpha)+\sqrt{ \Re(\alpha)^2+1}} + 
\frac 1 2 \left( 1 -\frac {\Re(\alpha)}{\sqrt{ \Re(\alpha)^2+1}}  \right) \delta_{- \big(\sqrt{ \Re(\alpha)^2+1}  -\Re(\alpha) \big)}.$$
\end{corollary}

Note that we have  $M_N =\tilde M_N+ 2 \Re(\alpha) \mathbb J_N$ where $\tilde M_N= \big( \frac{ (X_j-\alpha) +(\overline{ X_j}-\bar \alpha)}N \big)_{i,j}$ so the limiting distribution of $M_N$ with respect to $\Psi_N$ is the Boolean convolution of $\frac{\delta_{-1} + \delta_1}2$ with $\delta_\alpha$. 

\begin{proof}
Note that $\Psi_N(A_N \mathbb J_N B_N) = \Psi_N(A_N)\times  \Psi_N(B_N)$ for any matrices $A_N,B_N$. Hence with the above notation we have
 	\eq
		\Psi_N(M_N^{K+2}) = \Psi_N(\tilde M_N^2 M_N^K) + 2\Re(\alpha) \Psi_N(\tilde M_N) \Psi_N( M_N^K)  + 2\Re(\alpha)  \Psi_N(M_N^{K+1}).
	\qe
By the expression of the limiting distribution of $\tilde M_N$, we have $\Psi_N(\tilde M_N^2 M_N^K) = \Psi_N( M_N^K) +o(1)$ and $\Psi_N(\tilde M_N) =o(1)$. Hence the sequence of moments $(m_K)_{K\geq 0}$ of the limiting distribution of $M_N$ with respect to $\Psi_N$ satisfies the recurrence relation $m_{K+2} = \beta m_{K+1} + m_{K}$, $m_0=1$ and $m_1=\beta$, where $\beta = 2 \Re(\alpha)$. We get that for any $K\geq 0$,
	$$m_K = \frac 12 \Big(1+\frac \beta \gamma \Big) \Big( \frac{\beta+\gamma }2\Big)^K + \frac 12 \Big(1-\frac \beta \gamma\Big) \Big( \frac{\beta-\gamma }2\Big)^K , \ \beta = 2 \Re(\alpha), \gamma = \sqrt{\beta^2 +4}.$$
Hence the distribution has two Dirac mass, it is characterized by its moments, and we get the result after simplification.
\end{proof}

\chapter{Limit theorems for independent traffics}\label{Sec:DefTrafficInd}

We state the law of large number and the central limit theorem in the context of traffic independence. We see in both cases that the situation is much richer than for the classical notions of independence.

\section{Constant traffics and law of large numbers}

In a non commutative probability space, a constant non commutative random variable is an element distributed as a multiple of the identity, or equivalently an element freely independent with itself.

Recall that $\mathbb J_N$ denotes the matrix whose entries are $\frac 1 N$. We have seen in Section \ref{Sec:LinkBool} that $\mathbb J_N$ converges in traffic distribution. Denote by $\mathbb J$ a traffic distributed as the limit of $\mathbb J_N$. Recall that $deg(\cdot)$ is the graph monomial with two vertices $in=out$ and $v$ and one edge from $in$ to $v$.

\begin{proposition}\label{Prop:LGN} Let $(\mathcal A, \tau)$ be an algebraic traffic space, let $\Phi$ denote the trace associated to $\tau$ and set $\Psi = \Phi \circ deg$.
\begin{enumerate}
	\item An element $a$ of $\mathcal A$ is traffic independent from itself if and only if it has the same distribution as $\Phi(a) \mathbb I + \big( \Psi(a) - \Phi(a) \big) \mathbb J$. 
	\item Law of large numbers: Let $(a_n)_{n\geq 1}$ be a sequence of identically distributed independent traffics in $\mathcal A$ and let $a$ be distributed as the $a_n$'s. For each $n\geq 1$, denote $m_n:= \frac{a_1 + \dots + a_n}n$. Then, as $n$ goes to infinity, $m_n $ converges to $\Phi(a)  \mathbb  I + \big( \Psi(a) - \Phi(a) \big)  \mathbb J$ in traffic distribution.
\end{enumerate}
\end{proposition}
\begin{proof}[Proof of Proposition \ref{Prop:LGN}] 1. Assume that $a$ is traffic independent from itself. Let $T$ be a test graph in one variable $x$ and let $e_1\etc e_K$ be an enumeration of its edges. Let $\tilde T$ be the test graph in $K$ variables $x_1 \etc x_K$ obtained from $T$ by replacing label $x$ of the $k$-th edge of $T$ by $x_k$. Let $a_1 \etc a_K$ be independent copies of $a$. By associativity of traffic independence, $(a_1 \etc a_K)$ has the same distribution as $(a \etc a)$ and so $\tau^0\big[ T(a)\big] = \tau^0\big[ \tilde T(a_1 \etc a_K)\big]$. By definition of traffic independence, this quantity is nonzero only if the graph of colored components of $\tilde T$ with respect to $x_1 \etc x_K$ is a tree. Since the labels of the edges are pairwise distinct, this means that one obtains a tree when removing the self loops of $T$. For such a graph $T$, denote by $\ell$ its number of self loops and by $m$ its number of simple edges (edges that are not loops). Then we get by definition of traffic independence $\tau^0\big[ T(a)\big] = \tau^0\big[ \circlearrowleft^a \big]^{\ell}    \tau^0\big[ \cdot \overset{a}\leftarrow \cdot \big]^{m }$. 

But $\tau^0\big[ \circlearrowleft^a \big] = \tau\big[ \circlearrowleft^a \big] = \Phi(a)$ and $ \tau^0\big[ \cdot \overset{a}\leftarrow \cdot \big] =  \tau\big[ \cdot \overset{a}\leftarrow \cdot \big] - \tau\big[ \circlearrowleft^a \big]  = \Psi(a) - \Phi(a) = \Psi\big(a - \Delta(a)\big)$, where we used \eqref{Eq:ExRelation} in the last equality. Hence 
	$$\tau^0\big[ T(a)\big] =\Phi(a)^{\ell } \Psi\big(a - \Delta(a)\big)^m.$$

Let now prove that $a$ has the same distribution as $b=\Phi(a) \mathbb I + \Psi\big(a - \Delta(a)\big) \mathbb J$. Let $T$ be a test graph in one variable with $K$ edges $e_1\etc e_K$. Denoting $[K] =\{1\etc K\}$, for any map $\gamma : [K] \to \{1,2\}$, let $T_\gamma$ be the test graph in two variables $i$ and $j$ (for $\mathbb I$ and $\mathbb J$ respectively) obtained from $T$ by putting label $i$ for edges $e_k$ with $\gamma(k)=1$ and label $j$ otherwise. Denote by $\ell_{\gamma}$ (resp. $m_{\gamma}$) the number of edges labeled $i$ (resp. $j$) in $T_\gamma$. Then, the multi-linearity of $\tau^0$ w.r.t. the edges of the graphes implies that 
	\eq
		 \tau^0\big[T(b)\big]  = \sum_{\gamma:[K]\to \{1,2\}} \Phi(a)^{\ell_{\gamma}} \Psi\big(a - \Delta(a)\big)^{m_{\gamma}} \tau^0\big[ T_\gamma(  \mathbb I,   \mathbb J)\big].
	\qe
Denoting by $\tilde T_\gamma$ the graph obtained from $T_\gamma$ by erasing self-loops labeled $\mathbb I$, by Lemma \ref{Lem:PremilIndp} we get $ \tau^0\big[ T_\gamma(  \mathbb I,   \mathbb J)\big] =  \tau^0\big[ \tilde T_\gamma(\mathbb J)\big]$ if all edges labeled $i$ are simple loops and zero otherwise. Finally, recall from Section \ref{Sec:ExBool} that the injective distribution of $\mathbb J$ is the indicator of simple trees. Hence $\tau^0\big[ T(b)\big]$ vanishes if $T$ is not a graph consisting in a tree decorated with simple loops. Otherwise, $\tau^0\big[ T_\gamma(  \mathbb I,   \mathbb J)\big]$ is non zero only for the map $\gamma_0$ sending self-loops to 1 and simple edge to 2, for which $\tau^0\big[ T_\gamma(  \mathbb I,   \mathbb J)\big]$ is one. Hence $\tau^0\big[ T(b)\big] = \Phi(a)^{\ell } \Psi\big(a - \Delta(a)\big)^{m }$, with $\ell =\ell_{ \gamma_0}$ and $m=m_{\gamma_0}$, so we get as expected that $a$ is distributed as $b$.

It remains to prove that $b=\Phi(a) \mathbb I + \Psi\big(a - \Delta(a)\big) \mathbb J$ is traffic independent from itself. Let $T$ be a test graph in two variables $x$ and $y$. One the one hand we know that $\tau^0\big[ T(b,b)\big] = 0$ if $T$ is not a simple tree decorated with loops. This is equivalent to say that the graph of colored component of $T$ with respect to $x$ and $y$ is a tree and the colored components are trees decorated with loops. Therefore, we get $\tau^0\big[ T(b,b)\big] = \Phi(a)^{\ell } \Psi\big(a - \Delta(a)\big)^{m }$ with the same notation as above, and for each colored components $S$ of $T$ we have $\tau^0\big[S(b)\big] = \Phi(a)^{\ell_{S}} \Psi\big(a - \Delta(a)\big)^{m_{S}}$ with similar notations. Since $\ell = \sum_S \ell_{S}$ and $m= \sum_S m_{S}$, we get the result.
\\
\par 2. We now prove the law of large numbers. For any test graph $T$ in a single variable $x$ with edge set $E$ and any map $\gamma: E \to [n]:=\{1\etc n\}$, denote by $T_\gamma$ the test graph in $n$ variables $x_1 \etc x_n$ obtained from $T$ by putting the label $x_{\gamma(e)}$ on edge $e$. Then by multi-linearity w.r.t. the edges for $\tau^0$, we get
	$$\tau^0\Big[ T \Big( \frac{a_1 + \dots + a_n}n\Big) \Big] = \sum_{\gamma: E \to [n]} n^{-|E|} \tau^0\big[T_\gamma(a_1 \etc a_n)\big].$$
For any $\gamma$ let $\pi_\gamma$ be the partition of $E$ such that two edges belong to the same block whenever they have same label. Since the $a_i$ are identically distributed and independent, $\tau^0\big[T_\gamma(a_1 \etc a_n)\big]$ depends only on $\pi_\gamma$ and we denote by $\eta(\pi_\gamma)$ this quantity. Denote by $\mathcal P(E)$ the set of partitions of $E$ and by $|\pi|$ the number of blocks of an element $\pi$ of $\mathcal P(E)$. We then get
	\eq
		\tau^0\Big[ T \Big( \frac{a_1 + \dots + a_n}n\Big) \Big] & = & \sum_{\pi \in \mathcal P(E)} n(n-1) \dots (n-|\pi |+1) \times n^{-|E|} \eta(\pi) \\
		&=& \sum_{\pi \in \mathcal P(E)} n^{|\pi|-|E|}\big( 1 + o(1) \big) \eta(\pi).
	\qe
The only partition that contributes is the partitions consisting in singletons. Hence $\tau^0\Big[ T \Big( \frac{a_1 + \dots + a_n}n\Big) \Big] = \tau^0\big[\tilde T(a_1 \etc a_{|E|})\big] +o(1)$ where $\tilde T$ is obtained by putting different labels for its edges. The $a_i$ are independent and identically distributed so we have seen in the proof of the previous point that $ \tau^0\big[\tilde T(a_1 \etc a_{|E|})\big] = \tau^0\big[T\big( \Phi(a)  \mathbb  I + \Psi\big(a - \Delta(a)\big)  \mathbb J\big)\big]$.
\end{proof}

\section{Central limit theorem}

We recall the classical CLTs and then state the ''traffic`` version.

\begin{theorem}[The non commutative central limit theorems] Let consider the three situations (1,2,3), of a non-commutative $^*$-probability space $(\mathcal A, \Phi)$ for (1) and (3) and of a $^*$-algebra $\mathcal A$ endowed with a state $\Phi$ for (2). Consider a sequence $(a_n)_{n\geq 1}$ of identically distributed, self adjoint elements of $\mathcal A$, either free (1), Boolean (2) or tensor (3) independent. Assume that $\Phi(a_n)=0$ and $\Phi(a_n^2)=1$. Then $m_n = \frac{a_1 + \dots + a_n}{\sqrt n}$ converges in distribution to 
\begin{enumerate}
	\item[(1)] a semicircular variable $x$, i.e.  $\Phi(x^k) = \frac 1 {2\pi} \int t^k \sqrt{4-t^2}_+\mathrm dt$,
	\item[(2)] a Rademacher variable $y$ (also called a Bernouilli symmetric variable), i.e. $\Phi(y^k) = \one(k \mathrm{ \ is \ even})$.
	\item[(3)] a Gaussian variable $z$, i.e. $\Phi(z^k) = \frac 1 {\sqrt{2\pi}} \int t^k e^{-\frac {t^2}2}\mathrm dt$,
\end{enumerate}
\end{theorem}

The limit in CLT for traffics will be written as a sum of three terms that represent each of these variables.

\begin{definition}[Central traffics]\label{Def:CentralVar} Let $(\mathcal A, \tau)$ be an algebraic traffic space  with associated trace $\Phi$ and anti-trace $\Psi$ such that $\mathcal A$ is a $\mathcal G^*$-algebra. We say that an element $a\in \mathcal A$ is self-adjoint if $a^*=a$ and off-diagonal whenever $\Delta(a)=0$.
\begin{enumerate}
	\item A centered semicircular traffic $x$ with parameter 
		$$(\alpha,\beta) =\big(\Phi(x^2), \Phi(xx^t)\big), \ |\beta|\leq \alpha,$$
is a self-adjoint off-diagonal traffic with distribution given as follow: for any test graph $T$
	\eqa
		  \tau^0\big[ T(x) \big]  =   \one \big( T \textrm{ is a double tree}\big)
		  								\alpha^{ \ell(T)}	\beta^{k(T)},
	\qea
where $\ell(T)$ (respectively $k(T)$) is the number of double edges with opposing (respectively similar) orientation (see Proposition \ref{Prop:WigDistr}).
	\item A centered simple Boolean traffic $y$ with parameter given by a symmetric non-negative matrix
	$$\big(\alpha_{\ell,k}\big)_{\ell, k\geq 0} = \Big(\Phi\big( deg( y)^\ell deg^t( y)^k \big) \big)_{\ell,k\geq 0}, $$
 is a self-adjoint off-diagonal traffic with distribution given as follow:  for any test graph $T$
	\eqa\label{eq:SimpleBoolDistr}
		  \tau^0\big[ T(y) \big]  =   \one \big( T \textrm{ is a tree}\big)
		  								\sum_{\pi\in \mathcal P(E|V)} \prod_{B_v\in \pi} \alpha_{\ell(v) , k(v)},
	\qea
where $\mathcal P(E|V)$ is the set of partitions of edges of $T$ into blocks $B_v$ of edges having in common a same vertex $v\in V$, and $\ell(v)$ (respectively $k(v)$) is the number of edges of $B_v$ for which $v$ is the source (respectively the target), see Lemma \ref{Lem:BoolGaus}. 

A centered Gaussian Boolean traffic $y$ with parameter
	$$(\alpha, \beta) = \big(\Psi(y^2), \Psi(yy^t)\big),  \ |\beta|\leq \alpha,$$
is a centered simple Boolean traffic with parameter $(\esp[Y^\ell \bar Y^k])_{k,\ell\geq 0}$ for a complex Gaussian random variable $Y$.
 	\item A real centered Gaussian diagonal traffic $z$ with parameter $\alpha=\Phi(z^2)$ is a self-adjoint diagonal traffic whose distribution with respect to $\Phi$ is the centered real Gaussian law with variance $\alpha$.
\end{enumerate}
\end{definition}

A semicircular traffic $x$ is limit of a Wigner matrix by Lemma \ref{Prop:WigDistr}. A centered Gaussian Boolean traffic $y$ is limit of a matrix of the form $\big(\frac{X_i+\overline{ X_j}}N\big)$ as in Lemma \ref{Lem:BoolGaus}. A real centered Gaussian diagonal traffic $z$ is the limit of a diagonal matrix with real Gaussian entries by Lemma \ref{Lem:DiagMatrices}. The non commutative distributions of $x$ and $z$ with respect to $\Phi$ are respectively the semicircular and the Gaussian distributions, the non commutative distribution of $y$ with respect to the $\Psi$ is the Rademacher distribution.

\begin{theorem}[Traffic central limit theorem 1/2]\label{CLT} Let $(\mathcal A,\tau)$ be an algebraic traffic space with associated trace $\Phi$ and anti-trace $\Psi$ which is a $\mathcal G^*$-algebra. Let $(a_n)_{n\geq 1}$ be a sequence of self-adjoint, independent and identically distributed traffics in $\mathcal A$ and let $a$ be distributed as the $a_n$'s. Assume that $\Phi(a) =\Psi(a)=0$. Then $m_n =  \frac{a_1+ \dots + a_n}{\sqrt n}$ converges to a traffic $m$. If $\Phi$ is a state, then $m=x + y +z$ is the sum of a semicircular traffic $x$, a Gaussian Boolean traffic $y$ and a Gaussian diagonal traffic $z$ which in general are not traffic independent.

Nevertheless, $x$ and $z$ are traffic independent. Hence, seen as a non commutative random variable in $(\mathcal A, \Phi)$, $x$ and $z$ are free independent ($y$ has the null distribution w.r.t. $\Phi$). Denoting $\rho = \Phi\big( \Delta(a)  \Delta(a) \big) = \Phi\big( a \circ a \big)$ and assuming $\Phi(a^2)=1$, then the distribution of $m=x+z$ w.r.t. $\Phi$ is the free convolution of the Gaussian law with mean zero and variance $\rho$ and the semicircular law with mean zero and variance $1-\rho$. 
\end{theorem}

\begin{remark} \begin{itemize}
	\item Since we assume that $\Phi$ is a state, $\rho \geq 0$ and so the free convolution is well defined. It is the limit in distribution w.r.t. $\Phi$ of $\sqrt \rho D_N + \sqrt{1-\rho} X_N$, where $D_N$ is a diagonal matrix with i.i.d. Gaussian entries independent from a Wigner matrix $X_N$ (by Corollary \ref{Cor:ApplClass} and Proposition \ref{RigFree}).
We give below in Proposition \ref{CLTMatrix} a matrix model for the sum $m=x+y+z$.
	\item One must be careful with the statement that $x,y$ and $z$ are not independent, since there is no unicity in the decomposition of $m$ as a sum of semicircular, Gaussian Boolean and Gaussian diagonal traffics, as we will see later. The are situations where the variables given by the theorem are not independent but $m$ can be written as a sum of independent variables, see Example \ref{Ex:CLT2}.
\end{itemize}
\end{remark}

\begin{example}\label{Ex:CLT}Let us state an application which is due to Au in \cite{AU16}, from a preliminary version of Theorem \ref{CLT} which contains only the second part of the theorem. We present his result and reasoning which is an interesting use of the traffic CLT in order to find a new proof of a result on Wigner matrices. Let $m$ be a limit in traffic distribution as $N$ goes to infinity of $M_N = p X_N +q\deg(X_N)$, where $p,q\in \mathbb R$ and $X_N$ is a real Wigner matrix. The traffic $m$ is also the limit of $M'_N = p X'_N +q\deg(X'_N)$ for a real Gaussian Wigner matrix $X'_N$. By stability of Gaussian variables, for each $n$ we can write $X'_N$ as a normalized sum $\frac 1 {\sqrt n}  \sum_{i=1}^n {X'_N}^{(i)}$ of i.i.d. real Gaussian Wigner matrices. Hence $m$ has the distribution of a limit as in Theorem \ref{CLT} and so it distribution w.r.t. $\Phi$ is a free convolution of a semicircular distribution and a Gaussian distribution.
 \end{example}

Let $m$ be the limit in traffic distribution of $m_n=\frac{a_1 + \dots + a_n}{\sqrt n}$ in Theorem \ref{CLT}. 
 We still denote $(\mathcal A, \tau)$ the algebraic traffic space where $m$ lives. The good formal setting to describe the decomposition $m=x+y+z$ is given when there exists an element $\mathbb J$ in $\mathcal A$ as in the law of large numbers (Proposition \ref{Prop:LGN}).

\begin{definition}\label{Def:uconit} In an algebraic traffic space $(\mathcal A, \tau)$ with trace $\Phi$ and anti-trace $\Psi$, when it exists we denote by $\mathbb J \in \mathcal A$ a distinguished element independent from all $\mathcal A$ and such that $\Phi(\mathbb J)=0$, $\Psi(\mathbb J)=1$.
\end{definition}

Remark that in such a space $\mathcal A$, one has $a = \Phi(a) \mathbb I +  \big(\Psi(a)-\Phi(a)\big) \mathbb J + b$ where $\Phi(b)=0$ and $\Psi(b)=0$. The independence of $\mathbb J$ with everyone and the formula for its traffic distribution implies the following lemma.

\begin{lemma}\label{Lem:PropJ} In the setting of the previous definition, for any traffics $b_1,b_2$, one has $\Phi( b_1 \mathbb J)=0$ and $\Psi(b_1 \mathbb J b_2) = \Psi(b_1) \Psi(b_2)$.
\end{lemma}

We can now gives an explicit description of the limit in traffic distribution of the central limit theorem.

\begin{proposition}[Traffic central limit theorem 2/2]\label{CLTMatrix} Let $a$, $m_n= \frac{a_1 + \dots + a_n}{\sqrt n}$ and $m$ be as in Theorem \ref{CLT} and denote $\tilde a = a - \Delta(a)$. Assume that in the traffic space generated by $a$, the trace $\Phi$ is a state. We now consider three independent centered variables in some space containing $\mathbb J$ as in Definition \ref{Def:uconit}:
\begin{enumerate}
	\item a semicircular traffic $x'$ with parameter 
		$$(\alpha_1,\beta_1)=\big( \Phi(\tilde a^2), \Phi(\tilde a\tilde a^t) \big),$$
	\item a Gaussian Boolean traffic $y'$ with parameter 
		$$(\alpha_2,\beta_2)=\big( \Psi(\tilde a^2), \Psi(\tilde a\tilde a^t) \big),$$
	\item a Gaussian diagonal traffic $z'$ with parameter  
		$$\alpha_3= \Phi\big[ \Delta(a)^2\big] - \gamma  \Psi\big[ \Delta(a)  \frac{ \tilde a + \tilde a^t}{\sqrt 2} \big],$$
		 where $\gamma = \frac{ \Psi\big[ \Delta(a)  \frac{ \tilde a + \tilde a^t}{\sqrt 2} \big] }{\alpha_2 + \Re(\beta_2)}$ ($\gamma=0$ if $\alpha_2 + \Re(\beta_2)=0$).
\end{enumerate}
We denote 
	$$q(x') = x'\mathbb J + \mathbb J x' \mathrm{ \ and \ } r(y') = \frac{ deg(y') + deg^t(y')}{\sqrt{2}}.$$
 Then $m$ has the same distribution as 
 	$$m'=\big( x'- q(x') \big) + y' + \big( \gamma r(y') + z'\big).$$
 The centered traffics $x=x', y = -q(x')+y'$ and $z=\gamma r(y') + z'$ are respectively semicircular, Gaussian Boolean and Gaussian diagonal traffics.
\end{proposition}

\begin{remark}\label{Rk:TCL2} \begin{itemize}
	\item Since $\Phi$ is a state, the variables are well defined. Indeed, the Cauchy-Schwarz inequality $|\Phi(ab)|^2\leq \Phi(aa^*) \Phi(bb^*)$ implies that $\alpha_1\geq |\beta_1|$ and $x'$ is well defined. As well, since $\Psi(ab) = \Phi[ deg(ab)] = \Phi[deg^t(a) deg(b)]$ for any $a,b$, one has $\alpha_2\geq |\beta_2|$ and $y'$ is well defined. To see that $\alpha_3$ is non-negative, note that $\Psi( \Delta(a) b) =  \Phi[ a \, deg(b)]$ for any $a,b$ and so
	\eq
		\lefteqn{\gamma  \Psi\Big[ \Delta(a)  \frac{ \tilde a + \tilde a^t}{\sqrt 2} \Big] }\\
		& = & \frac 1 {\alpha_2 + \Re(\beta_2)} \Phi\big[ \Delta(a) r(\tilde a) \big]^2 \leq \frac 1 {\alpha_2 + \Re(\beta_2)} \Phi\big[ \Delta(a)^2 \big] \Phi \big[ { r(\tilde a)^2}  \big].
	\qe
 But since $\Psi(ab) = \Phi( deg(a) deg^t(b))$ then 
 	$$\Phi \big[ { r(\tilde a)^2} \big] =\Psi\big[\tilde a \tilde a^*\big] +  \Phi \Big[\frac{ deg(\tilde a)^2 +deg^t(\tilde a)^2}2 \Big].$$
Moreover, by Lemma \ref{Lem:TraceAdj} one has  
	$$\overline {\Phi \big[deg^t(\tilde a)^2\big] } = \Phi \big[(deg^t(\tilde a)^2)^*\big] = \Phi \big[(deg(\tilde a^*)^2)\big],$$
 and since $\tilde a^*=\tilde a$ we get $\Phi \big[ { r(\tilde a)^2}  \big]  = \alpha_2 + \Re(\beta_2)$ and then $\alpha_3\geq 0$.
 	\item Let $X_N'$ be a Wigner matrix with parameter $(\alpha_1,\beta_1)$, $Y_N'$ be a matrix as in Lemma \ref{Lem:BoolGaus} with parameter $(\alpha_2,\beta_2)$, $Z_N'$ be a diagonal matrix with i.i.d. real centered Gaussian entries with variance $\alpha_3$, and $\mathbb J_N$ be the matrix whose all entries are $\frac 1N$, the random matrices being independent. By Theorem \ref{MainTh}, the matrices are asymptotically traffic independent and $m'$ is distributed as the limit of $M_N = X_N' - q_N(X_N')+ Y_N' + \gamma r(Y_N') + Z_N'$, where $q_N(X_N')=X_N'\mathbb J_N + \mathbb J_N X_N'$ and $r(Y_N') = \frac{ deg(Y_N') + deg^t(Y_N')}{\sqrt{2}}$. In particular there exists a traffic space as in the proposition.
	\end{itemize}
\end{remark}

\begin{example}\label{Ex:CLT2} 
\begin{enumerate}
	\item Let $m_n$ be the normalized sum $\frac 1 {\sqrt n}  \sum_{i=1}^n   x_i$ of independent semicircular traffics with parameter $(1, \eta)$. Since $(x_i)_{i\geq 1}$ is the limit of a family of independent Gaussian Wigner matrices, $m_n$ and its limit $m$ as $n$ goes to infinity are also semicircular traffics with same parameter. Yet, the answer given by the above proposition looks more complicated. Using Definition \ref{Def:CentralVar}, we compute the parameters in order to apply Proposition \ref{CLTMatrix} is as follow. Let $x$ be distributed as the $x_i$.
\begin{itemize}
	\item $\Phi(x) =0$ by extra-diagonality of $x$ (that is $\Delta(x)=0$) by noting that $\Phi(x) = \Phi\big(\Delta(x)\big)$.
	\item $\Psi(x) = \tau\big[ \cdot \overset{x} \rightarrow \cdot  ] =   
		\tau^0\big[ \cdot \overset{x} \rightarrow \cdot  ] +
		 \tau^0\big[ \circlearrowleft^x ] =0$.
	\item $\Phi(x^2) = 1$ and $\Phi(xx^t) =\eta$ by definition of the parameters.
	\item $\Psi(x^2) =  \tau\big[ \cdot  \overset{x}\leftarrow \cdot \overset{x}\leftarrow \cdot\big]$. The only quotient of the latter graph that is a double tree (a graph which becomes a tree when the multiplicity of edges is forgotten) is $(\cdot  \underset{x}{\overset x{\leftrightarrows}}\cdot)$, and so $\Psi(x^2) = \tau^0 \big[(\cdot  \underset{x}{\overset x{\leftrightarrows}}\cdot)\big] = 1$.
	\item $\Psi(xx^t)=  \tau\big[ \cdot  \overset{x}\leftarrow \cdot \overset{x^t}\leftarrow \cdot\big]= \tau\big[ \cdot  \overset{x}\leftarrow \cdot \overset{x}\rightarrow \cdot\big]$. Similarly we have $\Psi(xx^t) = \tau^0 \big[(\cdot  \underset{x}{\overset x{\leftleftarrows}}\cdot)\big] = \eta$.
	\item $\Phi(\Delta(x)^2) = \Phi(\Delta(x)r(x))=0$ by extra-diagonality of $x$.
\end{itemize}

	 Hence $m$ has the distribution of $\big(x' -q(x')\big) + y'$, where $x'$ and $y'$ are independent, semicircular and Gaussian Boolean respectively, with same parameter as the $x_i$. 
	
	\begin{corollary}\label{Cor:SCplusGB} Let $x', y'$ be independent traffics, semicircular and Gaussian Boolean respectively, with same parameter. Then $x=x' - q(x')+y'$ is a semicircular traffic with same parameter.
	\end{corollary}
	
	Although $m$ as the same distribution as $x'$, the proposition invites us to first remover a part of $x'$ with the term $-q(x')$, and then to sample it again by adding $y'$. The reason of this approach will be clear with the third example below. Note that when the $x_i$ are nonzero, then the semicircular part $x'$ and the Gaussian Boolean part $-q(x') + y'$ are not independent. Indeed by Lemmas \ref{Lem:Elementary} and \ref{Lem:PropJ} one can see that $\Psi\big( x'(-q(x') + y')\big) = -\Psi\big( x'q(x') \big) = -\Phi\big((x')^2\big)<0$, which should vanish if they were independent. 	
	\item Let $m_n$ be the limit as $N$ goes to infinity of the normalized sum 
		$$\frac 1 {\sqrt {2n}} \sum_{i=1}^n (U_N^{(i)} + U_N^{(i)*})$$
 of independent Haar unitary matrices and their adjoint. By Theorem \ref{MainTh} and Proposition \ref{Prop:HaarDist}, $m_n$ is the normalized sum of i.i.d. self-adjoint centered traffics. Then the central limit theorem implies that the limit $m$ of $m_n$ is a semicircular traffic with parameter $(1,0)$, thanks to the same detour as in the previous example.
	
	Indeed, let us compute the parameters we need in the central limit theorem, using the limiting distribution of a Haar unitary matrix $(U_N, U_N^*)$ given in Proposition \ref{Prop:HaarDist}. We denote the traffics $(u,u^*) = \Nlim \big(U_N , U_N^{ *}\big)$ and $a=u+u^*$. Note that $a$ has the distribution of an extra-diagonal traffic since for a test graph $T$ labeled in $(u,u^*)$, $\tau^0[T] =0$ if it has a single loop. In particular $a$ has the same distribution as $\tilde a = a -\Delta(a)$.
	
	\begin{itemize}
		\item $\Phi(a) =0$ by extra-diagonality.
		\item $\Psi(a) =2 \Re \, \Psi(u) = 2 \Re \, \tau\big[ \cdot \overset{u} \rightarrow \cdot  ] =  2 \Re \, \big( 
		\tau^0\big[ \cdot \overset{u} \rightarrow \cdot  ] +
		 \tau^0\big[ \circlearrowleft^u ] \big)=0$ by Lemma \ref{Lem:TraceAdj}.
		 \item $\Phi( a^2) =2 \Re \Phi(u^2) + 2.$ But  $\Phi(u^2) = \tau\big[ \cdot  \underset{u}{\overset u{\leftrightarrows}}\cdot\big] =  0$ 	since $\tau^0[T]\neq 0$ only if $T$ has the same number of edges labeled $u$ and $u^*$ and that the quotient of a graph has the same number of edges labeled $u$ and $u^*$.
		 \item $\Phi( aa^t) = \Phi( uu^t + uu^{*t} + u^*u^t + u^*u^{*t}) = 0$ since 
	$$ \Phi(uu^t) = \tau \big[\cdot  \underset{u^t}{\overset u{\leftrightarrows}}\cdot \big] = \tau \big[\cdot  \underset{u}{\overset u{\leftleftarrows}}\cdot \big] = \tau^0 \big[\cdot  \underset{u}{\overset u{\leftleftarrows}}\cdot \big]  + \tau^0\big[   \,^u\circlearrowleftRotBis\cdot \circlearrowleftRot^u \big]=0$$
	and  $ \Phi(u^*u^{*t})=0$ with the same computation, and the two other terms vanishes as the number of edges labeled $u$ and $u^*$ are not equal.
	\item $\Psi( a^2) = 2 \Re \Psi(u^2) + 2 = 2$ for the same reason as before.
	\item $\Psi( aa^t) = \Psi( uu^t + uu^{*t} + u^*u^t + u^*u^{*t})=0$ since  as before 
		$$\Psi( uu^t) = \Psi(u^*u^{*t})=0$$ and 
	\eq
	\Psi(uu^{*t}) & = & \tau\big[ \cdot  \overset{u}\leftarrow \cdot \overset{(u^*)^t}\leftarrow \cdot\big] =  \tau\big[ \cdot  \overset{u}\leftarrow \cdot \overset{u^*}\rightarrow \cdot\big] =  \tau^0\big[ \cdot  \overset{u}\leftarrow \cdot  \overset{u^*}\rightarrow  \cdot\big] +\tau^0\big[  \,^u\circlearrowleftRotBis  \cdot  \overset{u^*}\rightarrow  \cdot  \big] 
	\\ & &+ \tau^0\big[   \cdot \overset{u}\leftarrow \cdot  \circlearrowleftRot^{u^*}\big]  +\tau^0\big[\cdot  \underset{u}{\overset {u^*}{\leftleftarrows}}\cdot \big]+ \tau^0\big[   \,^u\circlearrowleftRotBis \cdot  \circlearrowleftRot^{u^*} \big]   = 0
\qe
and the same computation yields $\Psi(u^*u^t )=0$
	\item $\Phi(\Delta(a)^2) = \Phi(\Delta(a)r(a))=0$ since $a$ has the distribution of an extra-diagonal traffic. 
	\end{itemize}
	Hence by Proposition \ref{CLTMatrix} and Corollary \ref{Cor:SCplusGB}, $m$ is a semicircular traffic with parameter $(1,0)$.
	\item In the previous examples, one can wonder why we do not chose the Gaussian Boolean part $y'$ in such a way $\Psi({y'}^2) = (\Psi - \Phi)(\tilde a^2)$, in order to not remove a part that is sampled again. The reason is that the quantity in the r.h.s. term is possible negative, as in the following example.
	
		Let $m_n$ be the limit as $N$ goes to infinity of the standardized sum $\frac 1 {\sqrt {2n}} \sum_{i=1}^n (V_N^{(i)} + V_N^{(i)^t} - 2\mathbb J_N)$ of independent uniform permutation matrices $V_N^{(i)}$ and their transpose, where $\mathbb J_N$ denotes the matrix whose all entries are $\frac 1 N$. By Theorem \ref{MainTh}, Proposition \ref{Prop:PermDistr} and the convergence of $\mathbb J_N$, $m_n$ is the normalized sum of i.i.d. traffics. Then the central limit theorem implies that the limit $m$ of $m_n$ has the distribution of $x'-q(x')$ where $x'$ is a semicircular traffic with parameter $(1,1)$. 
	
	Indeed, let us denote $v = \Nlim V_N^{(i)}$ and $a=v+v^t-2\mathbb J$. As in the previous case, $a$ has the distribution of an extra-diagonal traffic. 

\begin{itemize} 
		\item $\Phi(a) =0$  by extra-diagonality.
		\item $ \Psi(v) = \tau\big[ \cdot \overset{v} \rightarrow \cdot  ] =    
		\tau^0\big[ \cdot \overset{v} \rightarrow \cdot  ] +
		 \tau^0\big[ \circlearrowleft^v ] =1+0$ and so $\Psi(a)=0$.
		 \item By Lemma \ref{Lem:PropJ}, $\Psi(v^\varepsilon\mathbb J) = \Psi(\mathbb J v^\varepsilon)=0$ for $\varepsilon\in \{1,t\}$, so one has $\Phi( a^2) =  \Phi\big( (v+v^t)^2\big) = 2+2 \Phi(v^2)=2$ since 
	$$\Phi(v^2) = \tau\big[ \cdot  \underset{v}{\overset v{\leftrightarrows}}\cdot\big] = \tau^0\big[ \cdot  \underset{v}{\overset v{\leftrightarrows}}\cdot\big] + \tau^0\big[   \,^v\circlearrowleftRotBis \cdot  \circlearrowleftRot^v \big]=0 $$
		 \item Since $a^t=a$ we have $\Phi( aa^t) = \Phi( a^2)=2$.
	\item We can anticipate that $\Psi( a^2) =0$. Indeed, the matrix $V_N^{(i)} + V_N^{(i)t}$ is the adjacency matrix of a random graph (the graph of cycles of the associated permutation). The degree of each vertex of the graph (the number of neighbors) is two, and $\Psi\big( (V_N^{(i)} + V_N^{(i)^t} - 2\mathbb J_N)\big)$ is nothing else than the variance of the degree of a vertices uniformly chosen at random. More generally, we have the following fact.
	
	\begin{lemma} Let $(\mathcal A, \tau)$ be an algebraic traffic space with anti-trace $\Psi$, such that $\mathbb J \in \mathcal A$ as in Definition \ref{Def:uconit} exists. Let $a\in \mathcal A$ having the same distribution as $b - (b \mathbb J + \mathbb J b)$ for some $b\in \mathcal A$ such that $\Psi(b)=0$. Then one has $\Psi(a^2)=0$.
	\end{lemma}
	
	\begin{proof} Recall that $\Psi(b_1\mathbb J b_2) = \Psi(b_1) \Psi(b_2)$ for any $b_1,b_2\in \mathcal A$ by Lemma \ref{Lem:PropJ}. Hence we get
	\eq
		\Psi(a^2) & = & \Psi(b^2) - \Psi\big( b(b\mathbb J + \mathbb J b) + (b\mathbb J + \mathbb J b)b \big) + \Psi\big( (b\mathbb J + \mathbb J b)^2 \big)\\
	& = & \Psi(b^2) - 2 \Psi(b^2) + \Psi(b^2) = 0
	\qe
\end{proof}
	
	For the matrix $A_N = V_N^{(i)} + V_N^{(i)t} - 2\mathbb J_N$, one has $\Psi(A_N)=0$ and $(A_N \mathbb J_N + \mathbb J_N A_N)=0$. Hence the limit $a$ of $A_N$ satisfies the assumption of the lemma. 
	\item $\Psi( aa^t) = \Psi( a^2) = 0$.
	\item $\Phi(\Delta(a)^2) = \Phi(\Delta(a)r(a))=0$ since $a$ has the distribution of an extra-diagonal traffic. 
\end{itemize}
	Hence the expected result: the limit $m$ in the central limit theorem has the distribution of $x'-q(x')$ for a semicircular traffic with parameter $(1,1)$. The traffics $x'$ and $-q(x')$ are not independent since $\Psi\big( x' q(x') \big) = \Psi({x'}^2)=1$, and there is no known way to decompose $x'-q(x')$ as a sum of independent semicircular and Gaussian Boolean variables. 
\end{enumerate}
\end{example}

We first prove in Theorem \ref{CLT} the convergence in traffic distribution of $(m_n)_{n\geq 1}$, giving an expression for the distribution of its limit $m$ in formula \eqref{Eq:ProofCLT} below. Then we compute the distribution w.r.t. $\Phi$ to prove that it is a convolution of a Gaussian and of a semicircular distribution. The rest of the proof is dedicated to Proposition \ref{CLTMatrix}

\begin{proof}[Proof of the convergence in Theorem \ref{CLT}]
 Let $T=(V,E)$ be a test graph in one variable. With the same computation as for the law of large numbers, the multilinearity of $\tau^0$ with respect to the edges gives
 	$$\tau^0\Big[T\big(\frac{a_1+ \dots + a_n}{\sqrt n}\big)\Big] = \sum_{\pi \in \mathcal P(E)} n^{|\pi|-\frac{|E|}2}\big( 1+o(1)\big) \eta(\pi).$$
Here $\eta(\pi)$ equals $\tau^0\big[T_\gamma(a_1 \etc a_{|\pi|})\big]$, where $\gamma:E\to \{1\etc {|\pi|}\}$ is such that $e\sim_\pi f$ if and only if $\gamma(e) = \gamma(f)$, and $T_\gamma$ is obtained from $T$ by putting for each edge $e$ the label corresponding to $a_{\gamma(e)}$.
Assume that $\pi$ has a block of size one. Then, either $\mathcal G\mathcal C\mathcal C(T_\gamma)$ is not a tree and $\eta(\pi)=0$, or by the factorization property of traffic independence, one can factorize in $\eta(\pi)$ the term $\tau^0\big[ T(a)\big]$ where $T$ has a single edge. If $T=\circlearrowleft$ is a simple loop, then $\tau^0\big[\circlearrowleft(a)\big] = \tau\big[\circlearrowleft(a)\big] = \Phi(a) = 0$. 
If $T$ is a simple edge, since $ \tau^0 \big[T(a)\big]+\tau^0\big[\circlearrowleft(a)\big] =\tau\big[T(a)\big] = \Phi\big( deg(a)\big)=0$ we get as well $\tau^0\big[T(a)\big] =0$.
Hence the only partitions that contribute are those for which $\pi$ is a pair partition. This proves the convergence of $m_N$ in traffic distribution.  
 
Moreover, by definition of traffic independence, the partitions $\pi$ that contribute are those for which the graph of colored components of $T_\gamma$ is a tree. Since $\Phi(a) = \Phi\big( deg(a)\big) = 0$ and by the factorization property, the partitions must pair adjacent edges (that share at least one vertex). 

Hence, the only test graphs $T$ for which $\tau^0\big[T(m_n)\big]$ does not vanish at infinity are graphs that become trees if we forget the multiplicity of the edges and delete the loops, such that the multiplicity of simple edges (that are not self loops) is one or two. Denote by $\mathcal T_0$ the set of such graphs. Denote by $\mathcal P_{0}(T)$ the set of pair partitions $\pi$ of adjacent edges of $T$ such that $\{e_1,e_2\} \in \pi$ for any twin edges $e_1,e_2$ of $T$. For $\pi \in \mathcal P_{0}(T)$, let us identify $S\in \pi$ with the subgraph of $T$ consisting in the edges of $S$. We then get
	\eqa\label{Eq:ProofCLT}
		\tau^0\big[ T(m_n)\big] \underset{n\rightarrow \infty} \longrightarrow  \one( T \in \mathcal T_0) \sum_{\pi \in \mathcal P_{0}(T)} \prod_{S\in \pi} \tau^0[S].
	\qea
\end{proof}
\begin{proof}[Limiting distribution w.r.t. $\Phi$]
Let now compute the distribution w.r.t. $\Phi$ of $m$. Let $T$ consisting in a cycle of length $k$. Denoting by $V$ its vertex set, we have $\Phi( m^k) = \tau\big[ T(m)\big] = \sum_{\sigma\in \mathcal P(V)} \tau^0\big[T^\sigma(m)\big]$. For any partition $\sigma$ of $V$, the graph $T^\sigma$ has no cutting edge (edge whose removal disconnect the graph). Hence, a graph $T^\sigma$ in $\mathcal T_0$ for some $\sigma$ in $\mathcal P(V)$ consists in a double tree $T_0$ for which at each vertex $v\in V$ is attached an ensemble of self-loops $F_v$. Computing $\tau^0\big[   T^\sigma(m)\big]$ with the above formula, the pair partitions $\pi\in \mathcal P_{0}(  T^\sigma)$ in the sum must gather twin edges of $T_0$ and pair of loops attached at a same vertex. Hence $|F_v|$ must be even for this term to not vanish. In this case, denote by $2\ell_v$ the number of loops attached to a vertex $v\in V$. By Lemma \ref{lem:EdgesVertices} that gives the relation between the number of vertices and edges in a tree, the number of edges of $T_0$ is $2(|V|-1)$.

 We get, computing for the graph $S$ consisting in a double edge $\tau^0[S(a)]=\Phi( a^2)-\Phi(\Delta(a)^2) = 1 -\rho$,
\begin{eqnarray*}
	\tau^0\big[ \tilde T(m) \big]  & = &   (1-\rho)^{|V|-1} \prod_{v\in V} \rho^{\ell_v} \textrm{Card } \mathcal P_2(2m_v),
\end{eqnarray*}

\noindent where $ \mathcal P_2(2m)$ denotes the set of pair partitions of $2m$ elements. But 
		$$\textrm{Card }  \mathcal P_2(2m) = (2m-1)\times (2m-3) \dots 5 \times 3 \times  1 = \esp[X^{2m}]$$
where $X$ is a random variable distributed according to the gaussian measure centered of unit variance (by a basic enumeration and by integration by part respectively).
\\
\\Now, let $x$ be a centered semicircular traffic with parameter $(1,0)$, traffic independent from $z$ a centered Gaussian diagonal traffic with parameter $(1)$. Let prove that $m$ has the same distribution w.r.t. $\Phi$ as $\bar m = \sqrt{ \rho} z + \sqrt{1-\rho} x$. Let $T$ consisting in a cycle of length $k$ and for any partition $\sigma$ of its vertex set,
\begin{eqnarray*}
	\tau^0\big[ T^\sigma(\bar m) \big] & = & \sum_{\gamma: E \to \{1,2\}} \tau^0\big[ T^\sigma_\gamma(\sqrt \rho z,\sqrt{1-\rho} x) \big],
\qe
where in $T^\sigma_\gamma$ an edge $e$ has label $\gamma(e)$. By the definition of traffic independence, the support of the injective distribution of $(z,x)$ consists in double trees $T_0$ with a bunch of self loops $F_v$ attached at each vertex $v$. If $T^\sigma$ is such a test graph, the only map $\gamma$ which makes $ \tau^0\big[ \tilde T(\sqrt \rho z,\sqrt{1-\rho} x) \big]$ possibly non zero consists in labeling the edge of $T_0$ with labels 1 (for $x$) and the self-loops by label 2 (for $z$). By multi-linearity for $\tau$ we get
\eq
	\tau^0\big[ T(\sqrt \rho z + \sqrt{1-\rho} x) \big] & = &  (1-\rho)^{|V|-1} \prod_{v\in V} \rho^{m_k} \esp[X^{2m_v}],
\qe
as expected. By the first case in Theorem \ref{Th:ThreeIndepTraf}, since $x$ is unitarily invariant and $x$ and $z$ are traffic independent, they are free independent.
\end{proof}
 
\begin{proof}[Proof of Proposition \ref{CLTMatrix}] In order to prove that $m$ and $m'$ are equal in distribution, we can compute directly the distribution of $m'=\big( x'- q(x') \big) + y' + \big( \gamma r(y') + z'\big)$ using the distributions of $x',y',z'$ and their traffic independence. But it is much simpler to use Au's argument \cite{AU16}. The variable $m'$ is the limit of the matrix model $M_N = X_N' - q(X_N')+ Y_N' + \gamma r(Y_N') + Z_N'$ as in the above remark where the matrices $X_N',Y_N',Z_N'$ have Gaussian entries. Since $M_N$ is linear in these matrices, it is also Gaussian and so it is stable: it can be written $M_N=M_{N,n} = \frac 1 {\sqrt n} \sum_{i=1}^n M^{(n)}_{N}$ where the $M^{(n)}$ are independent copies of $M_N$. By Theorem \ref{MainTh} of asymptotic traffic independence, the limit $m'$ of $M_N$ is distributed as a limit of the central limit theorem (Theorem \ref{CLT}).

To prove that the variables $m = \underset{n\rightarrow \infty}\lim \frac{a_1 + \dots + a_n}{\sqrt n}$ and $m'=\big( x'- q(x') \big) + y' + \big( \gamma r(y') + z'\big)$ of Proposition \ref{CLTMatrix} have the same distribution, by Formula \eqref{Eq:ProofCLT} it remains to show that for any test graph $T$ with two edges, $\tau^0\big[ T(m')\big] = \tau^0\big[ T(m)\big]$. Note that from \eqref{Eq:ProofCLT} we have $\tau^0\big[ T(m)\big] = \tau^0\big[ T(a)\big]$ for test graphs with two edges, where we recall that $a$ is distributed as $a_1 \etc a_n$. Moreover, it is sufficient to prove the equality $\tau\big[ T(m')\big] = \tau\big[ T(a)\big]$ for these graphs (the number of edges of a test graph is unchanged when identifying vertices). There is a total of eight connected directed graphs with two edges, but since the variables are self-adjoint two pairs of quantities are related each other: by Lemma \ref{Lem:TraceAdj}, we have 
	$$\tau[ \cdot \rightarrow \cdot  \leftarrow \cdot ]= \tau \big[ (\cdot \leftarrow \cdot\rightarrow   \cdot)^* \big]= \overline{\tau \big[ \cdot \leftarrow \cdot\rightarrow \cdot   \big]},$$ 
 and
	$$\tau[ \cdot \leftarrow \cdot\circlearrowleft ] = \tau\big[ (\cdot \rightarrow\cdot \circlearrowleft )^* \big ] = \overline{ \tau[ (\cdot \rightarrow\cdot \circlearrowleft   ] }.$$

Equivalently, it is then sufficient to prove that the six quantities
\eqa
	\Phi\big( \tilde m'  (\tilde m')^\varepsilon\big), \Psi\big( \tilde m'  (\tilde m')^\varepsilon\big), \Psi\Big( \Delta(m') \frac{\tilde m' + (\tilde m')^t}{\sqrt 2}\Big), \Phi\big( \Delta(m' )^2\big),
\qea
where $ \tilde m' := m' -\Delta(m')$ and $\varepsilon \in \{1,t\}$, are equal to the same quantities where $m'$ is replaced by $a$.

Note that by Lemma \ref{Lem:PropJ}, for any traffics $b_1, b_2$ one has $\Phi(b_1 q(b_2))=0$ and $\Psi( b_1 q(b_2)) = \Psi(b_1b_2) + \Psi(b_1)\Psi(b_2)$. Recall form Lemma \ref{Lem:Elementary} that if $b_1$ and $b_2$ are independent then $\Phi(b_1b_2)=\Phi(b_1) \Phi(b_2)$ and $\Psi(b_1b_2) = \Psi(b_1)\Psi(b_2)$. Hence, using the fact that $\Phi\big( y' (y')^t \big)=0$ since it is of Boolean type and that the variables are centered, we get for $\varepsilon \in \{1,t\}$
	\eq
		\Phi( \tilde m'  (\tilde m')^\varepsilon) & = & \Phi\Big( \big(x' - q(x') \big) \big(x' - q(x') \big)^\varepsilon \Big) + \Phi\big( y' (y')^t \big) \\
		& = & \Phi\big( x' (x')^\varepsilon\big) = \Phi(\tilde a \tilde a^\varepsilon),\\
		\Psi( \tilde m'  (\tilde m')^\varepsilon) & = & \Psi\Big( \big(x' - q(x') \big) \big(x' - q(x') \big)^\varepsilon \Big) + \Psi\big( y' (y')^t \big)\\
		&  =& \Psi\big( y' (y')^t \big) = \Psi(\tilde a \tilde a^\varepsilon).
	\qe
Moreover, we have
	\eq
		\Psi\big( \Delta(m')\frac{\tilde m' + (\tilde m')^t}{\sqrt 2}\big) &  = & \Psi\big( \gamma r(y')   \frac{y' + (y')^t}{\sqrt 2} \big) =  \Psi\big(   r(a)   \frac{\tilde a + \tilde a^t}{\sqrt 2} \big) \frac{ \eta}{\alpha_2 + \Re(\beta_2)},
	\qe
where, using the definition of $deg$ in terms of graph operations and Lemma \ref{Lem:TraceAdj}, we have
	\eq
		\eta & = & \frac 1 2   \Psi\big( deg(y')y' + deg(y')(y')^t + deg^t(y')y' + deg^t(y')(y')^t \big) \\
		& = & \frac 1 2   \Psi\big( (y')^t y' + (y')^t (y')^t + y' y' + y' (y')^t\big) = \alpha_2 + \Re(\beta_2).
	\qe
	
At last, we have 
	\eq
		\Phi\big( \Delta(m')^2 \big) & = &  \gamma^2 \Phi\big( r(y')^2\big) + \Phi\big( (z')^2\big)\\
		& = & \gamma \frac{ \Psi\big( \Delta(a) \frac{ \tilde a + \tilde a^t}{\sqrt 2} \big) }{\big( \Psi(\tilde a^2) + \Re \Psi (\tilde a \tilde a^t) \big)} \big( \Psi(\tilde a^2) + \Re \Psi(\tilde a \tilde a^t) \big)\\
		& & \ \  + \Phi\big( \Delta(a)^2 \big) - \gamma  \Psi\big( \Delta(a) \frac{ \tilde a + \tilde a^t}{\sqrt 2} \big) \\
		& = & \Phi\big( \Delta(a)^2 \big).
	\qe
This finishes the proof that $m$ is distributed as $m'$.

It remains to proof the last statement of the proposition, telling that $-q(x') + y'$ is a Gaussian Boolean traffic and that $\gamma r(y') + z'$ is a Gaussian diagonal traffic. By the model of large random matrices, the sum of two independent Gaussian Boolean (respectively Gaussian diagonal) traffics are Gaussian Boolean (respectively Gaussian diagonal). So the following lemma concludes the proof.
\end{proof}

\begin{lemma}\label{Lem:CLT} 
\begin{enumerate}
	\item If $x$ is a centered semicircular traffic with parameter $(\alpha,\beta)$ then $q(x) = x \mathbb J + \mathbb J x$ is Gaussian Boolean with parameter $(\alpha, \beta)$.
	\item If $y$ is a centered Gaussian Boolean with parameter $(\alpha, \beta)$, then  
		$$r(y) = \frac{ deg(y) + deg^t(y)}{\sqrt 2}$$
	 is centered Gaussian diagonal random variable with parameter $(\alpha + \Re (\beta))$.
\end{enumerate}
\end{lemma}

\begin{proof} {\bf 1.} Let $X_N = \left( \frac{x_{i,j}}{\sqrt N}\right)_{i,j}$ be a Wigner matrix with centered Gaussian entries such that $\esp[ |x_{i,j}|^2]=\alpha$ and $\esp[x_{i,j}^2]=\beta$. Then $q_N(X_N):= X_N\mathbb J_N + \mathbb J_N X_N =  \left (\frac{Y_i + \bar Y_j}N\right)_{i,j}$ where $Y_i = \sum_j \frac{x_{i,j}}{\sqrt N}$ is a Gaussian variable such that $\esp[ |Y_i|^2]=\alpha$ and $\esp[Y_i^2]=\beta$. Then $X_N$ converges to the semicircular variable $x$ and $q_N(X_N)$ converges to a centered Gaussian Boolean traffic $y$. Hence $q(x) = y$ in distribution.

{\bf 2.} Firstly, $r(y)$ is diagonal since $deg(y)$ and $deg^t(y)$ are diagonal. Moreover $deg(y)^* = deg^t(y)$ since $y$ is self-adjoint and so $r(y)$ is self-adjoint. It remains to prove that $\Phi\big(r(y)^K\big) = \esp[ Z^K]$ for a real centered Gaussian variable $Z$ such that $\esp[Z^2]=(\alpha + \Re (\beta))$, for any $K\geq 1$. But, with $Y$ denoting a centered Gaussian random variable such that $\esp[|Y|^2] = \alpha$, $\esp\big[ Y^2\big] = \beta$, one has
	\eq
		\Phi\big( r(y)^K\big) & = & \Phi\left( \frac{ deg(y) + deg^t(y)}{\sqrt 2} \right) =2^{-\frac k 2} \sum_{k=0}^K \binom{K}{k}  \Phi\big( deg(y)^k deg^t(y)^{K-k}\big) \\
		& = &  \Phi\left( \frac{ deg(y) + deg^t(y)}{\sqrt 2} \right) =2^{-\frac k 2} \sum_{k=0}^K \binom{K}{k} \esp[Y^k \bar Y^{K-k}] \\
		& = & \esp\Big[ \big(\frac {Y + \bar Y}{\sqrt 2}\big)^K \Big] = \esp[Z^K].
	\qe
\end{proof}

\begin{center}
	{\bf Acknowledgments:}
\end{center}
This work has been done over several years and has received the support of many people, in particular Mikael de la Salle, Alice Guionnet, Roland Speicher, James Mingo, Guillaume C\'ebron, Antoine Dahlqvist, Benson Au, Florent Benaych-Georges, Sandrine P\'ech\'e, Djalil Chafa\"i, Muriel Livernet and Gregory Ginot. The author would like to gratefully thank them for their suggestions and advices during the preparation of this paper. The author thanks Antoine Dahlqvist and Guillaume C\'ebron for their help in correcting an error in Proposition \ref{Prop:HaarDist}, Benson Au for his careful proofreading, and the reviewers for their suggestions.

\backmatter
\bibliographystyle{amsplain}
\providecommand{\bysame}{\leavevmode\hbox to3em{\hrulefill}\thinspace}
\providecommand{\MR}{\relax\ifhmode\unskip\space\fi MR }
\providecommand{\MRhref}[2]{%
  \href{http://www.ams.org/mathscinet-getitem?mr=#1}{#2}
}
\providecommand{\href}[2]{#2}

\end{document}